\documentclass[times,sort&compress,3p,authoryear]{elsarticle}
\usepackage[labelfont=bf]{caption}

\usepackage{amsmath,amsfonts,amssymb,amsthm,booktabs,color,epsfig,graphicx,hyperref,url}
\usepackage{natbib}
\usepackage{subcaption}
\usepackage{caption}
\usepackage{float} 
\usepackage{multirow}
\usepackage{diagbox}
\usepackage{etoc}
\usepackage{minitoc}
\usepackage{titletoc}
\usepackage{appendix}
\theoremstyle{plain}
\newtheorem{theorem}{Theorem}

\newtheorem{proposition}{Proposition}
\newtheorem{lemma}{Lemma}

\theoremstyle{definition}

\newtheorem{Remark}{\bf Remark}

\newtheorem{Assumpa}{Assumption}
\newtheorem{Assumpb}{Assumption}
\newtheorem{Assumpc}{Assumption}
\newtheorem{Assumpd}{Assumption}
\newtheorem{Assumpaa}{Assumption}
\newtheorem{Assumpbb}{Assumption}

\newcommand{\bA}{\boldsymbol{A}}
\newcommand{\bB}{\boldsymbol{B}}

\newcommand{\bD}{\boldsymbol{D}}

\newcommand{\bG}{\boldsymbol{G}}
\newcommand{\bH}{\boldsymbol{H}}
\newcommand{\bI}{\boldsymbol{I}}

\newcommand{\bP}{\boldsymbol{P}}
\newcommand{\bQ}{\boldsymbol{Q}}
\newcommand{\bR}{\boldsymbol{R}}
\newcommand{\bS}{\boldsymbol{S}}
\newcommand{\bT}{\boldsymbol{T}}
\newcommand{\bU}{\boldsymbol{U}}
\newcommand{\bV}{\boldsymbol{V}}
\newcommand{\bW}{\boldsymbol{W}}
\newcommand{\bX}{\boldsymbol{X}}
\newcommand{\bY}{\boldsymbol{Y}}
\newcommand{\bZ}{\boldsymbol{Z}}

\newcommand{\bg}{\boldsymbol{g}}

\newcommand{\bx}{\boldsymbol{x}}

\newcommand{\Prob}{\mathbb{P}}
\newcommand{\Expe}{\mathbb{E}}
\newcommand{\Var}{\mathrm{Var}}

\newcommand{\tr}{\mathrm{tr}}

\newcommand{\tsum}{\mathrm{sum}}
\newcommand{\com}{\mathrm{com}}

\renewcommand{\bar}{\overline}

\usepackage{dsfont}

\makeatletter
\newcommand*{\transpose}{%
	{\mathpalette\@transpose{}}%
}
\newcommand*{\@transpose}[2]{%
	\raisebox{\depth}{$\m@th#1\intercal$}%
}
\begin{document}

\begin{frontmatter}

\title{On testing mean of high dimensional compositional data\tnoteref{t1}}

\author[a]{Qianqian Jiang}

\author[a]{Wenbo Li}

\author[a]{Zeng Li\corref{cor1}}
\ead{liz9@sustech.edu.cn}

\cortext[cor1]{Corresponding author}

\affiliation[a]{organization={Department of Statistics and Data Science, Southern University of Science and Technology},
            city={Shenzhen},
            country={China}}

\begin{abstract}
We investigate one/two-sample mean tests for high-dimensional compositional data when the number of variables is comparable with the sample size, as commonly encountered in microbiome research. Existing methods mainly focus on max-type test statistics which are suitable for detecting  sparse signals. However, in this paper, we introduce a novel approach using sum-type test statistics which are capable of detecting weak but dense signals. By establishing the asymptotic independence between the max-type and sum-type test statistics, we further propose a combined max-sum type test to cover both cases. We derived the asymptotic null distributions and power functions for these test statistics. Simulation studies demonstrate the superiority of our max-sum type test statistics which exhibit robust performance regardless of data sparsity.
\end{abstract}

\begin{keyword} 
High-dimensional compositional data \sep
One/Two-sample mean test \sep
Max-Sum type \sep 
Sum-type \sep 
Max-type.
\MSC[2020] Primary 62H15 \sep
Secondary 62F05
\end{keyword}

\end{frontmatter}

\section{Introduction}
Compositional data, which refers to observations whose sum is a constant, such as proportions or percentages, play a pivotal role in various scientific fields such as geology, economics, and genomics. In microbiome studies, for instance, the varying amounts of DNA-generating material among different samples often lead to the normalization of sequencing read counts to relative abundances, resulting in compositional data \cite{Li2015MicrobiomeMA}. In the context of microbiome and metagenomic research, a fundamental challenge in microbiome data analysis is to test whether two populations share the same microbiome composition. This issue can be conceptualized as a two-sample testing scenario for compositional data. Existing methods for two-sample mean test on compositional data are typically designed for low dimensions, where the number of variables is less than  the sample size. For instance, \cite{Jones2003TheSA} discussed the generalized likelihood ratio test in their book. However, compositional data vectors are inherently high-dimensional due to the diverse array of microbe types in the microbiome. Classical testing methods may exhibit limited power or may not be applicable when the number of variables approaches or exceeds the sample size. Consequently, new statistical methods are required to effectively analyze high-dimensional compositional data.

We specifically focus on mean tests tailored for high-dimensional compositional data.  In current literature, \cite{Cao2018TwosampleTO} first introduced a max-type two-sample mean test based on the maximum of normalized mean difference,  demonstrating significant testing power for sparse alternative hypotheses. However, this test statistic exhibits limited power in identifying  dense yet weak signals where the differences between the means are individually weak and span across multiple locations. To address this issue, we first propose a sum-type Hotelling's statistic to test weak but dense signals in high-dimensional compositional data. In fact, sum-type test statistics for dense alternative hypotheses have been extensively studied in many works including \cite{Hu2016ARO}, \cite{Chen2011ARH}, \cite{Chen2010ATT}, and \cite{Srivastava2008ATF}.  However, these methods can not be directly applied to high-dimensional compositional data due to the violation of regularity conditions. Therefore, we carry out a systematic examination over the test statistic proposed by \cite{Srivastava2008ATF}, \cite{Srivastava2009ATF} and provide an alternative derivation of the null distribution under a different set of assumptions tailored to compositional data.

Moreover, we establish the asymptotic independence between the max-type and sum-type test
statistics for high-dimensional compositional data. Building upon this finding, we further propose a combined max-sum type test statistic to encompass both sparse and dense situations. To establish the asymptotic independence, we draw inspiration and make use of the results in \cite{Feng2022AsymptoticIO}, which establish the asymptotic independence between the sums and maxima of dependent Gaussian variables. We extend the result of \cite{Feng2022AsymptoticIO} to the degenerate case with semi-definite covariance matrix and derive the null distribution of the proposed max-sum type statistics for high-dimensional compositional data. The central ideas of the proof for our max-sum type test lies in the concentration phenomenon observed in high-dimensional compositional data.

To sum up, in this paper we address both the one-sample and two-sample mean test problems for high-dimensional compositional data. We introduce two testing approaches, namely the sum-type test and max-sum type test (combo-type test), to deal with both dense weak and sparse signals. We derive the asymptotic null distributions and power functions of these test statistics. Notably, we are the first to adopt the sum-type test statistics to high-dimensional compositional data and provide a new way to derive the null distribution which can be linked with the study of max-sum type test statistic. Detailed results can be found in Theorem \ref{one-sample-sum} and Theorem \ref{two-sample-sum}.  
In the context of the max-type test, we supplement the one-sample test and provide a different proof from \cite{Cao2018TwosampleTO}, which enables us to establish its asymptotic independence with the sum-type test statistics.
Most importantly, our proposed max-sum type test represents a new and innovative contribution to current literature (see Theorem \ref{one-sample-max-sum}, \ref{two-sample-max-sum}). Simulation studies demonstrate the superiority
of our max-sum type test statistics which exhibit robust performance regardless of data sparsity. 

This paper is organized as follows. Section \ref{prel} covers preliminaries, encompassing the centred log-ratio transformation of the original data and testable hypotheses. Section \ref{sum} delves into the sum-type test for one-sample/two-sample mean tests for high-dimensional compositional data. Section \ref{max} explores the max-type test and  Section \ref{max-sum-test} concentrates on the asymptotic independence and the max-sum type test. Section \ref{simu} presents numerical studies and detailed proofs are relegated to appendix.

Before moving forward, let us introduce some notations that will be used throughout this paper. We adhere to the convention of using regular letters for scalars and bold-face letters for vectors or matrices. Let $\bA$ be any matrix, with its $(i,j)$-th entry denoted as $A_{ij}$, its transpose as $\bA^T$, its trace as $\tr(\bA)$, and its $j$-th largest eigenvalue as $\lambda_j(\bA)$. The diagonal matrix of $\bA$ is denoted by $\bD_{\bA}$. The norms of $\bA$ are denoted as $\|\bA\|_2 = \sqrt{\lambda_{1}(\bA\bA^T)}$ and  $\|\bA\|_1=\max_{1\leq j\leq p}\sum_{i=1}^p|a_{ij}|$. We write $X_n\stackrel{d}{\to}a$ if $X_n$ converges in distribution to $a$. Furthermore, constants, denoted as $C$ and $K$, may vary between different lines of our analysis.

\section{Preliminaries}\label{prel}
We posit that compositional variables stem from a vector of latent variables referred to as the ``basis". While these basis components are conceived to encapsulate the genuine abundances of bacterial taxa within a microbial community, it is important to note that they remain unobserved. Denote by  $\{\bX_{i}^{(k)}\}_{i=1}^{n_k}$ the observed $n_{k}$ sample  data  for group  $k$ $(k=1,2)$, where the  $\bX_{i}^{(k)}$  represent compositions that lie in the  $(p-1) $-dimensional simplex 
$$\mathcal{S}^{p-1}=   \left\{\left(x_{1}, \ldots, x_{p}\right): x_{j}>0 \text{ ($j=1, \ldots, p$) }, \sum_{j=1}^{p} x_{j}=1\right\}.$$
Let  $\{\bW_{i}^{(k)}\}_{i=1}^{n_k}$  denote the unobserved $n_{k}$ sample basis data, which generate the observed compositional data via the normalization
\begin{equation}\label{coda}
X_{i j}^{(k)}=W_{i j}^{(k)} / \sum_{\ell=1}^{p} W_{i \ell}^{(k)}, \quad\left(i=1, \ldots, n_{k} ; j=1, \ldots, p ; k=1,2\right),
\end{equation}
where $X_{i j}^{(k)}$ and $W_{ij}^{(k)}>0$ are the $j$th components of $\bX_{i}^{(k)}$ and $\bW_{i}^{(k)}$, respectively. 

The unit-sum constraint dictates that compositional variables must not vary independently, rendering many covariance-based multivariate analysis methods inapplicable. To address this, \cite{Jones1982TheSA} proposed a relaxation of the constraint by conducting statistical analysis through log ratios. Among various log-ratio transformations, the centered log-ratio transformation stands out due to its appealing features and widespread adoption. The centred log-ratio transformation of the observed compositional data $\boldsymbol{X}^{(k)}$  are defined by $
\log \left\{X_{ij}^{(k)} / (\prod_{j=1}^{p} X_{ij}^{(k)})^{1 / p}\right\}$.
Considering the scale-invariance of the centered log ratios, we can substitute $X_{ij}^{(k)}$ with $W_{ij}^{(k)}$ and obtain
\begin{equation}\label{YZeq}
       \bG \log \bX_{i}^{(k)} = \bG\log \bW_{i}^{(k)},\quad\left(i=1,\ldots, n_k; k=1,2\right),
\end{equation}
where  $\bG=\bI_{p}-p^{-1} \boldsymbol{1}_{p} \boldsymbol{1}_{p}^{\mathrm{T}}$, $\boldsymbol{1}_p$ is a $p$-dimensional vector of all ones. Suppose that the log basis vectors $\{\log \bW_{i}^{(k)} \}_{i=1}^{n_k}$ $(k=1,2)$ are two independent samples, each from a distribution with mean  $\boldsymbol{\mu}_{k}^W=\left(\mu_{k 1}^W, \ldots, \mu_{k p}^W\right)^{\mathrm{T}} $ and having common covariance matrix  $\boldsymbol{\Sigma}_p=\left(\sigma_{i j}\right)$.
One might consider testing the hypothesis
\begin{equation}\label{hypotest2}
    H_{0}: \boldsymbol{\mu}_{1}^W=\boldsymbol{\mu}_{2}^W  \text { versus } H_{1}: \boldsymbol{\mu}_{1}^W \neq \boldsymbol{\mu}_{2}^W.
\end{equation}
Given that the compositional data $\{\bX_{i}^{(k)}\}_{i=1}^{n_k}$ are observed while the basis data $\{\bW_{i}^{(k)}\}_{i=1}^{n_k}$ are unobserved, rather than testing the hypothesis in \eqref{hypotest2}, we propose testing the following:
\begin{equation}\label{hypotest1}
    H_{0}: \boldsymbol{\mu}_{1}^W=\boldsymbol{\mu}_{2}^W+c \boldsymbol{1}_{p} \text { for some } c \in \mathbb{R} \text { versus } H_{1}: \boldsymbol{\mu}_{1}^W \neq \boldsymbol{\mu}_{2}^W+c \boldsymbol{1}_{p} \text { for any } c \in \mathbb{R}.
\end{equation}
Let $\boldsymbol{\mu}_{k}^X=\Expe(\bG\log\bX_i^{(k)})$ for $k=1,2$. Taking the expectation on both sides of the aforementioned equality \eqref{YZeq}, we obtain $\boldsymbol{\mu}_{k}^X = \mathbf{G}\boldsymbol{\mu}_k^W$.  The matrix  $\bG$  has rank  $p-1$  and hence has a null space of dimension 1. As a result, $ \boldsymbol{\mu}_1^W=\boldsymbol{\mu}_2^W+c \boldsymbol{1}_{p}$  for some  $c \in \mathbb{R}$ if and only if $ \boldsymbol{\mu}_1^X = \boldsymbol{\mu}_2^X$. Hence, the hypothesis in \eqref{hypotest1} is equivalent to the following hypothesis:

\begin{equation}\label{hypotest3-two}
    H_{0}: \boldsymbol{\mu}_1^X=\boldsymbol{\mu}_2^X \quad \text { versus } \quad H_{1}: \boldsymbol{\mu}_1^X \neq \boldsymbol{\mu}_2^X.
\end{equation}
Therefore, we advocate for testing hypothesis in \eqref{hypotest1} by examining the hypothesis in \eqref{hypotest3-two}. We advocate for this approach due to the following reasons:
\begin{itemize}
  \item While the hypotheses in \eqref{hypotest1} and \eqref{hypotest3-two} are equivalent, those in \eqref{hypotest1} hold significance only in the presence of existing bases. In contrast, the hypothesis in \eqref{hypotest3-two} specifically pertains to compositions through the centered log ratios, showcasing their inherent scale-invariance and the ability to be tested using observed compositional data. Therefore, we focus on testing the hypothesis in \eqref{hypotest3-two}, and consequently, those in \eqref{hypotest1}. This approach allows us to test the hypothesis in \eqref{hypotest1} using the observed compositional data $\{\bX_i^{(k)}\}_{i=1}^{n_k}$ for $k = 1, 2$, while the hypothesis in \eqref{hypotest2} are not testable using the observed compositional data.
 \item A basis is determined by its composition only up to a multiplicative factor, and the set of bases giving rise to a composition $\bx \in \mathcal{S}^{p-1}$ forms the equivalence class
$\bW(\bx) = \{(tx_1, \ldots , tx_p) : t > 0\}$ (\cite{Jones2003TheSA}). As an immediate consequence, a log basis vector is determined by the resulting composition only up to an additive constant, and the set of log basis vectors corresponding to $\bx$ constitutes the equivalence class
$\boldsymbol{\log} \bW(\bx) = \{(\log x_1 + c, \ldots , \log x_p + c) : c \in \mathbb{R}\}$. For these reasons, we focus on testing hypothesis in \eqref{hypotest1}.
\item Moreover, $H_0$ in \eqref{hypotest2} implies $H_0$ in \eqref{hypotest1}, meaning that rejecting \eqref{hypotest1} would also result in the rejection of \eqref{hypotest2}.
\end{itemize}
The nature test statistics for testing the hypothesis in \eqref{hypotest3-two}, and hence in \eqref{hypotest1}, would be based on the test statistics designed for $\bG \log \bX_{i}^{(k)}$. For ease of notation, we denote by
\begin{equation}\label{YXeq}
\bY_{i}^{(k)}=
\bG \log \bX_{i}^{(k)}, \quad\left(i=1,\ldots, n_k; k=1,2\right).
\end{equation}
In the end, our goal is to test the hypothesis in \eqref{hypotest1}, which is equivalent to testing the hypothesis in \eqref{hypotest3-two} using the observed compositional data. Note that the hypothesis in \eqref{hypotest1} provide a foundation for deriving biological interpretations related to the genuine abundances.
 
 In the context of one-sample test, let $\left\{\bX_{i}\right\}_{i=1}^n$ be the observed sample compositional data, $\left\{\bW_{i}\right\}_{i=1}^n$ be the unobserved sample basis data, which generate the observed compositional data via the normalization
\begin{equation*}
X_{i j}=W_{i j} / \sum_{\ell=1}^{p} W_{i \ell}, \quad\left(i=1, \ldots, n ; j=1, \ldots, p \right).
\end{equation*}
Define
\begin{equation*}
\bY_{i}=
\bG \log \bX_{i}, \quad\left(i=1,\ldots, n\right).
\end{equation*}
Let
\[\boldsymbol{\mu}^W=\Expe (\log \bW_i),\ \ \boldsymbol{\mu}^X=\Expe (\bG\log \bX_i).\] 
We consider a test of the hypothesis
\begin{equation}\label{hypotest1-one}
     H_{0}: 
   \boldsymbol{\mu}^W=\boldsymbol{\mu}_{0}+c \boldsymbol{1}_{p} \text { for some } c \in \mathbb{R} \text { versus } H_{1}: \boldsymbol{\mu}^W \neq \boldsymbol{\mu}_{0}+c \boldsymbol{1}_{p} \text { for any } c \in \mathbb{R},
\end{equation}
which is equivalent to testing
\begin{equation}\label{hypotest3-one}
    H_{0}: 
  \boldsymbol{\mu}^X=\bG\boldsymbol{\mu}_{0} \quad \text { versus } \quad H_{1}: \boldsymbol{\mu}^X \neq \bG\boldsymbol{\mu}_{0}.
\end{equation}
The test statistics for testing the hypothesis in \eqref{hypotest3-one}, and hence in \eqref{hypotest1-one}, would be based on the test statistics designed for $\bY_i$.

\section{Sum-type mean test}\label{sum}

In this section, our focus is on statistical inference issues concerning high-dimensional compositional data. Currently, there is a lack of a sum-of-squares-type test statistic suitable for dense-type alternative hypotheses in both one-sample and two-sample mean tests for high-dimensional compositional data. To address this gap, we initially present a sum-type test statistic in the form of a sum of squares. More importantly, the theoretical analysis of the asymptotic distribution for the sum-type test statistic is intended for use in proving our novel max-sum type test proposed in Section \ref{max-sum-test}. As we aim to test the latent basis structures, we will impose conditions directly on the log basis variables. Assuming a common basis covariance matrix, the two populations share a common centered log-ratio covariance matrix $\boldsymbol{\Gamma}_p=\operatorname{cov}\left(\bY_i^{(k)}\right)(k=1,2)$, which, according to \eqref{YZeq}, is expressed as
$$
\boldsymbol{\Gamma}_p=\bG \boldsymbol{\Sigma}_p \bG^{T} .
$$

Let $\bR_p=\left(\rho_{i j}\right)$ denote the correlation matrices of $\log\bW_i^{(k)}$, $i=1,\ldots,n_k$, $k=1,2$, respectively. Let $\lambda_i(\bR_{p})$, $i=1,\ldots,p$ be the eigenvalues of the correlation matrix $\bR_p$. We initially impose  conditions on the covariance structures of the log basis variables as follows:
\begin{Assumpa}\label{distriassume}
    $\log\bW_i^{(k)}\stackrel{\text{i.i.d}}{\sim}N(\boldsymbol{\mu}^W_k,\boldsymbol{\Sigma}_p)$ for $i=1,\ldots,n_k$, $k=1,2$. (In one-sample case, $\log\bW_i\stackrel{\text{i.i.d}}{\sim}N(\boldsymbol{\mu}^W,\boldsymbol{\Sigma}_p)$ for $i=1,\ldots,n$)
\end{Assumpa}

\begin{Assumpb}\label{sumassume}
     $\lim_{p\to\infty}\frac{\tr(\bR_{p}^{i})}{p}<\infty$ for  $i=2,4$, $\max_{1\leq i<j\leq p }|\rho_{ij}|=O(p^{-1/2-\alpha}(\log p)^C)$ for $\alpha>0$ and some constant $C>0$; $N=O(p^{\epsilon})$ for $\epsilon \in (\frac{1}{2}, 1]$, where $N=n_1+n_2$. (In one-sample case, $n=O(p^{\epsilon})$ for $\epsilon \in (\frac{1}{2}, 1]$.)
\end{Assumpb}

\begin{Assumpc}\label{addassume1}
    $1/\kappa_1\leq\sigma_{ii}\leq \kappa_1$ for $1\leq i\leq p$ and some constant $\kappa_1>0$.
\end{Assumpc}

\subsection{One sample test: Gaussian case}
For the one-sample test for high-dimensional compositional data, we propose the sum-type test statistics
\begin{align}\label{sumst1}
    T_{\tsum}=\frac{n\bar{\bY}^{T} {\bD}_{\hat{\boldsymbol{\Gamma}}_p}^{-1}\bar{\bY}-(n-1)p/(n-3)}{\sqrt{2[\tr({\bD}_{\hat{\boldsymbol{\Gamma}}_p}^{-1/2}\hat{\boldsymbol{\Gamma}}_p{\bD}_{\hat{\boldsymbol{\Gamma}}_p}^{-1/2})^2-p^2/(n-1)]}},
\end{align}
where $\bar{\bY}=\frac{1}{n}\sum_{i=1}^{n}\bY_i$,  ${\bD}_{\hat{\boldsymbol{\Gamma}}_p}$ is the diagonal matrix of the sample covariance matrix $\hat{\boldsymbol{\Gamma}}_p=\frac{1}{n}\sum_{i=1}^{n}(\bY_i-\bar{\bY})^T(\bY_i-\bar{\bY})$, and ${\bD}_{\hat{\boldsymbol{\Gamma}}_p}^{-1/2}\hat{\boldsymbol{\Gamma}}_p{\bD}_{\hat{\boldsymbol{\Gamma}}_p}^{-1/2}$ is the sample correlation matrix. The main component of $T_{\text{sum}}$ can be expressed as a sum of random variables. Note that, similar to the proof of Theorem 2.1  in \cite{Srivastava2008ATF}, we have

$$T_{\tsum}=\tilde{\tilde{T}}_{\tsum}+o_p(1):=\frac{n\boldsymbol{\bar{Y}}^T\bD_{\boldsymbol{\Gamma}_p}^{-1}\boldsymbol{\bar{Y}}-p}{\sqrt{2\tr({\bD}_{{\boldsymbol{\Gamma}}_p}^{-1/2}{\boldsymbol{\Gamma}}_p{\bD}_{{\boldsymbol{\Gamma}}_p}^{-1/2})^2}}+o_p(1),$$
where  $\bD_{\boldsymbol{\Gamma}_p}$ is the diagonal matrix of $\boldsymbol{\Gamma}_p$. Moreover, as proved in Section \ref{pfone-sample-sum} we have
\begin{align*}
\tilde{\tilde{T}}_{\tsum}=\tilde{T}_{\tsum}+o_p(1):=\frac{\bV^T\bV-tr(\bG\bR_p\bG)}{\sqrt{2\tr(\bG\bR_p\bG)^2}}+o_p(1),
\end{align*}
where 
\begin{align}\label{eqV}
    \bV:=\sqrt{n}\bG\bD_{\boldsymbol{\Sigma}_p}^{-1/2}(\bar{\log\bW}-c\boldsymbol{1}_p)\sim N(\boldsymbol{0},\bG\bR_p\bG)\ \ \text{under} \ \ H_0,
\end{align}
and $\bar{\log\bW}=\frac{1}{n}\sum_{i=1}^n\log\bW_i$ and $\bD_{\boldsymbol{\Sigma}_p}$ is the diagonal matrix of $\boldsymbol{\Sigma}_p$. To derive the asymptotic distribution for $\tilde{T}_{\tsum}$, we establish the following technical lemma and its proof is postponed to Section \ref{pftheorem_sum}.

\begin{lemma}\label{theorem_sum}  
If $\bZ\sim N(c\boldsymbol{1}_p,\boldsymbol{\Psi}_p)$, $\bH=\bG \bZ=(H_1,\ldots,H_p)^T$, and $\lim_{p\to \infty}\frac{\tr(\boldsymbol{\Psi}_{p}^{2})}{p}<\infty$, then as $p\to \infty$,
$$\frac{H_1^2+\cdots + H_p^2-\tr(\bG\boldsymbol{\Psi}_p\bG)}{\sqrt{2\tr(\bG\boldsymbol{\Psi}_p\bG)^2}} \stackrel{d}{\to} N(0, 1).$$
\end{lemma}
\begin{Remark}
In Lemma \ref{theorem_sum},  our study focuses on Gaussian data with a non-negative definite covariance matrix $\bG\boldsymbol{\Psi}_p\bG$. In contrast, the results presented in \cite{Feng2022AsymptoticIO} pertain to Gaussian data with a strictly positive definite covariance matrix $\boldsymbol{\Psi}_p$. It's important to note that our findings are applicable to high-dimensional compositional data, whereas the results in \cite{Feng2022AsymptoticIO} are not suited for such data. Additionally, our theoretical analysis is different from \cite{Feng2022AsymptoticIO}.
\end{Remark} 
By Lemma \ref{theorem_sum}, we conclude that as $p\to \infty$
$$\tilde{T}_{\tsum}:=\frac{\bV^T\bV-tr(\bG\bR_p\bG)}{\sqrt{2tr(\bG\bR_p\bG)^2}}\to N(0,1).$$
Thus, we derive Theorem \ref{one-sample-sum} as follows. The detailed proof of Theorem \ref{one-sample-sum} is postponed to Section \ref{pfone-sample-sum}.
\begin{theorem}\label{one-sample-sum}
If Assumptions \ref{distriassume}-\ref{addassume1} holds, under the null hypothesis in \eqref{hypotest1-one}, as $p\to\infty$ then $$T_{\tsum}\stackrel{d}{\to} N(0, 1).$$ 
\end{theorem}
\begin{Remark}
    \begin{enumerate}
        \item[(1)] The sum-type test $T_{\tsum}$ discussed in \cite{Feng2022AsymptoticIO} can be extended to high-dimensional compositional data. Our theoretical analysis differs, as we only require $\lim_{p\to\infty}\frac{\text{tr}(\bR_{p}^{i})}{p}<\infty$ for $i=2,4$, while \cite{Feng2022AsymptoticIO} requires it for $i=2,3,4$.
        \item[(2)] It is important to note that, although $T_{\text{sum}}$ shares the same formulation as the test statistic defined in \cite{Srivastava2008ATF}, our theoretical analysis differs. Our analysis relies on central limit theory for random quadratic forms in \cite{wang2014jointclt}. Furthermore, we extend the result of \cite{wang2014jointclt} to our situation. This extension is necessary because, to apply the results in \cite{wang2014jointclt}, assumptions must be imposed on the centered log-ratio transformation of the observed compositional data. This is not directly applicable to our situation, where we impose assumptions directly on the log-basis variables, as our objective is to test the latent basis structures.  Crucially, our theoretical analysis in the derivation of the asymptotic distribution for $T_{\text{sum}}$ is intended to be used for the proof of our max-sum type test in Section \ref{max-sum-test}. Furthermore, it is pertinent to highlight that \cite{Srivastava2009ATF} have explored the behavior of $T_{\tsum}$ even under non-Gaussian data. 
    \end{enumerate}
\end{Remark}
Next, we consider the asymptotic power of $T_{\tsum}$ under local
alternatives 
\begin{align}\label{onesumH1}
   H_1: \boldsymbol{\mu}^W = \boldsymbol{\delta}_{\tsum}+c\boldsymbol{1}_p,
\end{align}
where $\boldsymbol{\delta}_{\tsum}=\frac{1}{\sqrt{n(n-1)}}\boldsymbol{\delta}_{s}$, $\boldsymbol{\delta}_{s}$ is a vector of constants, $c \in \mathbb{R}$, $\boldsymbol{1}_p$ is a $p$-dimensional vector of all ones.
We further assume that 
\begin{Assumpd}\label{dddm}
 $\frac{\boldsymbol{\delta}_{s}'\bG\bD_{\boldsymbol{\Gamma}_p}^{-1}\bG\boldsymbol{\delta}_{s}}{p}\leq M$ for every $p$, where $M$ does not depend on $p$ and $\bD_{\boldsymbol{\Gamma}_p}$ is the diagonal matrix of $\boldsymbol{\Gamma}_p=\bG\boldsymbol{\Sigma}_p\bG$.   
\end{Assumpd}
The following theorem gives the asymptotic distribution of $T_{\tsum}$ under the local
alternative in \eqref{onesumH1}, which is derived similarly to Theorem 2.2 in \cite{Srivastava2008ATF}, then its proof omitted.
\begin{proposition}\label{onesumal}
    Assume that  Assumptions \ref{distriassume}-\ref{addassume1} and \ref{dddm} hold. Under the local alternative in \eqref{onesumH1}, the power function of $T_{sum}$ is
    \begin{align*}
      \beta_S(\boldsymbol{\mu}^W,\alpha)=  \lim_{p\to\infty}\Phi(-z_{1-\alpha}+\frac{n(\boldsymbol{\mu}^W)^T\bG\bD_{\boldsymbol{\Gamma}_p}^{-1}\bG\boldsymbol{\mu}^W}{\sqrt{2\tr({\bD}_{\boldsymbol{\Gamma}_p}^{-1/2}\boldsymbol{\Gamma}_p{\bD}_{\boldsymbol{\Gamma}_p}^{-1/2})^2}}),
    \end{align*}
where $z_{1-\alpha}$ is the $(1-\alpha)-$quantile of $N(0,1)$.
\end{proposition}
\begin{Remark}
 The asymptotic distribution of $T_{\tsum}$ for high-dimensional compositional data under the local alternative in \eqref{onesumH1} is the same as Theorem 2.2 in \cite{Srivastava2008ATF}.

\end{Remark}
\subsection{Two sample test: Gaussian case}
For the two-sample test for high-dimensional compositional data, we propose the sum-type test statistics
\begin{align}\label{sumst2}
    T_{\tsum,2}=\frac{\frac{n_1n_2}{n_1+n_2}(\bar{\bY}^{(1)}-\bar{\bY}^{(2)})^{T} \bD_{\hat{\bar{\boldsymbol{\Gamma}}}_p}^{-1}(\bar{\bY}^{(1)}-\bar{\bY}^{(2)})-\frac{(n_1+n_2-2)p}{(n_1+n_2-4)}}
{\sqrt{2\big[\tr(  {\bD}_{\hat{\bar{\boldsymbol{\Gamma}}}_p}^{-1/2} \hat{\bar{\boldsymbol{\Gamma}}}_p{\bD}_{\hat{\bar{\boldsymbol{\Gamma}}}_p}^{-1/2})^2-\frac{p^2}{(n_1+n_2-2)}\big]c_{p,N}}},
\end{align}
where $\bar{\bY}^{(k)}=\frac{1}{n_k}\sum_{i=1}^{n_k}\bY_i^{(k)}$,  $\bD_{\hat{\bar{\boldsymbol{\Gamma}}}_p}$ is the diagonal matrix of sample covariance matrix 
\begin{align}\label{twosamcov}
\hat{\bar{\boldsymbol{\Gamma}}}_p=\frac{1}{n_1+n_2}\Big[\sum_{i=1}^{n_1}
(\bY_{i}^{(1)}-\bar{\bY}^{(1)})(\bY_{i}^{(1)}-\bar{\bY}^{(1)})^T+\sum_{i=1}^{n_2}
(\bY_{i}^{(2)}-\bar{\bY}^{(2)})(\bY_{i}^{(2)}-\bar{\bY}^{(2)})^T\Big],
\end{align}
and  $c_{p,N}=1+\frac{\tr({\bD}_{\hat{\bar{\boldsymbol{\Gamma}}}_p}^{-1/2} \hat{\bar{\boldsymbol{\Gamma}}}_p{\bD}_{\hat{\bar{\boldsymbol{\Gamma}}}_p}^{-1/2})^2}{p^{3/2}}$.  Similarly to the discussion for the one-sample case above Theorem \ref{one-sample-sum}, we derive Theorem \ref{two-sample-sum} as follows and its proof is postponed to Section \ref{pftwo-sample-sum}.

\begin{theorem}\label{two-sample-sum}
If Assumptions \ref{distriassume}-\ref{addassume1} hold and $\lim_{p\to\infty}\frac{n_1}{n_2}=\kappa\in (0,\infty)$, under the null hypothesis in $\eqref{hypotest1}$, as $p\to\infty$ then \begin{equation*}
    T_{\mathrm{sum},2}\stackrel{d}{\to} N(0, 1).
\end{equation*}
\end{theorem}
Next, we deal with the behavior of $T_{\tsum,2}$ under local
alternatives 
\begin{align}\label{twosumH1}
   H_1: \boldsymbol{\mu}_1^W-\boldsymbol{\mu}_2^W = \boldsymbol{\delta}_{\tsum,2}+c\boldsymbol{1}_p,
\end{align}
where $\boldsymbol{\delta}_{\tsum,2}=(\frac{n_1+n_2}{(n_1+n_2-2)n_1n_2})^{1/2}\boldsymbol{\delta}_{s}$, $\boldsymbol{\delta}_{s}$ is a vector of constants and satisfies the condition \eqref{dddm}. 
The following theorem gives the asymptotic distribution of $T_{\tsum,2}$ under the local
alternative in \eqref{twosumH1}, which is derived similarly to Theorem 2.2 in \cite{Srivastava2008ATF} and its proof is omitted.

\begin{proposition}\label{twosumal}
    Assume that Assumptions $\ref{distriassume}$-\ref{addassume1} and $\ref{dddm}$ hold. Under the local alternative in \eqref{twosumH1}, the power function of $T_{\tsum,2}$ is
    \begin{align*}
        \beta_{S,2}(\boldsymbol{\mu}_1^W,\boldsymbol{\mu}_2^W,\alpha)=\lim_{p\to\infty}\Phi(-z_{1-\alpha}+\frac{\frac{n_1n_2}{n_1+n_2}(\boldsymbol{\mu}_1^W-\boldsymbol{\mu}_2^W)^T\bG\bD_{\boldsymbol{\Gamma}_p}^{-1}\bG(\boldsymbol{\mu}_1^W-\boldsymbol{\mu}_2^W)}{\sqrt{2\tr({\bD}_{\boldsymbol{\Gamma}_p}^{-1/2}\boldsymbol{\Gamma}_p{\bD}_{\boldsymbol{\Gamma}_p}^{-1/2})^2}}),
    \end{align*}
where $z_{1-\alpha}$ is the $1-\alpha-$quantile of $N(0,1)$.
\end{proposition}

\subsection{Non-Gaussian cases}
In the preceding sections, we consistently assume that the log basis vectors follow a Gaussian distribution. This assumption is necessary to ensure the independence of the sum-type and max-type test statistics for deriving the asymptotic distribution of the combo-type test statistics. However, if we exclusively use the sum-type test statistics, the same asymptotic normality can be achieved under more general assumptions, as outlined below. 
\begin{Assumpaa}\label{assumpa2}
        $\log\bW_i^{(k)}=\boldsymbol{\mu}^W_k + \boldsymbol{\Sigma}_p^{1/2}\bU^{(k)}_i$ where $\bU_i^{(k)} = (U_{i1}^{(k)},\dots,U_{ip}^{(k)})^T$ and $U_{ij}^{(k)}$, $i=1,\dots,n_k$, $j=1,\dots,p$, $k=1,2$,  are i.i.d with zero mean, variance 1 and finite fourth moment $\theta$. (In one sample case, $\log\bW_i=\boldsymbol{\mu}^W + \boldsymbol{\Sigma}_p^{1/2}\bU_i$ where $\bU_i = (U_{i1},\dots,U_{ip})^T$ and $U_{ij}$, $i=1,\dots,n$, $j=1,\dots,p$, are i.i.d with zero mean, variance 1 and finite fourth moment $\theta$.)
\end{Assumpaa} 
\begin{Assumpbb}\label{assumpbb1}
        $\lim\limits_{p\to\infty} \frac{\tr (\bR_p^i)}{p} < \infty$, $i=2,4$, $||\bR_p||_{1} = o(p)$, and $N=O(p^{\epsilon})$, $\frac{1}{2} <\epsilon \leq 1$.
\end{Assumpbb}

\subsubsection*{One sample test: non-Gaussian case}

In the one-sample case, we utilize the statistic $T_{\text{sum}}$ defined by \eqref{sumst1}. The proof of Theorem \ref{onesample-sum-nongaussian} is postponed to Section \ref{pfonesample-sum-nongaussian}.
\begin{theorem}\label{onesample-sum-nongaussian}
Under null hypothesis in \eqref{hypotest1-one}
$$
H_{0}: 
   \boldsymbol{\mu}^W=\boldsymbol{\mu}_{0}+c \boldsymbol{1}_{p} \text { for some } c \in \mathbb{R},
$$
if Assumptions \ref{assumpa2}, \ref{assumpbb1} and \ref{addassume1} hold, as $p \to\infty$
    $$T_{\tsum}\stackrel{d}{\to} N(0, 1).$$ 
\end{theorem}
\begin{Remark}
\begin{enumerate}
    \item[(i)] Although  Theorem \ref{onesample-sum-nongaussian} looks formally the same as Theorem 3.1 in \cite{Srivastava2009ATF} , the theoretical analysis is radically different, in that our test statistic is designed based on transformed log basis vector $\bG\log W_i$ in the limitations of the component data. Moreover, conditions of Theorem \ref{onesample-sum-nongaussian} are weaker than the conditions of Theorem 3.1 in \cite{Srivastava2009ATF}.
    \item[(ii)] In contrast to Assumption \ref{distriassume}, Assumption \ref{assumpa2} discards the restriction of the Gaussian distribution and only requires the existence of fourth-order moments. Compared to Assumption \ref{sumassume}, Assumption \ref{assumpbb1} also has a weaker requirement for correlation matrix $\bR_p$.
\end{enumerate}

\end{Remark}
Next we consider the local alternative in \eqref{onesumH1}, and we have the following power analysis,  which is derived similarly to Theorem 4.1 in \cite{Srivastava2009ATF} and its proof is omitted.
\begin{proposition}
     Assume that Assumptions \ref{assumpa2}, \ref{assumpbb1}, \ref{addassume1} and \ref{dddm} hold. Under the local alternative in \eqref{onesumH1}, the power function of $T_{\tsum}$ is 
    \begin{align*}
      \beta_S(\boldsymbol{\mu}^W,\alpha)=  \lim_{p\to\infty}\Phi(-z_{1-\alpha}+\frac{n(\boldsymbol{\mu}^W)^T\bG\bD_{\boldsymbol{\Gamma}_p}^{-1}\bG\boldsymbol{\mu}^W}{\sqrt{2\tr({\bD}_{\boldsymbol{\Gamma}_p}^{-1/2}\boldsymbol{\Gamma}_p{\bD}_{\boldsymbol{\Gamma}_p}^{-1/2})^2}}),
    \end{align*}
where $z_{1-\alpha}$ is the $(1-\alpha)-$quantile of $N(0,1)$.
\end{proposition}

\subsubsection*{Two sample test: non-Gaussian case}

In the two-sample case, we utilize the statistic $T_{\text{sum},2}$  defined by \eqref{sumst2}. We commence by extending the results of \cite{Srivastava2009ATF} in one sample-case to the two-sample case under more mild conditions, subsequently deriving results applicable to compositional data. The following lemma establishes the consistency of estimators, allowing them to be replaced by information from the population when deriving the asymptotic distribution. The proof of Lemma \ref{twosample-consistencycor} is postponed to Section \ref{pftwosample-consistencylemma}.
\begin{lemma}\label{twosample-consistencycor}
    Let $\hat{\gamma}_{ii}$ be diagonal elements of sample covariance $ \hat{\bar{\boldsymbol{\Gamma}}}_p$, where $ \hat{\bar{\boldsymbol{\Gamma}}}_p$ is define by \eqref{twosamcov}, if Assumptions \ref{assumpa2}, \ref{assumpbb1} and \ref{addassume1} hold, then $\hat{\gamma}^{-1}_{ii}$ are consistent estimators of $\gamma_{ii}^{-1}$ and $\frac{1}{p}[\tr(  {\bD}_{\hat{\bar{\boldsymbol{\Gamma}}}_p}^{-1/2} \hat{\bar{\boldsymbol{\Gamma}}}_p{\bD}_{\hat{\bar{\boldsymbol{\Gamma}}}_p}^{-1/2})^2-(p^2/(N-2))]$ is a consistent estimator of $\tr ({\bD}_{\boldsymbol{\Gamma}_p}^{-1/2}\boldsymbol{\Gamma}_p{\bD}_{\boldsymbol{\Gamma}_p}^{-1/2})^2/p$.
\end{lemma}
\begin{Remark}
 Lemma \ref{twosample-consistencycor} establishes the consistency of the estimators in $T_{\tsum,2}$, allowing them to be replaced by their population counterparts when deriving the asymptotic distribution. As \cite{Srivastava2008ATF} has discussed the consistency of the estimators in both the one-sample and two-sample cases for the Gaussian scenario, we derive the consistency of the estimators in $T_{\tsum,2}$ for the non-Gaussian scenario similarly to \cite{Srivastava2009ATF} and its proof is postponed to Section \ref{pfone-sample-sum}.
\end{Remark}

The following theorem shows that the asymptotic distribution of the test statistic, constructed from the centered log-ratio vector $\bY_i^{(k)}$, specifically designed for compositional data. The proof of Theorem \ref{twosamplemaintheorem} is postponed to Section \ref{pftwosamplemaintheorem}.
\begin{theorem}\label{twosamplemaintheorem}
    Under the null hypothesis in $\eqref{hypotest1}$
    $$
    H_{0}: \boldsymbol{\mu}_{1}^W=\boldsymbol{\mu}_{2}^W+c \boldsymbol{1}_{p} \text { for some } c \in \mathbb{R},
$$
if Assumptions \ref{assumpa2}, \ref{assumpbb1} and \ref{addassume1} hold and $\lim_{p\to\infty}\frac{n_1}{n_2}=\kappa\in (0,\infty)$, as $p\to\infty$ then \begin{equation*}
    T_{\tsum,2}\stackrel{d}{\to} N(0, 1).
\end{equation*}
\end{theorem}
\begin{Remark}
    Theorem \ref{twosamplemaintheorem}, in which the test statistic is of the same form as \cite{Srivastava2008ATF}, holds under weaker conditions compared with \cite{Srivastava2009ATF}.
\end{Remark}
Next we consider the local alternative in \eqref{twosumH1}, we have the following power analysis which is derived similarly to Theorem 4.1 in \cite{Srivastava2009ATF} and its proof is postponed to \ref{pftwosample-power}.
\begin{proposition}\label{twosample-power}
    Assume that Assumptions \ref{assumpa2}, \ref{assumpbb1}, \ref{addassume1} and $\ref{dddm}$ hold. Under the local alternative in \eqref{twosumH1}, the power function of $T_{sum,2}$ is
    \begin{align*}
        \beta_{S,2}(\boldsymbol{\mu}_1^W,\boldsymbol{\mu}_2^W,\alpha)=\lim_{p\to\infty}\Phi(-z_{1-\alpha}+\frac{\frac{n_1n_2}{n_1+n_2}(\boldsymbol{\mu}_1^W-\boldsymbol{\mu}_2^W)^T\bG\bD_{\boldsymbol{\Gamma}_p}^{-1}\bG(\boldsymbol{\mu}_1^W-\boldsymbol{\mu}_2^W)}{\sqrt{2\tr({\bD}_{\boldsymbol{\Gamma}_p}^{-1/2}\boldsymbol{\Gamma}_p{\bD}_{\boldsymbol{\Gamma}_p}^{-1/2})^2}}),
    \end{align*}
where $z_{1-\alpha}$ is the $1-\alpha-$quantile of $N(0,1)$.
\end{proposition}

\section{Max-type mean test}\label{max}

Following the preceding section, for one-sample and two-sample mean tests on high-dimensional compositional data, we introduce a max-type statistic designed specifically for sparse-type alternative hypotheses. It's worth noting that, there has been no previous development of a max-type statistic for one-sample mean tests of high-dimensional compositional data. In contrast, \cite{Cao2018TwosampleTO} has previously explored max-type statistics for two-sample mean tests in his work. It is important to highlight that, although our max-type statistics for two-sample mean tests are the same as the test statistic defined in Cao's work, our theoretical analysis differs significantly and is intricately linked with the proof for our novel max-sum type test in Section \ref{max-sum-test}.  To proceed, we need more notations. For two sequences of numbers $\{a_p\geq 0;\, p\geq 1\}$ and $\{b_p>0;\, p\geq 1\}$, we write $a_p\ll b_p$ if $\lim_{p\to\infty}\frac{a_p}{b_p}=0$.  The following assumption will be imposed:
\begin{Assumpb}\label{assumption_max0}
    For some $\varrho\in (0,1)$, assume $|\rho_{ij}|\leq \varrho $  for all $1\leq i<j \leq p$ and $p\geq 2$, and $\max_{1\leq i \leq p}\sum_{j=1}^p\rho_{ij}^2\leq C$ for some constant $C> 0$, $\lambda_{\max}(\bR_p)\ll \sqrt{p}(\log p)^{-1}$. And $\log p =o(N^{1/3})$. (In one-sample case, $\log p =o(n^{1/3})$.)
\end{Assumpb}

\subsection{One sample test}

In this section, we propose the max-type test statistics
\begin{align}\label{maxst1}
    T_{\max}=n\max_{1\leq i\leq p}\frac{\bar{Y}_i^2}{\widehat{\gamma}_{ii}},
\end{align}
where $\bar{Y}_i$ is the ith coordinate of $\bar{\bY}=\frac{1}{n}\sum_{i=1}^n\bY_i$, $\widehat{\gamma}_{ii}$ is the $i$-th sample variance in $\hat{\boldsymbol{\Gamma}}_p$. Note that, according to the proof of Theorem 4 in \cite{Feng2022AsymptoticIO},we have
$$T_{\max}=\tilde{\tilde{T}}_{\max}+o_p(1):=n\max_{1\leq i\leq p}\frac{\bar{Y}_i^2}{\gamma_{ii}}+o_p(1),$$
where $\gamma_{ii}$ is the $i$-th population variance in $\boldsymbol{\Gamma}_p$. Moreover, as proved in Section \ref{pfone-sample-max} we have
$$\tilde{\tilde{T}}_{\max}=\tilde{T}_{\max}+o_p(1),$$
where $\tilde{T}_{\max}=\max_{1\leq i \leq p}V_i^2$ and $\bV=(V_1,\cdots,V_p)^T$ is in \eqref{eqV}. To derive the asymptotic distribution for $\tilde{T}_{\max}$, we establish the following technical lemma and its proof is postponed to Section \ref{pftheorem_max}.
\begin{lemma}\label{theorem_max} 
If $\bZ\sim N(c\boldsymbol{1}_p,\boldsymbol{\Psi}_p)$, $\bH=\bG\bZ$, and the following condition hold:
\begin{itemize}
    \item
 Let $\boldsymbol{\Psi}_p=(\psi_{ij})_{1\le i,j\leq p}$. For some $\varrho\in (0,1)$, assume $|\psi_{ij}|\leq \varrho $  for all $1\leq i<j \leq p$ and $p\geq 2$. Suppose $\{\delta_p;\, p\geq 1\}$  and  $\{\kappa_p;\, p\geq 1\}$ are  positive constants with 
 $\delta_p=o(1/\log p)$ and $\kappa=\kappa_p\to 0$
as $p\to\infty$. For $1\leq i \leq p$, define 
$ B_{p,i}=\big\{1\leq j \leq p;\, |\psi_{ij}|\geq \delta_p\big\}$ and  $C_p=\big\{1\leq i \leq p;\, |B_{p,i}|\geq p^{\kappa}\big\}$. Assume that  $|C_p|/p\to 0$ as $p \to\infty$. And $\lambda_{\max}(\boldsymbol{\Psi}_p)\ll \sqrt{p}(\log p)^{-1}$.
\end{itemize}
Then as $p\to \infty$
\begin{align*}
    \max_{1\leq i \leq p}H_i^2-2\log p +\log\log p
\end{align*}
converges to a Gumbel distribution with cdf $F(x)=\exp\{-\frac{1}{\sqrt{\pi}}e^{-x/2}\}$.
\end{lemma}

\begin{Remark}
     Our study focuses on Gaussian data characterized by a non-negative definite covariance matrix $\bG\boldsymbol{\Psi}_p\bG$. In contrast, the results presented in \cite{Feng2022AsymptoticIO} apply to Gaussian data with a strictly positive definite covariance matrix $\boldsymbol{\Psi}_p$. It is crucial to highlight that our findings are specifically tailored for high-dimensional compositional data, while the results in \cite{Feng2022AsymptoticIO} are not suitable for such data. In addition to this distinction, we also require an additional condition $\lambda_{\text{max}}(\boldsymbol{\Psi}_p) \ll \sqrt{p}(\log p)^{-1}$ when dealing with high-dimensional compositional data.
\end{Remark}
By Lemma \ref{theorem_max} we have under Assumptions \ref{distriassume}, \ref{assumption_max0}, $$\tilde{T}_{\max}-2\log p +\log\log p$$
converges to a Gumbel distribution with cdf $F(x)=\exp\{-\frac{1}{\sqrt{\pi}}e^{-x/2}\}$ as $p\to \infty$. Thus we derive Theorem \ref{one-sample-max} as follows and the detailed proof of Theorem \ref{one-sample-max} is postponed to Section \ref{pfone-sample-max}.
\begin{theorem}\label{one-sample-max}
If Assumptions \ref{distriassume}, \ref{assumption_max0},  \ref{addassume1} hold, under the null hypothesis in \eqref{hypotest1-one}, then as $p\to\infty$
$$T_{\max}-2\log p +\log\log p$$
converges to a Gumbel distribution with cdf $F(x)=\exp\{-\frac{1}{\sqrt{\pi}}e^{-x/2}\}$.
\end{theorem}

\begin{Remark}

\begin{enumerate}
    \item[(1)] The max-type test statistics $T_{\max}$ used for one-sample testing in the context of high-dimensional compositional data, exhibits an identical distribution to the case discussed in \cite{Feng2022AsymptoticIO}.
    \item[(2)] To date, there is a dearth of research on one-sample tests specifically tailored for high-dimensional compositional data.
\end{enumerate}

\end{Remark}

Next, we deal with the behavior of $T_{\max}$ under local
alternatives 
\begin{align}\label{onemaxH1}
   H_1: \mu_j^W\neq c, \ \ j\in\bS; \ \ \mu_j^W=c,\ \ j\in\bS^c
\end{align}
for some $c\in\mathbb{R}$ and $\bS\subset \{1,\cdots,p\}$ with cardinality $s$, where $\bS^c$ denotes the complement of $\bS$. We first define the asymptotic $\alpha-$ level test
\begin{align*}
    \Phi_{\alpha}=I(T_{\max}\geq q_{\alpha}+2\log p-\log\log p),
\end{align*}
where $q_{\alpha}=-\log\pi -2\log\log (1-\alpha)^{-1}$ is the $(1-\alpha)$-quantile of the Gumbel distribution. Define the signal vector $\boldsymbol{\delta}_{\max}=(\frac{\log p}{n})^{1/2}(\sigma_{11}\delta_{m,1},\cdots,\sigma_{pp}\delta_{m,p})^T$ by
\begin{align*}
     \boldsymbol{\mu}^W=\boldsymbol{\delta}_{\max}+c\boldsymbol{1}_p.
\end{align*}
The following theorem gives the asymptotic power of $T_{\max}$, which is derived similar to Theorem 2 in \cite{Cao2018TwosampleTO}. Thus the proof is omitted here. 
\begin{proposition}\label{onemaxal}
   Assume that Assumptions \ref{distriassume}, \ref{assumption_max0}, \ref{addassume1} hold. Under the local alternative in \eqref{onemaxH1}, if $\sum_{j=1}^p|\delta_{m,j}|=o(p)$ and $\max_{j\in\bS}|\delta_{m,j}|\geq \sqrt{2}+\epsilon$ for some constant $\epsilon>0$, then $pr(\Phi_\alpha=1)\to1$ as $n,p\to\infty$.
\end{proposition}

\subsection{Two sample test}
For the two-sample test, we propose the max-type test statistics
\begin{align}\label{maxst2}
    T_{\max,2}=\frac{n_1n_2}{n_1+n_2}\max_{1\leq i\leq p}\frac{(\bar{Y}_{i}^{(1)}-\bar{Y}_{i}^{(2)})^2}{\widehat{\gamma}_{ii}},
\end{align}
where $\bar{Y}_{i}^{(k)}$ is the $i$th coordinate of $\bar{\bY}^{(k)}\in \mathbb{R}^p$, and $\bar{\bY}^{(k)}=\frac{1}{n_k}\sum_{i=1}^{n_k}\bY_i^{(k)}$,  $\hat{\gamma}_{ii}$ is the $i$th sample variance in $\hat{\bar{\boldsymbol{\Gamma}}}_p$ defined in (\ref{twosamcov}). Similar to the discussion above Theorem \ref{one-sample-max}, we derive Theorem \ref{two-sample-max} as follows and the detailed proof of Theorem \ref{two-sample-max} is postponed to Section \ref{pftwo-sample-max}. 

\begin{theorem}\label{two-sample-max}
Under the null hypothesis in \eqref{hypotest1} and the same conditions of Theorem \ref{one-sample-max}, then as $p\to\infty$
$$T_{\max,2}-2\log p +\log\log p$$ 
converges to a Gumbel distribution with cdf $F(x)=\exp\{-\frac{1}{\sqrt{\pi}}e^{-x/2}\}$.
\end{theorem}
\begin{Remark}
Note that, although the $T_{\max,2}$ is the same to the test statistic defined in \cite{Cao2018TwosampleTO}, our theoretical analysis is radically different and is intricately linked with the proof for our novel max-sum type test.

The behavior of $T_{\max,2}$ under local
alternatives 
\begin{align}\label{twomaxH1}
  H_1: \mu_{1j}^W\neq \mu_{2j}^W+c, \ \ j\in\bS; \ \ \mu_{1j}^W=\mu_{2j}^W+c,\ \ j\in\bS^c,
\end{align}
is the same as Theorem 2 in \cite{Cao2018TwosampleTO}. We first define the asymptotic $\alpha-$ level test 
\begin{align*}
    \Phi_{\alpha,2}=I(T_{\max,2}\geq q_{\alpha}+2\log p-\log\log p),
\end{align*}
where $q_{\alpha}=-\log\pi -2\log\log (1-\alpha)^{-1}$ is the $(1-\alpha)$-quantile of the Gumbel distribution. Define the signal vector $\boldsymbol{\delta}_{\max,2}=(\frac{\log p}{N})^{1/2}(\sigma_{11}\delta_{m,1},\cdots,\sigma_{pp}\delta_{m,p})^T$ by
\begin{align*}
     \boldsymbol{\mu}_{1}^W-\boldsymbol{\mu}_{2}^W=\boldsymbol{\delta}_{\max,2}+c\boldsymbol{1}_p, \ \ j=1,\cdots,p.
\end{align*}
The asymptotic power of $T_{\max,2}$ has been provided by \cite{Cao2018TwosampleTO}. For readers' convenience, we restate it below.
\begin{proposition}[\cite{Cao2018TwosampleTO}]\label{twomaxal}
    Assume that Assumptions \ref{distriassume}, \ref{assumption_max0},  \ref{addassume1}  hold. Under the local alternative in \eqref{twomaxH1}, if $\sum_{j=1}^p|\delta_{m,j}|=o(p)$ and $\max_{j\in\bS}|\delta_{m,j}|\geq \sqrt{2}+\epsilon$ for some constant $\epsilon>0$, then as $N,p\to\infty$
$$pr(\Phi_{\alpha,2}=1)\to1.$$ 
\end{proposition}
\end{Remark}

\section{Max-Sum type mean test}\label{max-sum-test}
Finally, to concurrently address both dense and sparse-type alternative hypotheses, we introduce a novel combo-type test statistic by combining sum-type and max-type tests. Furthermore, we apply these combo-type test statistics to conduct both one-sample and two-sample tests for high-dimensional compositional data. The following assumption will be imposed:
\begin{Assumpb}\label{assumption_max}
\mbox{There exist constants } $C,\alpha>0$ \mbox{ so that }  $\max_{1\leq i<j\leq p }|\rho_{ij}|=O(p^{-1/2-\alpha}(\log p)^C)$; $p^{-1/2}(\log p)^C \ll \lambda_{\min}(\bR_p)\leq \lambda_{\max}(\bR_p)\ll \sqrt{p}(\log p)^{-1}$\ \mbox{and}\ 
 $\lambda_{\max}(\bR_p)/\lambda_{\min}(\bR_p)=O(p^{\tau})$\ \mbox{for some}\ $\tau \in (0, 1/4)$; $\sup_{p\geq 2}\frac{\tr(\bR_{p}^{i})}{p}<\infty$, for  $i=2,4$; $N=O(p^{\epsilon})$ for $\epsilon \in (\frac{1}{2}, 1]$. (In one-sample case, $n=O(p^{\epsilon})$ for $\epsilon \in (\frac{1}{2}, 1]$.) 
\end{Assumpb}
\begin{Remark}
To derive the asymptotic distribution for our combo-type test statistics, we need the asymptotic independence of sum-type and max-type statistics, which requires weaker dependence in $\bR_p$. In other words, Assumption \ref{assumption_max} in Theorem \ref{one-sample-max-sum} is, in fact, stronger than both Assumption \ref{sumassume} in Theorem \ref{one-sample-sum} and Assumption \ref{assumption_max0} in Theorem \ref{one-sample-max}.
\end{Remark}

\subsection{One sample test}
Let's begin by recalling the definitions of $T_{\tsum}$ and $T_{\max}$ from \eqref{sumst1} and \eqref{maxst1}, respectively:
\begin{align*}
    T_{\tsum}&=\frac{n\bar{\bY}^{T} \bD_{\hat{\boldsymbol{\Gamma}}_p}^{-1}\bar{\bY}-(n-1)p/(n-3)}{\sqrt{2[\tr({\bD}_{\hat{\boldsymbol{\Gamma}}_p}^{-1/2}\hat{\boldsymbol{\Gamma}}_p{\bD}_{\hat{\boldsymbol{\Gamma}}_p}^{-1/2})^2-p^2/(n-1)]}},\nonumber\\
    T_{\max}&=n\max_{1\leq i\leq p}\frac{\bar{Y}_i^2}{\widehat{\gamma}_{ii}}.
\end{align*}
We propose a combo-type test statistic by combining the max-type and the sum-type tests, defined as
\begin{align*}
T_{\com}=\min\{P_S, P_M\},
\end{align*}
where
$
P_{S}=1-\Phi\left\{T_{\tsum}\right\}$ and $P_{M}=1-F(T_{\max}-2\log p+\log\log p).
$
Here $P_{S}$ and $P_{M}$ are the $p$-values for the tests by using statistics $T_{\tsum}$ and $T_{\max}$, separately, and $T_{\com}$ is defined by the smaller one. As demonstrated in the numerical simulations in the later Section \ref{simu}, our combo-type test performs exceptionally well, exhibiting robust performance irrespective of the sparsity of the alternative hypothesis.

In preceding sections, we explore the asymptotic distribution of sum-type and max-type statistics. To derive the asymptotic distribution for $T_{\com}$, we first demonstrate the asymptotic independence between sums and maxima of dependent Gaussian random variables with a non-negatively definite covariance matrix in Lemma \ref{theorem_max_sum}, whose proof is postponed to Section \ref{pftheorem_max_sum}.
\begin{lemma}\label{theorem_max_sum} 
If $\bZ\sim N(c\boldsymbol{1}_p,\boldsymbol{\Psi}_p)$, $\bH=\bG\bZ$ and the following condition hold:
\begin{itemize}
    \item 
\mbox{There exist } $C>0$ \mbox{ and } $\varrho\in (0, 1)$ \mbox{ so that }
$\max_{1\leq i<j \leq p}|\psi_{ij}|\leq \varrho$ \mbox{ and }  $\max_{1\leq i \leq p}\sum_{j=1}^p\psi_{ij}^2\leq (\log p)^C$
 \mbox{for all}\ $p\geq 3$;\ $p^{-1/2}(\log p)^C \ll \lambda_{\min}(\boldsymbol{\Psi}_p)\leq \lambda_{\max}(\boldsymbol{\Psi}_p)\ll \sqrt{p}(\log p)^{-1}$\ \text{and}\ \\
 $\lambda_{\max}(\boldsymbol{\Psi}_p)/\lambda_{\min}(\boldsymbol{\Psi}_p)=O(p^{\tau})$\ \mbox{for some}\ $\tau \in (0, 1/4)$. 
\end{itemize}
Then 
$\frac{H_1^2+\cdots + H_p^2-\tr(\bG\boldsymbol{\Psi}_p\bG)}{\sqrt{2\tr(\bG\boldsymbol{\Psi}_p\bG)^2}}$ and $\max_{1\leq i \leq p}H_i^2-2\log p +\log\log p$ are asymptotically independent.
\end{lemma}

\begin{Remark}
         Our study focuses on Gaussian data characterized by a non-negative definite covariance matrix $\bG\boldsymbol{\Psi}_p\bG$. In contrast, the results presented in \cite{Feng2022AsymptoticIO} apply to Gaussian data with a strictly positive definite covariance matrix $\boldsymbol{\Psi}_p$. It is crucial to highlight that our findings are specifically tailored for high-dimensional compositional data, while the results in \cite{Feng2022AsymptoticIO} are not suitable for such data. 
\end{Remark}
Under Lemma \ref{theorem_max_sum}, we can derive the asymptotically independence between $T_{\tsum}$ and $T_{\max}-2\log p +\log\log p$ as follows, whose proof is postponed to Section \ref{pfone-sample-max-sum}
\begin{theorem}\label{one-sample-max-sum}
Under the null hypothesis in \eqref{hypotest1-one}, if Assumptions \ref{distriassume}, \ref{assumption_max}, \ref{addassume1} hold, then $T_{\tsum}$ and $T_{\max}-2\log p +\log\log p$ are asymptotically independent.
\end{theorem}
According to the definition of $T_{\com}$, its asymptotic distribution can be characterized by the minimum of two standard uniform random variables. We consider the asymptotic power of $T_{\com}$ under local alternatives 
\begin{align}\label{onecomH1}
   H_1: \boldsymbol{\mu}^W = \boldsymbol{\delta}_{\com}+c\boldsymbol{1}_p,
\end{align}
where $\boldsymbol{\delta}_{\com}=\boldsymbol{\delta}_{\tsum}$ or $\boldsymbol{\delta}_{\max}$ and $\boldsymbol{\delta}_{\tsum}$, $\boldsymbol{\delta}_{\max}$ satisfy the conditions in Propositions \ref{onesumal} and \ref{onemaxal}. As indicated in (1) of Theorem \ref{onecoro}, the proposed combo-type test enables a level-$\alpha$ test by rejecting the null hypothesis when $T_{\com} < 1 - \sqrt{1 - \alpha} \approx \frac{\alpha}{2}$ for small $\alpha$. The power function of our combo-type test is
\begin{align}
\beta_C(\boldsymbol{\mu}^W,\alpha)&=P(T_{\com} < 1-\sqrt{1-\alpha})\geq \max\left\{\beta_S(\boldsymbol{\mu}^W,\alpha/2),\beta_M(\boldsymbol{\mu}^W,\alpha/2)\right\}, \label{powerms0}
\end{align}
when $\alpha$ is small, where $\beta_M(\boldsymbol{\mu}^W,\alpha)$ and $\beta_S(\boldsymbol{\mu}^W,\alpha)$ are the power functions of $T_{\max}$ and $T_{\tsum}$  with significant level $\alpha$, respectively. Additionally, we discuss the power functions of $T_{\com}$ in (2) of Theorem \ref{onecoro}. 

\begin{theorem}\label{onecoro}
Assume that Assumptions \ref{distriassume}, \ref{assumption_max}, \ref{addassume1} hold,
    \begin{enumerate}
        \item[(1)] Under the null hypothesis in \eqref{hypotest1-one}, then as $p\to\infty$
\begin{align*}
    \text{$T_{\com}$  converges weakly to a distribution with density  $G(w)=2(1-w)I(0\leq w \leq 1)$. }
\end{align*}

\item[(2)] 
Furthermore, under the local alternative in \eqref{onecomH1}, the power function of our combo-type test
$$\beta_C(\boldsymbol{\mu}^W,\alpha)\to1,\ as \ p\to\infty.$$
       \end{enumerate}
\end{theorem}
\begin{Remark}
    The max-sum type test $T_{\com}$ used for one-sample testing in the context of high-dimensional compositional data, exhibits an identical distribution to the case discussed in \cite{Feng2022AsymptoticIO}.
\end{Remark}

\subsection{Two sample test}
Let's begin by recalling the definitions of $T_{\tsum,2}$ and $T_{\max,2}$ from \eqref{sumst2} and \eqref{maxst2}, respectively:
\begin{align*}
T_{\tsum,2}&=\frac{\frac{n_1n_2}{n_1+n_2}(\bar{\bY}^{(1)}-\bar{\bY}^{(2)})^{T} \bD_{\hat{\bar{\boldsymbol{\Gamma}}}_p}^{-1}(\bar{\bY}^{(1)}-\bar{\bY}^{(2)})-\frac{(n_1+n_2-2)p}{(n_1+n_2-4)}}
{\sqrt{2\big[\tr( {\bD}_{\hat{\bar{\boldsymbol{\Gamma}}}_p}^{-1/2} \hat{\bar{\boldsymbol{\Gamma}}}_p{\bD}_{\hat{\bar{\boldsymbol{\Gamma}}}_p}^{-1/2})^2-\frac{p^2}{(n_1+n_2-2)}\big]c_{p,N}}},\nonumber\\
    T_{\max,2}&=\frac{n_1n_2}{n_1+n_2}\max_{1\leq i\leq p}\frac{(\bar{Y}_{i}^{(1)}-\bar{Y}_{i}^{(2)})^2}{\widehat{\gamma}_{ii}}.
\end{align*}
Under Lemma \ref{theorem_max_sum}, we can derive the asymptotically independence between $T_{\tsum,2}$ and $T_{\max,2}-2\log p +\log\log p$ as follows, whose proof is postponed to Section \ref{pftwo-sample-max-sum}. 
\begin{theorem}\label{two-sample-max-sum}
Under the null hypothesis in \eqref{hypotest1} and Assumptions \ref{distriassume}, \ref{assumption_max}, \ref{addassume1}, then $T_{\tsum,2}$ and $T_{\max,2}-2\log p +\log\log p$ are asymptotically independent.
\end{theorem}
Based on Theorem \ref{two-sample-max-sum}, we propose the following test statistic which utilizes the max-type and sum-type tests. Define
\begin{align}\label{twosamlpecombine}
T_{\com,2}=\min\{P_{M,2}, P_{S,2}\},
\end{align}
where
$P_{M,2}=1-F(T_{\max,2}-2\log p+\log\log p)$ with $F(y)= e^{-\pi^{-1/2}e^{-y/2}}$ and $P_{S,2}=1-\Phi(T_{\tsum,2})$, are the $p$-values of the two tests, respectively. We consider the asymptotic power of $T_{\com,2}$ under local alternatives 
\begin{align}\label{twocomH1}
   H_1: \boldsymbol{\mu}_1^W-\boldsymbol{\mu}_2^W = \boldsymbol{\delta}_{\com,2}+c\boldsymbol{1}_p,
\end{align}
where $\boldsymbol{\delta}_{\com,2}=\boldsymbol{\delta}_{\tsum,2}$ or $\boldsymbol{\delta}_{\max,2}$ and $\boldsymbol{\delta}_{\tsum,2}$, $\boldsymbol{\delta}_{\max,2}$ satisfy the conditions in Propositions \ref{twosumal} and \ref{twomaxal}. Similar to Theorem \ref{onecoro}, we promptly derive the following result through asymptotic independence of $T_{\tsum,2}$ and $T_{\max,2}-2\log p+\log\log p$. As stated in (1) of Theorem \ref{twocoro}, the proposed combo-type test facilitates a level-$\alpha$ test by rejecting the null hypothesis when $T_{\com,2} < 1 - \sqrt{1 - \alpha} \approx \frac{\alpha}{2}$ for small $\alpha$. The power function of our combo-type test is
\begin{align}
\beta_{C,2}(\boldsymbol{\mu}_1^W,\boldsymbol{\mu}_2^W,\alpha)&=P(T_{\com,2} < 1-\sqrt{1-\alpha})\geq \max\left\{\beta_{S,2}(\boldsymbol{\mu}_1^W,\boldsymbol{\mu}_2^W,\alpha/2),\beta_{M,2}(\boldsymbol{\mu}_1^W,\boldsymbol{\mu}_2^W,\alpha/2)\right\},\label{powerms02}
\end{align}
when $\alpha$ is small, where $\beta_{M,2}(\boldsymbol{\mu}_1^W,\boldsymbol{\mu}_2^W,\alpha)$ and $\beta_{S,2}(\boldsymbol{\mu}_1^W,\boldsymbol{\mu}_2^W,\alpha)$ are the power functions of $T_{\max,2}$ and $T_{\tsum,2}$  with significant level $\alpha$, respectively. The power functions are discussed in (2) of Theorem \ref{twocoro}. Again, as demonstrated in the numerical simulations in the later Section \ref{simu}, our combo-type test performs exceptionally well, exhibiting robust performance irrespective of the sparsity of the alternative hypothesis.

\begin{theorem}\label{twocoro}
    \begin{enumerate}
    Assume that Assumptions \ref{distriassume}, \ref{assumption_max}, \ref{addassume1} hold,
        \item[(1)] Under the null hypothesis in \eqref{hypotest1}, then as $p\to\infty$,
\begin{align}
    \text{$T_{\com,2}$ converges weakly to a distribution with density  $G(w)=2(1-w)I(0\leq w \leq 1)$.}
\end{align}
\item[(2)] 
Furthermore, under the local alternative in \eqref{twocomH1}, the power function of our combo-type test
$$\beta_{C,2}(\boldsymbol{\mu}_1^W,\boldsymbol{\mu}_2^W,\alpha)\to1, \ as\ p\to\infty.$$
    \end{enumerate} 
\end{theorem}
\begin{Remark}
     The combo-type test $T_{\com,2}$ employed for two-sample testing in the context of high-dimensional compositional data shares an identical distribution with the case discussed in \cite{Feng2022AsymptoticIO} under the null hypothesis. 
\end{Remark}

\section{Simulation}\label{simu}
In this section, we conduct simulation studies on the testing problems discussed in the previous section. The design of distributions and covariance structures comes from simulation part of \cite{Feng2022AsymptoticIO}. We compare different test statistics sum-type $T_{\tsum}$ proposed in \cite{Srivastava2009ATF}, max-type $T_{\max}$ proposed in \cite{Cao2018TwosampleTO} and combo-type $T_{\com}$ proposed in our paper and validate the advantages of the proposed combo-type tests, where
\begin{align*}
    &T_{\tsum}=\frac{n\bar{\bY}^{T} \bD_{\hat{\boldsymbol{\Gamma}}_p}^{-1}\bar{\bY}-(n-1)p/(n-3)}{\sqrt{2[\tr({\bD}_{\hat{\boldsymbol{\Gamma}}_p}^{-1/2}\hat{\boldsymbol{\Gamma}}_p{\bD}_{\hat{\boldsymbol{\Gamma}}_p}^{-1/2})^2-p^2/(n-1)]}}, \\
    &T_{\tsum,2}=\frac{\frac{n_1n_2}{n_1+n_2}(\bar{\bY}^{(1)}-\bar{\bY}^{(2)})^{T} \bD_{\hat{\bar{\boldsymbol{\Gamma}}}_p}^{-1}(\bar{\bY}^{(1)}-\bar{\bY}^{(2)})-\frac{(n_1+n_2-2)p}{(n_1+n_2-4)}}
{\sqrt{2\big[\tr( {\bD}_{\hat{\bar{\boldsymbol{\Gamma}}}_p}^{-1/2} \hat{\bar{\boldsymbol{\Gamma}}}_p{\bD}_{\hat{\bar{\boldsymbol{\Gamma}}}_p}^{-1/2})^2-\frac{p^2}{(n_1+n_2-2)}\big]c_{p,N}}},\\
    &T_{\max}=n\max_{1\leq i\leq p}\frac{\bar{Y}_i^2}{\widehat{\gamma}_{ii}}, \quad T_{\max,2}=\frac{n_1n_2}{n_1+n_2}\max_{1\leq i\leq p}\frac{(\bar{Y}_{i}^{(1)}-\bar{Y}_{i}^{(2)})^2}{\widehat{\gamma}_{ii}},\\
&T_{\com}=\min\{P_{S}, P_{M}\}, \quad T_{\com,2}=\min\{P_{M,2}, P_{S,2}\}.
\end{align*}
The rejection domains of the three types of test statistics are listed below:
\begin{align*}
    &R_{\tsum}/R_{\tsum,2} = \{T_{\tsum}\ or \ T_{\tsum,2}\geq z_{1-\alpha}\},\\
    &R_{\max}/R_{\max,2}  = \{T_{\max} \ or\ T_{\max,2} \geq -\log\pi-2\log\log(1-\alpha)^{-1} + 2\log p -\log\log p \},\\
    &R_{\com}/R_{\com,2} = \{T_{\com} \ or\ T_{\com,2} < 1-\sqrt{1-\alpha}\},
\end{align*}
where $z_{1-\alpha}$ is the $(1-\alpha)$ quantile of standard normal distribution $N(0,1)$.

We consider the product structure $\log \bW_i = \boldsymbol{\mu}^W + \boldsymbol{\Sigma}^{1/2}\bU_i$, where each component of $\bU_i$ is independently generated from three distributions:

\textbf{A1.} Standard normal distribution $N(0,1)$;

\textbf{A2.} t distribution $t(3)/\sqrt{3}$.

\textbf{A3.} The mixture normal random variable $U$, where $U$ has density function $0.1f_1(x) +0.9f_2(x)$ with $f_1(x)$ and $f_2(x)$ being the densities of $N (0, 9)$ and $N (0, 1)$, respectively.

The following three scenarios of covariance matrices will be considered.

\textbf{B1.} $\boldsymbol{\Sigma} = \left(0.5^{|i-j|}\right)_{1\leq i, j \leq p}$.

\textbf{B2.} $\boldsymbol{\Sigma}=\mathbf{D}^{1 / 2} \bR_p \bD^{1 / 2}$ with $\mathbf{D}=\operatorname{diag}\left(\sigma_{1}^{2}, \cdots, \sigma_{p}^{2}\right)$ and $\mathbf{R}_p=\mathbf{I}_{p}+\boldsymbol{b} \boldsymbol{b}^{T}-\check{\mathbf{B}}$, where $\sigma_{i}^{2}$ are generated independently from Uniform $(1,2), \boldsymbol{b}=\left(b_{1}, \cdots, b_{p}\right)^{T}$ and $\check{\mathbf{B}}=\operatorname{diag}\left(b_{1}^{2}, \cdots, b_{p}^{2}\right)$. The first $\left[p^{0.3}\right]$ entries of $\boldsymbol{b}$ are independently sample from Uniform(0.7, 0.9), and the remaining entries are set to be zero, where $[.]$ denotes taking integer part.

\textbf{B3.} $\boldsymbol{\Sigma}=\boldsymbol{\gamma} \boldsymbol{\gamma}^{T}+\left(\mathbf{I}_{p}-\rho_{\epsilon} \boldsymbol{W}\right)^{-1}\left(\mathbf{I}_{p}-\rho_{\epsilon} \boldsymbol{W}^{T}\right)^{-1}$, where $\boldsymbol{\gamma}=\left(\gamma_{1}, \cdots, \gamma_{\left[p^{\delta_{\gamma}}\right]}, 0,0, \cdots, 0\right)^{T}$. Here $\gamma_{i}$ with $i=1, \cdots,\left[p^{\delta_{\gamma}}\right]$ are generated independently from $\operatorname{Uniform}(0.7,0.9)$. Let $\rho_{\epsilon}=0.5$ and $\delta_{\gamma}=0.3$. Let $\boldsymbol{W}=\left(w_{i_{1} i_{2}}\right)_{1 \leq i_{1}, i_{2} \leq p}$ have a so-called rook form, i.e., all elements of $\boldsymbol{W}$ are zero except that $w_{i_{1}+1, i_{1}}=w_{i_{2}-1, i_{2}}=0.5$ for $i_{1}=1, \cdots, p-2$ and $i_{2}=3, \cdots, p$, and $w_{1,2}=w_{p, p-1}=1$.

Then $\bX^{(k)}$ are generated by transformation 
\begin{equation*}
X_{i j}^{(k)}=W_{i j}^{(k)} / \sum_{\ell=1}^{p} W_{i \ell}^{(k)}, \quad\left(i=1, \ldots, n_{k} ; j=1, \ldots, p ; k=1,2\right).
\end{equation*}
We will work on sample size $n=200$ in one-sample case and $n_1=100, n_2=100$ in two-sample case, and three different dimension with $p=200,400,600$.

\subsection*{One-sample case}
Under the null hypothesis, we set $\boldsymbol{\mu}^W =0$ and $\alpha=0.05$. Table \ref{t1}demonstrate that the empirical sizes of the three test statistics can be well-controlled under various scenarios. 

We also examine the empirical power of each test. Define $\boldsymbol{\mu}^W = (\mu_1^W,\dots, \mu_p^W)^T$.  In order to show the advantages of three statistics with different structures of local alternatives, the following alternative are considered.
For different number of nonzero-mean variables $m = 1,\dots, 20$, we consider $\mu_j^W = \delta$ for $0<j\leq m$ and $\mu_j^W  =0$ for $j>m$. The parameter $\delta $ is chosen as $||\boldsymbol{\mu}^W||^2 = m\delta^2 =0.5$.
\begin{table}[htb]
\centering
\scalebox{0.7}{\begin{tabular}{ccccccccccccc}
\toprule
           &      &         & \textbf{A1} &       &  &       & \textbf{A2}     &       &  &       & \textbf{A3}   &\\
           &      & $p = 200$ & $p=400$  & $p=600$ &  & $p=200$ & $p=400$ & $p=600$ &  & $p=200$ & $p=400$ & $p=600$\\
\midrule
           & $\mathbf{T_{sum}}$ & 0.072   & 0.056  & 0.053 &  & 0.046 & 0.051 & 0.058 &  & 0.043 & 0.054 & 0.055\\
\textbf{B1}& $\mathbf{T_{max}}$ & 0.058   & 0.069  & 0.065 &  & 0.025 & 0.036 & 0.033 &  & 0.047 & 0.038 & 0.049\\
           & $\mathbf{T_{com}}$ & 0.058   & 0.050  & 0.052 &  & 0.041 & 0.046 & 0.047 &  & 0.037 & 0.049 & 0.051\\
\midrule
           & $\mathbf{T_{sum}}$ & 0.072   & 0.061  & 0.047 &  & 0.041 & 0.042 & 0.057 &  & 0.054 & 0.054 & 0.036\\
\textbf{B1}& $\mathbf{T_{max}}$ & 0.055   & 0.068  & 0.058 &  & 0.031 & 0.040 & 0.033 &  & 0.040 & 0.036 & 0.050\\
           & $\mathbf{T_{com}}$ & 0.054   & 0.069  & 0.056 &  & 0.044 & 0.049 & 0.048 &  & 0.046 & 0.045 & 0.039\\
\midrule
           & $\mathbf{T_{sum}}$ & 0.065   & 0.060   & 0.067 &  & 0.055 & 0.054 & 0.050  &  & 0.044 & 0.049 & 0.049\\
\textbf{B3}& $\mathbf{T_{max}}$ & 0.064   & 0.048  & 0.066 &  & 0.033 & 0.043 & 0.047 &  & 0.053 & 0.043 & 0.042\\
           & $\mathbf{T_{com}}$ & 0.053   & 0.056  & 0.063 &  & 0.045 & 0.049 & 0.045 &  & 0.064 & 0.053 & 0.040\\
\hline
\bottomrule
\end{tabular}}
\caption{\textbf{One-sample case:} Sizes of $\mathbf{T_{sum}}$ proposed in \cite{Srivastava2009ATF}, $\mathbf{T_{max}}$ proposed in \cite{Cao2018TwosampleTO}, $\mathbf{T_{com}}$ proposed in our paper in one-sample case with $\alpha = 0.05, n =200$ based on 1000 replications under three different distributions denoted as A1-A3 and three different covariance structure denoted as B1-B3.}\label{t1}
\end{table}
\begin{figure}[htb]
    \centering
    \includegraphics[width=0.75\linewidth]{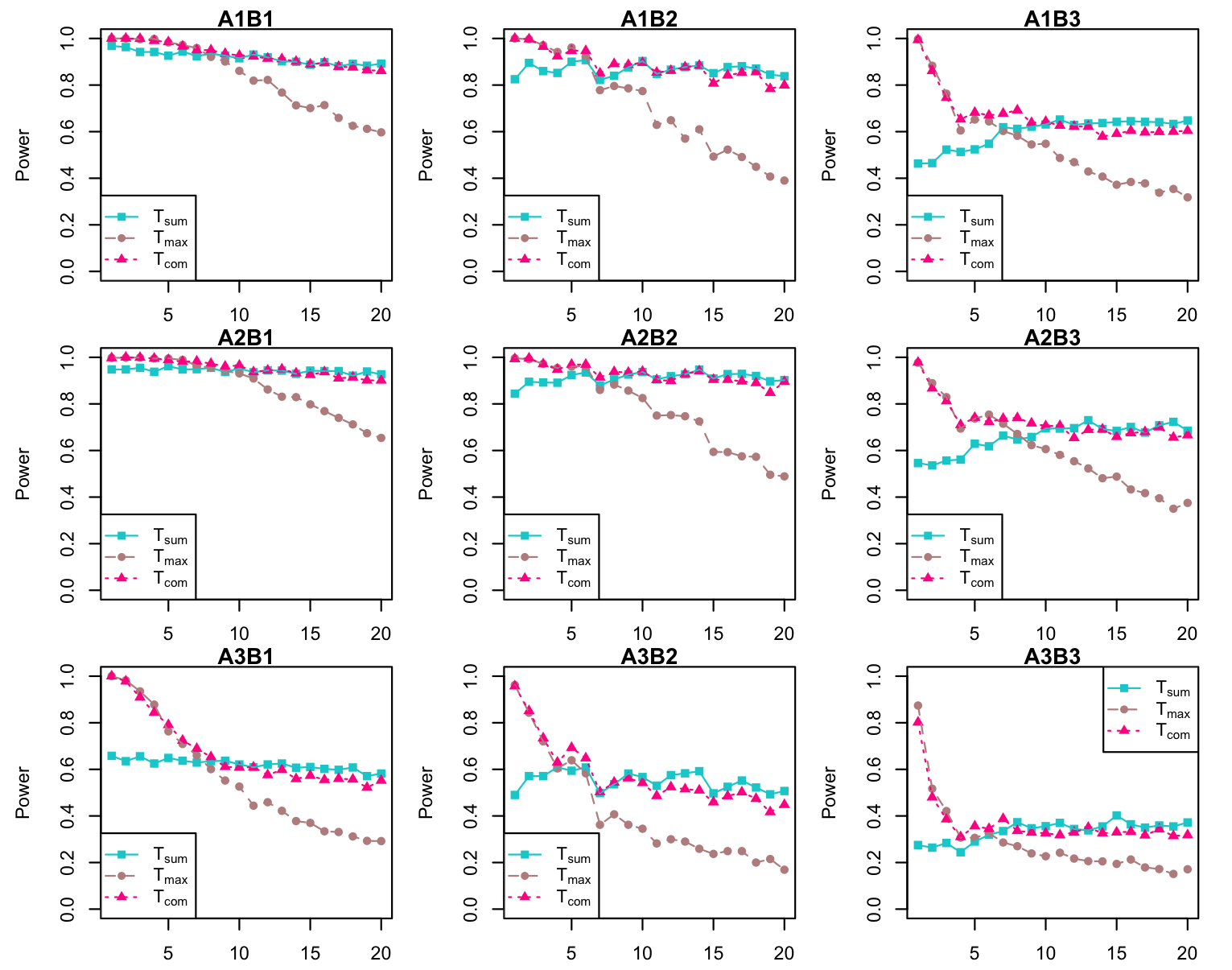}
    \caption{\textbf{One-sample case:} Powers of $\mathbf{T_{sum}}$ proposed in \cite{Srivastava2009ATF}, $\mathbf{T_{max}}$ proposed in \cite{Cao2018TwosampleTO} and $\mathbf{T_{com}}$ proposed in our paper vs. number of variables with nonzero means with $\alpha = 0.05, n=200, p=200$ based on 1000 replications under three different distributions denoted as A1-A3 and three different covariance structure denoted as B1-B3.}
    \label{fig:onesamplepower}
\end{figure}
Figure \ref{fig:onesamplepower} illustrate the empirical power of the three statistics under both sparse and dense structures. As anticipated, the statistic $\mathbf{T_{max}}$ excels in sparse structures but exhibits a notable decrease in efficacy in dense structures. Moreover, we observe that the statistic $\mathbf{T_{sum}}$ outperforms $\mathbf{T_{max}}$ significantly under dense structures and is not sensitive to sparsity. Our proposed statistic $\mathbf{T_{com}}$ demonstrates comparable efficacy to $\mathbf{T_{max}}$ under sparse structures, and its performance closely aligns with $\mathbf{T_{sum}}$ in the case of large $m$. Since, in practice, it is challenging to ascertain whether the true underlying model is sparse or not, our proposed combo-type test, with its robust performance, should be a more favorable choice over the competing approaches.
\newpage
\subsection*{Two-sample case}
\begin{figure}[htb]
    \centering
    \includegraphics[width=0.75\linewidth]{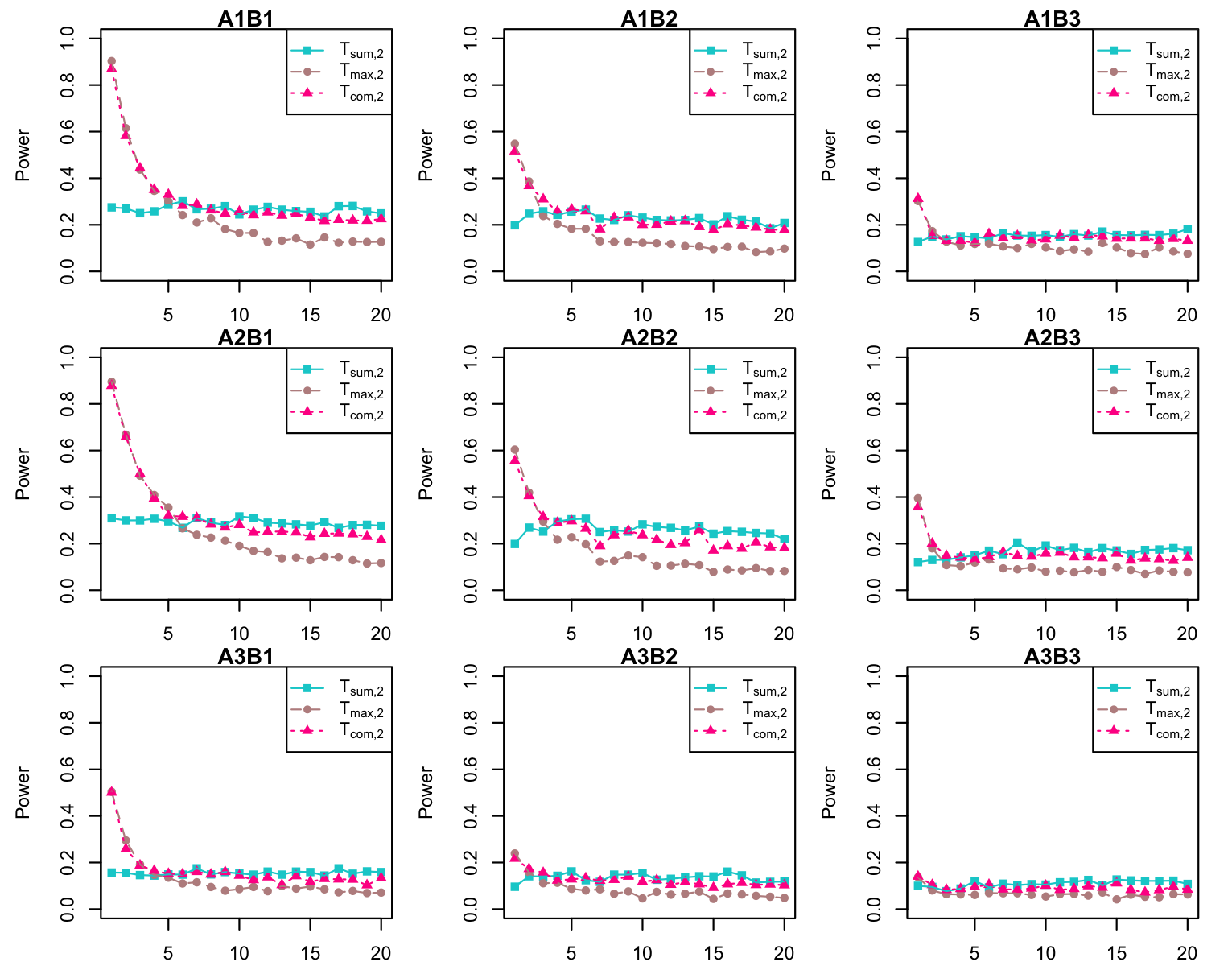}
    \caption{\textbf{Two-sample case:} Powers of $\mathbf{T_{sum,2}}$ proposed in \cite{Srivastava2009ATF}, $\mathbf{T_{max,2}}$ proposed in \cite{Cao2018TwosampleTO} and $\mathbf{T_{com,2}}$ proposed in our paper vs. number of variables with nonzero means with $\alpha = 0.05, n_1=n_2=100, p=200$ based on 1000 replications under three different distributions denoted as A1-A3 and three different covariance structure denoted as B1-B3.}
    \label{fig:twosamplepower}
\end{figure}

Under the same settings as those for the one-sample null and alternative hypotheses, Table \ref{t2} illustrates that the empirical sizes of the three two-sample test statistics can be well-controlled across various scenarios. Additionally, the power depicted in Figure \ref{fig:twosamplepower} suggests that $T_{\com,2}$ remains robust even with data exhibiting different levels of sparsity.
\begin{table}[htb]
\centering
\scalebox{0.7}{\begin{tabular}{ccccccccccccc}
\toprule
           &      &      &\textbf{A1}&       &  &       &\textbf{A2}&       &  &       &\textbf{A3}&\\
           &      &$p = 200$&$p=400$ &$p=600$&  &$p=200$&$p=400$&$p=600$&  &$p=200$&$p=400$&$p=600$\\
\midrule
           & $\mathbf{T_{sum,2}}$ & 0.050    & 0.067  & 0.055 &  & 0.059 & 0.056 & 0.053 &  & 0.065 & 0.069 & 0.054\\
\textbf{B1}& $\mathbf{T_{max,2}}$ & 0.044   & 0.052  & 0.066 &  & 0.036 & 0.040 & 0.029 &  & 0.041 & 0.049 & 0.048\\
           & $\mathbf{T_{com,2}}$ & 0.056   & 0.067  & 0.052 &  & 0.04  & 0.049 & 0.052 &  & 0.044 & 0.063 & 0.053\\
           \midrule
           & $\mathbf{T_{sum,2}}$ & 0.062   & 0.049  & 0.070  &  & 0.055 & 0.058 & 0.048 &  & 0.057 & 0.057 & 0.058\\
\textbf{B2}& $\mathbf{T_{max,2}}$ & 0.042   & 0.054  & 0.061 &  & 0.034 & 0.027 & 0.034 &  & 0.038 & 0.044 & 0.042\\
           & $\mathbf{T_{com,2}}$ & 0.063   & 0.056  & 0.066 &  & 0.044 & 0.037 & 0.050  &  & 0.044 & 0.048 & 0.047\\ 
           \midrule
           & $\mathbf{T_{sum,2}}$ & 0.059   & 0.062  & 0.042 &  & 0.048 & 0.051 & 0.055 &  & 0.060 & 0.059 & 0.049\\
\textbf{B3}& $\mathbf{T_{max,2}}$ & 0.051   & 0.059  & 0.049 &  & 0.036 & 0.024 & 0.046 &  & 0.033 & 0.047 & 0.049\\
           & $\mathbf{T_{com,2}}$ & 0.057   & 0.058  & 0.061 &  & 0.058 & 0.056 & 0.049 &  & 0.054 & 0.050 & 0.059\\
\hline
\bottomrule
\end{tabular}}
\caption{\textbf{Two-sample case:} Sizes of $\mathbf{T_{sum,2}}$ proposed in \cite{Srivastava2009ATF}, $\mathbf{T_{max,2}}$ proposed in \cite{Cao2018TwosampleTO} and $\mathbf{T_{com,2}}$ proposed in our paper in two-sample case with $\alpha = 0.05, n_1=n_2=100$ based on 1000 replications under three different distributions denoted as A1-A3 and three different covariance structure denoted as B1-B3.} \label{t2}
\end{table}

\clearpage
\section{Conclusions}\label{sec5}
High-dimensional compositional data is prevalent in fields such as microbiome research, where the number of variables is comparable to the sample size. In this paper, we focus on mean tests of high-dimensional compositional data for both one-sample and two-sample cases. We propose sum-type and max-sum type tests, which can detect both dense weak and sparse signals. Thorough theoretical and numerical analyses are conducted on these test statistics, and they exhibit satisfactory numerical performance in simulation studies.

\bibliographystyle{elsarticle-num-names} 
\bibliography{references.bib}

\clearpage
\section*{Appendix}
\appendix

Before beginning the proof, we first review the notation as follows. Denote by 
\begin{align*}
    \bQ_i&=\log \bW_i \ (i=1,\ldots,n),\ \
    \bQ_i^{(k)} =\log \bW_i^{(k)} \ (i=1,\ldots,n_k, k=1,2,)\\
    \bar{\bQ}&=\frac{1}{n}\sum_{i=1}^n\bQ_i=(\bar{Q}_1,\ldots,\bar{Q}_p)^T,\ \, \bar{\bQ}^{(k)}=\frac{1}{n_k}\sum_{i=1}^{n_k}\bQ_i^{(k)} \ (k=1,2).
\end{align*}
Recall that 
\begin{align*}
    \bZ&=(Z_1,\ldots,Z_p)^T\sim N(c\boldsymbol{1}_p,\boldsymbol{\Psi}_p),\ \ \bH=\bG\bZ=(H_1,\ldots,H_p)^T\sim N(\boldsymbol{0},\bG\boldsymbol{\Psi}_p\bG)
    \end{align*}
and
\begin{align*}
    \bY_i&=\bG\log \bW_i=\bG\bQ_i,\ \
    \bY_i^{(k)}=\bG\log \bW_i^{(k)}=\bG \bQ_i^{(k)} \ (i=1,\ldots,n_k, k=1,2),\\
    \bar{\bY}&=\frac{1}{n}\sum_{i=1}^n\bY_i=(\bar{Y}_1,\ldots,\bar{Y}_p)^T,\ \, \bar{\bY}^{(k)}=\frac{1}{n_k}\sum_{i=1}^{n_k}\bY_i^{(k)} \ (k=1,2).
\end{align*}
In the one-sample test, $\bV$ is given by
\begin{align*}
\bV
:=\sqrt{n}\bG\bD_{\boldsymbol{\Sigma}_p}^{-1/2}(\bar{\bQ}-c\boldsymbol{1}_p)=(V_1,\cdots,V_p)^T.
\end{align*}
Meanwhile, in the two-sample test, $\bV_2$ is given by
\begin{align*}
    \bV_2:=\Big(\frac{n_1n_2}{n_1+n_2}\Big)^{1/2}\bG\bD_{\boldsymbol{\Sigma}_p}^{-1/2}(\bar{\bQ}^{(1)}-\bar{\bQ}^{(2)}-c\bI_p)
    =(V_{1,2},\cdots,V_{p,2})^T.
\end{align*}
$\boldsymbol{\Sigma}_p=(\sigma_{ij})$ and $\boldsymbol{\Gamma}_p=(\gamma_{ij})$ are the covariance matrices of $\log \bW_i$ and $\bY_i$, respectively. $\bR_p=\bD_{\boldsymbol{\Sigma}_p}^{-1/2}\boldsymbol{\Sigma}_p\bD_{\boldsymbol{\Sigma}_p}^{-1/2}=(\rho_{ij})$ and $\bT_p=\bD_{\boldsymbol{\Gamma}_p}^{-1/2}\boldsymbol{\Gamma}_p\bD_{\boldsymbol{\Gamma}_p}^{-1/2}=(\tau_{ij})$ are the correlation matrices of $\log \bW_i$ and $\bY_i$, respectively. For the reader's convenience, we present the table of contents for Appendix as below:

\startcontents[sections]
\printcontents[sections]{}{1}{\setcounter{tocdepth}{1}}

\section{Proof of Section 3: sum-type test}
\subsection{Proof of Lemma \ref{theorem_sum}}\label{pftheorem_sum}

\begin{proof}
    We aim to prove the asymptotic normality of $H_1^2+\cdots + H_p^2$. It is composed of the following two major steps.\\
    \textbf{Step 1:}
    Since $H_1,\ldots,H_p$ are dependent, we first transfer this problem to independent scenario. We first show that $H_1^2+\cdots + H_p^2$ has the same distribution as that of $$(\xi_1, \cdots, \xi_p)\bG\boldsymbol{\Psi}_p\bG(\xi_1, \cdots, \xi_p)^T,$$ where $\xi_1, \xi_2,\cdots,\xi_{p-1}$ be i.i.d. $N(0,1)$-distributed random variables.\\
    \textbf{Step 2:} We derive the asymptotic distribution of $(\xi_1, \cdots, \xi_p)\bG\boldsymbol{\Psi}_p\bG(\xi_1, \cdots, \xi_p)^T$.

\textbf{Step 1:} Since $\bZ\sim N(c\boldsymbol{1}_p,\boldsymbol{\Psi}_p)$, $\bH=\bG\bZ\sim N(\boldsymbol{0},\bG\boldsymbol{\Psi}_p\bG)$. Let $(\bG\boldsymbol{\Psi}_p\bG)^{1/2}$ be a non-negative definite matrix such that $(\bG\boldsymbol{\Psi}_p\bG)^{1/2}\cdot (\bG\boldsymbol{\Psi}_p\bG)^{1/2}=\bG\boldsymbol{\Psi}_p\bG$. Let $\xi_1, \xi_2,\cdots,\xi_p$ be i.i.d. $N(0,1)$-distributed random variables. Then $(H_1, \cdots, H_p)^T$ and  $(\bG\boldsymbol{\Psi}_p\bG)^{1/2}(\xi_1, \cdots, \xi_p)^T$ have the same distribution. As a consequence, $H_1^2+\cdots + H_p^2$ has the same distribution as that of
\begin{align} \label{An_Ru}
(\xi_1, \cdots, \xi_p)(\bG\boldsymbol{\Psi}_p\bG)^{1/2}\cdot(\bG\boldsymbol{\Psi}_p\bG)^{1/2}(\xi_1, \cdots, \xi_p)^T
=(\xi_1, \cdots, \xi_p)\bG\boldsymbol{\Psi}_p\bG(\xi_1, \cdots, \xi_p)^T.
\end{align}

\textbf{Step 2.} Note that $\bH=\bG\bZ=(\bg_1,\dots,\bg_p)^T\bZ=(H_1,\dots,H_p)^T$, where $H_i=\bg_i^T\bZ$. Thus,
\begin{align*}
    \Var(H_i)=\bg_i^T\Var(\bZ)\bg_i=\bg_i^T\boldsymbol{\Psi}_p\bg_i.
\end{align*}
Hence, 
\begin{align}\label{trGamma}
    E(H_1^2+\cdots + H_p^2)&=\sum_{i=1}^p\Var(H_i)=\sum_{i=1}^p\tr(\boldsymbol{\Psi}_p\bg_i\bg_i^T)=\tr(\boldsymbol{\Psi}_p\bG)=\tr(\bG\boldsymbol{\Psi}_p\bG).
\end{align}
Let $\lambda_{p, 1}, \lambda_{p,2}, \cdots, \lambda_{p,p}$ be the eigenvalues of $\bG\boldsymbol{\Psi}_p\bG$ and $\boldsymbol{O}$ be a $p\times p$ orthogonal matrix such that $\bG\boldsymbol{\Psi}_p\bG=\boldsymbol{O}^T\,\mbox{diag}(\lambda_{p,1}, \cdots, \lambda_{p,p-1},\lambda_{p,p})\boldsymbol{O}$, where $\lambda_{p,p}=0$.
By the orthogonal invariance of normal distributions, $\boldsymbol{O}(\xi_1, \cdots, \xi_p)^T$ and $(\xi_1, \cdots, \xi_p)^T$ have the same distribution. By \eqref{An_Ru}, $H_1^2+\cdots + H_p^2$ is equal to
\begin{align*}
    \big[\boldsymbol{O}(\xi_1, \cdots, \xi_p)^T\big]^T\,\mbox{diag}(\lambda_{p, 1}, \cdots, \lambda_{p,p-1},0)\,\big[\boldsymbol{O}(\xi_1, \cdots, \xi_p)^T\big],
\end{align*}
and hence has the same distribution as that of $\lambda_{p,1}\xi_1^2+\cdots +\lambda_{p, p-1}\xi_{p-1}^2$. 
Thus,
\begin{align}\label{Hvar}
\mbox{Var}(H_1^2+\cdots + H_p^2)=\lambda_{p,1}^2\mbox{Var}(\xi_1^2)+\cdots +\lambda_{p, p-1}^2\mbox{Var}(\xi_{p-1}^2)
=2\lambda_{p,1}^2+\cdots +2\lambda_{p, p-1}^2
 =  2\cdot \mbox{tr}(\bG\boldsymbol{\Psi}_p\bG)^2.
\end{align}
Next, we decompose $\bG=\bP^{T}diag\{1,\dots,1,0\}\bP$, where $\bP$ is a $p\times p$ orthogonal matrix. We have
$$ \bG\boldsymbol{\Psi}_p\bG
        =\bP^{T}diag\{1,\dots,1,0\}\bP\boldsymbol{\Psi}_p \bP^{T}diag\{1,\dots,1,0\}\bP.  $$
Hence by eigenvalue interlacing theorem, we have 
\begin{align}\label{eigval}
    \lambda_{p,1}(\boldsymbol{\Psi}_p)\geq \lambda_{p,1}(\boldsymbol{\bG\boldsymbol{\Psi}_p\bG}_p)\geq \dots\geq\lambda_{p,p-1}(\boldsymbol{\Psi}_p)\geq \lambda_{p,p-1}(\boldsymbol{\bG\boldsymbol{\Psi}_p\bG}_p)\geq\lambda_{p,p}(\boldsymbol{\Psi}_p)>0.
\end{align}
Moreover, $\lambda_{p,p}(\boldsymbol{\bG\boldsymbol{\Psi}_p\bG})=0$. Thus, by \eqref{eigval}
\begin{align}\label{trpsi}
    \lim_{p\to\infty}\frac{\tr(\bG\boldsymbol{\Psi}_p\bG)^2}{p}\leq \lim_{p\to\infty}\frac{\tr(\boldsymbol{\Psi}_p)^2}{p}<\infty.
\end{align}
Therefore, by \eqref{trGamma}, \eqref{Hvar}, \eqref{trpsi} and Theorem 2.1 in \cite{wang2014jointclt}, we have as $p\to\infty$,
\begin{align}\label{disequ}
    \frac{\bH^T\bH-\tr(\bG\boldsymbol{\Psi}_p\bG)}{\sqrt{2\tr(\bG\boldsymbol{\Psi}_p\bG)^2}}&\stackrel{d}{=}\frac{(\xi_1, \cdots, \xi_p)\bG\boldsymbol{\Psi}_p\bG(\xi_1, \cdots, \xi_p)^T-\tr(\bG\boldsymbol{\Psi}_p\bG)}{\sqrt{2\tr(\bG\boldsymbol{\Psi}_p\bG)^2}}\overset{d}{\to} N(0,1).
\end{align}
\end{proof}

\subsection{Proof of Theorem \ref{one-sample-sum}: sum-type test for one-sample case}\label{pfone-sample-sum}

\begin{proof}
We divide the proof of Theorem \ref{one-sample-sum} into the following steps.\\
    \textbf{Step 1:} We first show that
    \begin{align*}
    \tilde{\tilde{T}}_{\tsum}=\tilde{T}_{\tsum}+o_p(1),
\end{align*}
where $\tilde{\tilde{T}}_{\tsum}=\frac{n\boldsymbol{\bar{Y}}^T\bD^{-1}\boldsymbol{\bar{Y}}-p}{\sqrt{2tr\bT_p^2}}$ and $\tilde{T}_{\tsum}=\frac{\bV^T\bV-tr(\bG\bR_p\bG)}{\sqrt{2tr(\bG\bR_p\bG)^2}}$.\\
 \textbf{Step 2:} Secondly, similar to the proof of Theorem 2.1  in \cite{Srivastava2008ATF}, we have
\begin{align}\label{onett}
T_{\tsum}=\tilde{\tilde{T}}_{\tsum}+o_p(1).
\end{align}
Specifically, by Assumptions \ref{sumassume}-\ref{addassume1} and \eqref{gsbound}, and akin to the proof of \eqref{boundT1}, \eqref{boundT2}, we have for $i=2,4$ 
\begin{align}\label{GtrT}
    \lim_{p\to \infty}\frac{\tr(\bT_p^i)}{p}\leq K \lim_{p\to \infty}\frac{\tr(\bR_p^i)}{p}<\infty,
\end{align}
where $K>0$ is a constant. By \eqref{GtrT} and $n=O(p^{\epsilon})$ for $\epsilon \in (\frac{1}{2}, 1]$ in Assumption \eqref{sumassume}, $[\tr({\bD}_{\hat{\boldsymbol{\Gamma}}_p}^{-1/2}\hat{\boldsymbol{\Gamma}}_p{\bD}_{\hat{\boldsymbol{\Gamma}}_p}^{-1/2})^2-p^2/(n-1)]/p$ is a consistent estimator of $\tr(\bT_p^2)/p$ as $p\to \infty$. Hence, similar to the proof of Theorem 2.1  in \cite{Srivastava2008ATF}, \eqref{onett} is derived.\\
 \textbf{Step 3:} By Steps 1--2, we have
\begin{align}\label{onesum3}
T_{\tsum}=\tilde{T}_{\tsum}+o_p(1).
\end{align}
 Moreover, by Lemma \ref{theorem_sum} we have $\tilde{T}_{\tsum}=\frac{\bV^T\bV-\tr(\bG\bR_p\bG)}{\sqrt{2\tr(\bG\bR_p\bG)^2}}\overset{d}{\to} N(0,1)$ as $p\to \infty$. Hence, $T_{\tsum}\overset{d}{\to} N(0,1)$ as $p\to \infty$.

\textbf{Therefore, the remaining part of the proof is dedicated to proving Step 1.}

\textbf{Step 1:}
Since $\sqrt{n}(\bar{\bQ}-c\boldsymbol{1}_p)\sim N(\boldsymbol{0},\boldsymbol{\Sigma}_p)$ under $H_0$, $\sqrt{n}\bar{\bY}=\sqrt{n}\bG\bar{\bQ}\sim N(\boldsymbol{0},\boldsymbol{\Gamma}_p)$ under $H_0$. Thus, $\sqrt{n}\bD_{\boldsymbol{\Sigma}_p}^{-1/2}(\bar{\bQ}-c\boldsymbol{1}_p)\sim N(\boldsymbol{0},\bR_p)$, $\sqrt{n}\bD_{\boldsymbol{\Gamma}_p}^{-1/2}\bar{\bY}=\sqrt{n}\bD_{\boldsymbol{\Gamma}_p}^{-1/2}\bG\bar{\bQ}\sim N(\boldsymbol{0},\bT_p)$, $\bV:=\sqrt{n}\bG\bD_{\boldsymbol{\Sigma}_p}^{-1/2}(\bar{\bQ}-c\boldsymbol{1}_p):=(V_1,\cdots,V_p)\sim N(\boldsymbol{0},\bG\bR_p\bG)$ under $H_0$. Thus, 
\begin{align}\label{diffsum0}
    \tilde{\tilde{T}}_{\tsum}-\tilde{T}_{\tsum} = \frac{n\boldsymbol{\bar{Y}}^T\bD_{\boldsymbol{\Gamma}_p}^{-1}\boldsymbol{\bar{Y}}-p}{\sqrt{2\tr\bT_p^2}}-\frac{\bV^T\bV-\tr(\bG\bR_p\bG)}{\sqrt{2\tr(\bG\bR_p\bG)^2}}
    =I_1+I_2,
\end{align}
where
\begin{align*}
    I_1&:=\frac{\sqrt{2\tr(\bG\bR_p\bG)^2}[(n\boldsymbol{\bar{Y}}^T\bD_{\boldsymbol{\Gamma}_p}^{-1}\boldsymbol{\bar{Y}}-p)-(\bV^T\bV-\tr(\bG\bR_p\bG))]}{\sqrt{2\tr\bT_p^2}\sqrt{2\tr(\bG\bR_p\bG)^2}},\nonumber\\
    I_2&:=\frac{(\sqrt{2\tr(\bG\bR_p\bG)^2}-\sqrt{2\tr\bT_p^2})(\bV^T\bV-\tr(\bG\bR_p\bG))}{\sqrt{2\tr\bT_p^2}\sqrt{2\tr(\bG\bR_p\bG)^2}}.\nonumber
\end{align*}
We need to show
\begin{align}\label{J1J2limit0}
    I_1\to 0 \ \text{and}  \ I_2\to 0  \ \text{in probability  as} \ p\to\infty.
\end{align}
To do so, it suffices to prove that both their means and variances converge to 0. Write $\sqrt{n}(\bar{\bQ}-c\boldsymbol{1}_p)= \boldsymbol{\Sigma}_p^{1/2} \boldsymbol{\xi}$, where $\boldsymbol{\xi}\sim\boldsymbol{N}(0,\bI_p)$, $\boldsymbol{\Sigma}_p^{1/2}$ is a positive definite matrix such that $\boldsymbol{\Sigma}_p^{1/2}\boldsymbol{\Sigma}_p^{1/2}=\boldsymbol{\Sigma}_p$. Note that, 
\begin{align*}
    n\boldsymbol{\bar{Y}}^T\bD_{\boldsymbol{\Gamma}_p}^{-1}\boldsymbol{\bar{Y}}= \boldsymbol{\xi}^T \boldsymbol{\Sigma}_p^{1/2}\bG\bD_{\boldsymbol{\Gamma}_p}^{-1}\bG\boldsymbol{\Sigma}_p^{1/2}\boldsymbol{\xi}, \ \
    \bV^T\bV = \boldsymbol{\xi}^T \boldsymbol{\Sigma}_p^{1/2}\bD_{\boldsymbol{\Sigma}_p}^{-1/2}\bG\bD_{\boldsymbol{\Sigma}_p}^{-1/2}\boldsymbol{\Sigma}_p^{1/2}\boldsymbol{\xi}.
\end{align*}
Thus, 
\begin{align}\label{Enyy0}
    \Expe(\bV^T\bV) = \Expe(\boldsymbol{\xi}^T \boldsymbol{\Sigma}_p^{1/2}\bD_{\boldsymbol{\Sigma}_p}^{-1/2}\bG\bD_{\boldsymbol{\Sigma}_p}^{-1/2}\boldsymbol{\Sigma}_p^{1/2}\boldsymbol{\xi})
    =\tr( \boldsymbol{\Sigma}_p^{1/2}\bD_{\boldsymbol{\Sigma}_p}^{-1/2}\bG\bD_{\boldsymbol{\Sigma}_p}^{-1/2}\boldsymbol{\Sigma}_p^{1/2})
    =\tr( \bG\bR_p\bG),
\end{align}
and
\begin{align}\label{Enydynyy0}
\Expe[n\boldsymbol{\bar{Y}}^T\bD_{\boldsymbol{\Gamma}_p}^{-1}\boldsymbol{\bar{Y}}-\bV^T\bV]=\Expe[\boldsymbol{\Sigma}_p^{1/2}\bG\bD_{\boldsymbol{\Gamma}_p}^{-1}\bG\boldsymbol{\Sigma}_p^{1/2}]-\tr( \bG\bR_p\bG)
=\tr[\bT_p] - \tr[\bG\bR_p\bG].
\end{align}
Thus, by \eqref{Enyy0}-\eqref{Enydynyy0} and $\tr[\bT_p] =p$, 
\begin{align}\label{EI1}
    \Expe I_1 = \Expe\frac{\sqrt{2\tr(\bG\bR_p\bG)^2}[(n\boldsymbol{\bar{Y}}^T\bD_{\boldsymbol{\Gamma}_p}^{-1}\boldsymbol{\bar{Y}}-p)-(\bV^T\bV-\tr(\bG\bR_p\bG))]}{\sqrt{2\tr\bT_p^2}\sqrt{2\tr(\bG\bR_p\bG)^2}}=0,
\end{align}
and
\begin{align}\label{EI2}
    \Expe I_2 = \Expe \frac{(\sqrt{2\tr(\bG\bR_p\bG)^2}-\sqrt{2\tr\bT_p^2})(\bV^T\bV-\tr(\bG\bR_p\bG))}{\sqrt{2\tr\bT_p^2}\sqrt{2\tr(\bG\bR_p\bG)^2}} = 0.
\end{align}
Note that
\begin{align}\label{varnydynyy0}
\Var[n\boldsymbol{\bar{Y}}^T\bD_{\boldsymbol{\Gamma}_p}^{-1}\boldsymbol{\bar{Y}}-\bV^T\bV]
&=\Var[\boldsymbol{\xi}^T \boldsymbol{\Sigma}_p^{1/2}(\bG\bD_{\boldsymbol{\Gamma}_p}^{-1}\bG-\bD_{\boldsymbol{\Sigma}_p}^{-1/2}\bG\bD_{\boldsymbol{\Sigma}_p}^{-1/2})\boldsymbol{\Sigma}_p^{1/2}\boldsymbol{\xi}]\nonumber\\
&=2\tr[(\bG\bD_{\boldsymbol{\Gamma}_p}^{-1}\bG-\bD_{\boldsymbol{\Sigma}_p}^{-1/2}\bG\bD_{\boldsymbol{\Sigma}_p}^{-1/2})\boldsymbol{\Sigma}_p]^2\nonumber\\
&=\sum_{i,j=1}^p(a_{ij}-b_{ij})^2\sigma_{ij}^2\nonumber\\
&\leq \max_{1\leq i,j\leq p}|a_{ij}-b_{ij}|^2\tr(\boldsymbol{\Sigma}_p^2),
\end{align}
where $\bG\bD_{\boldsymbol{\Gamma}_p}^{-1}\bG:=(a_{ij})$, $\bD_{\boldsymbol{\Sigma}_p}^{-1/2}\bG\bD_{\boldsymbol{\Sigma}_p}^{-1/2}:=(b_{ij})$. Moreover,
\begin{align*}
    a_{ii}= \gamma_{ii}^{-1}-\frac{1}{p}\gamma_{ii}^{-1}-\frac{1}{p}\gamma_{ii}^{-1}+\frac{1}{p^2}\sum_{i=1}^p\gamma_{ii}^{-1},\ \
    a_{ij}=-\frac{1}{p}\gamma_{ii}^{-1}-\frac{1}{p}\gamma_{jj}^{-1}+\frac{1}{p^2}\sum_{i=1}^p\gamma_{ii}^{-1}, \ \ i\neq j,
\end{align*}
and 
\begin{align*}
    b_{ii}=\sigma_{ii}^{-1}-\frac{1}{p}\sigma_{ii}^{-1}, \ \ 
   b_{ij}=-\frac{1}{p}\sigma_{ii}^{-1/2}\sigma_{jj}^{-1/2}, \ \ i\neq j,
\end{align*}
thus
\begin{align}\label{abdiffbound}
    \max|a_{ij}-b_{ij}|=O(p^{-1/2-\alpha}(\log p)^{C/2}).
\end{align}
Furthermore, 
\begin{align}\label{trineqTG}
    \tr\bT_p^2 = \tr(\bD_{\boldsymbol{\Gamma}_p}^{-1/2}\boldsymbol{\Gamma}_p\bD_{\boldsymbol{\Gamma}_p}^{-1/2})^2
    \geq \min_{i,j}\gamma_{ii}^{-1}\gamma_{jj}^{-1}\sum_{i,j=1}^p\gamma_{ij}^2
    \geq K_1 \tr\boldsymbol{\Gamma}_p^2,
\end{align}
where $K_1>0$ is a constant and  the last inequality is derived from Assumption \ref{addassume1}, which implies that $1/\kappa_2<\gamma_{ii}<\kappa_2$ for some $\kappa_2>0$. Simlarly,
\begin{align}\label{trineqTG2}
    \tr\bT_p^2 
    \leq  K_2 \tr\boldsymbol{\Gamma}_p^2,
\end{align}
where $K_2>0$ is a constant. Therefore,

\begin{align}\label{boundT1}
    K_1 \tr \boldsymbol{\Gamma}_p^2 \leq \tr \bT_p^2 \leq K_2 \tr \boldsymbol{\Gamma}_p^2.
\end{align}
Note that, by eigenvalue interlacing theorem
\begin{align}\label{ineqGS2}
\frac{\tr(\bG\boldsymbol{\Sigma}_p\bG)^2}{\tr(\boldsymbol{\Sigma}_p^2)}&=\frac{\sum_{i=1}^p\lambda_{p,i}^2(\bG\boldsymbol{\Sigma}_p\bG)}{\sum_{i=1}^p\lambda_{p,i}^2(\boldsymbol{\Sigma}_p)}\leq 1,
\end{align}
and 
\begin{align}\label{ineqGS3}
    \frac{\tr(\bG\boldsymbol{\Sigma}_p\bG)^2}{\tr(\boldsymbol{\Sigma}_p^2)}&=\frac{\sum_{i=1}^p\lambda_{p,i}^2(\bG\boldsymbol{\Sigma}_p\bG)}{\sum_{i=1}^p\lambda_{p,i}^2(\boldsymbol{\Sigma}_p)}\geq 1- \frac{\lambda_{\max}^2(\boldsymbol{\Sigma}_p)}{\sum_{i=1}^p\lambda_{p,i}^2(\boldsymbol{\Sigma}_p)}
    > 0.
\end{align}
Hence by \eqref{ineqGS2}-\eqref{ineqGS3} and \eqref{boundT1},
\begin{align}\label{boundT2}
    K_3 \tr \boldsymbol{\Sigma}_p^2 \leq \tr \bT_p^2 \leq K_4 \tr \boldsymbol{\Sigma}_p^2.
\end{align}
where $K_3, K_4>0$ are constants.
By \eqref{varnydynyy0}, \eqref{abdiffbound} and \eqref{boundT2}, we get as $p\to\infty$
\begin{align}\label{VI1}
    \Var I_1 =\frac{\Var[n\boldsymbol{\bar{Y}}^T\bD_{\boldsymbol{\Gamma}_p}^{-1}\boldsymbol{\bar{Y}}-\bV^T\bV]}{2tr\bT_p^2}
    \leq C_2\frac{\Var[n\boldsymbol{\bar{Y}}^T\bD_{\boldsymbol{\Gamma}_p}^{-1}\boldsymbol{\bar{Y}}-\bV^T\bV]}{2tr\boldsymbol{\Sigma}_p^2}
    \leq C_2\max_{1\leq i,j\leq p}|a_{ij}-b_{ij}|^2\to 0.
\end{align}
Note that, 
\begin{align}\label{varnyy0}
    \Var(\bV^T\bV) = \Var(\boldsymbol{\xi}^T \boldsymbol{\Sigma}_p^{1/2}\bD_{\boldsymbol{\Sigma}_p}^{-1/2}\bG\bD_{\boldsymbol{\Sigma}_p}^{-1/2}\boldsymbol{\Sigma}_p^{1/2}\boldsymbol{\xi})
    =2\tr( \bG\bR_p\bG)^2.
\end{align}
By \eqref{varnyy0}, we get 
\begin{align}\label{VI2}
   \Var(I_2)
     \leq \Big(\frac{\tr(\bG\bR_p\bG)^2-tr\bT_p^2}{\tr\bT_p^2}\Big)^2.
\end{align}
Moreover, 
\begin{align}\label{grg2t2}
|\tr(\bG\bR_p\bG)^2-\tr\bT_p^2| &\leq|\tr(\bG\bR_p\bG)^2-\tr\bR_p^2|+|\tr\bR_p^2-\tr(\bD_{\boldsymbol{\Sigma}_p}^{-1/2}\bG\boldsymbol{\Sigma}_p\bG\bD_{\boldsymbol{\Sigma}_p}^{-1/2})^2| \nonumber\\
    &\ \ +|\tr(\bD_{\boldsymbol{\Sigma}_p}^{-1/2}\boldsymbol{\Gamma}_p\bD_{\boldsymbol{\Sigma}_p}^{-1/2})^2-\tr\bT_p^2|\nonumber\\
    &:= S_1+S_2+S_3.
\end{align}
Thus, to estimate the bound of \eqref{grg2t2}, we need to estimate the bounds of $S_1$, $S_2$ and $S_3$. 
Furthermore, by $\max_{1\leq i\neq j\leq p}|\rho_{ij}|=O(p^{-\frac{1}{2}-\alpha}(\log p)^{\frac{C}{2}})$ in Assumption \ref{sumassume},
\begin{align}\label{Rbound}
    ||\bR_p||_1 = \max_{1\leq j\leq p}\sum_{i=1}^p|\rho_{ij}|=1+\max_{1\leq j\leq p}\sum_{i\neq j}^p|\rho_{ij}|\leq 1+p\max_{1\leq i\neq j\leq p}|\rho_{ij}| =O(p^{\frac{1}{2}-\alpha}(\log p)^{\frac{C}{2}}).
\end{align}
Hence,
\begin{align}\label{sigma1norm}
    ||\boldsymbol{\Sigma}_p||_1=\max_{1\leq j\leq p}\sum_{i=1}^p|\sigma_{ij}|=\max_{1\leq j\leq p}\sum_{i=1}^p|\rho_{ij}\sqrt{\sigma_{ii}\sigma_{jj}}|
    \leq \kappa_1 ||\bR_p||_1=O(p^{\frac{1}{2}-\alpha}(\log p)^{\frac{C}{2}}).
\end{align}
Since $1/\kappa_1\leq\sigma_{ii}\leq \kappa_1$ for some constant $\kappa_1>0$,
\begin{align}\label{sigmaidot}
    |\sigma_{i\cdot}|=|\frac{1}{p}\sum_{j=1}^p\sigma_{ij}|
    \leq\frac{1}{p}\sum_{j=1}^p|\rho_{ij}\sqrt{\sigma_{ii}\sigma_{jj}}|
    \leq \kappa_1\frac{1}{p}\sum_{j=1}^p|\rho_{ij}|=\frac{1}{p}\kappa_1||\bR_p||_1=O(p^{-1/2-\alpha}(\log p)^{C/2}),
\end{align}
and similarly $|\sigma_{j\cdot}|=O(p^{-1/2-\alpha}(\log p)^{C/2})$, and 
\begin{align}\label{sigmadotdot}
    |\sigma_{\cdot\cdot}|
    =|\frac{1}{p^2}\sum_{i,j=1}^p\sigma_{ij}|
    &\leq \frac{1}{p}\sum_{i=1}^p \frac{1}{p}\max_{1\leq i\leq p}\sum_{j=1}^p|\sigma_{ij}|
    \leq \kappa_1 \frac{1}{p}\sum_{i=1}^p \frac{1}{p}\max_{1\leq i\leq p}\sum_{j=1}^p |\rho_{ij}|\nonumber \\
    &\leq \kappa_1 \frac{1}{p}||\bR_p||_1=O(p^{-1/2-\alpha}(\log p)^{C/2}).
\end{align}
Since $\gamma_{ij}=\sigma_{ij}-\sigma_{i\cdot}-\sigma_{j\cdot}+\sigma_{\cdot\cdot}$,  by \eqref{sigmaidot}-\eqref{sigmadotdot}, we get
\begin{align}\label{gsbound}
    \max_{1\leq i,j\leq p}|\gamma_{ij}-\sigma_{ij}|=O(p^{-1/2-\alpha}(\log p)^{C/2}).
\end{align}
Similarly,
\begin{align}\label{pwbound}
    \max_{1\leq i,j\leq p}|\rho_{ij}-w_{ij}|=O(p^{-1/2-\alpha}(\log p)^{C/2}),
\end{align}
where $\bG\bR_p\bG=(w_{ij})$. By \eqref{Rbound}, \eqref{pwbound} and $\tr \bT_p^2\geq p$,
\begin{align}\label{S1lim}
    \frac{S_1}{\tr \bT_p^2}&=\frac{|\tr(\bG\bR_p\bG)^2-\tr\bR_p^2|}{\tr \bT_p^2}
    =\frac{\sum_{i,j=1}^p(\rho_{ij}-w_{ij})^2+2\sum_{i,j=1}^pw_{ij}(\rho_{ij}-w_{ij})}{\tr \bT_p^2}\nonumber\\
     &\leq O(p^{-2\alpha}(\log p)^C)+O(p^{-1/2-\alpha}(\log p)^{C/2})||\bG\bR_p\bG||_1\nonumber\\
      &\leq O(p^{-2\alpha}(\log p)^C)+O(p^{-1/2-\alpha}(\log p)^{C/2})||\bG||_1||\bR_p||_1||\bG||_1\nonumber\\
     &\leq O(p^{-2\alpha}(\log p)^C) \to 0,\ \
    \text{as $p\to\infty$}.
\end{align}
Therefore, with $1/\kappa_1\leq\sigma_{ii}\leq \kappa_1$ for some constant $\kappa_1>0$ and $1/\kappa_2\leq\gamma_{ii}\leq \kappa_2$ for some constant $\kappa_2>0$, \eqref{sigma1norm} and \eqref{gsbound}, we have
\begin{align}\label{r2dgd2}
   S_2&=|\tr\bR_p^2-\tr(\bD_{\boldsymbol{\Sigma}_p}^{-1/2}\bG\boldsymbol{\Sigma}_p\bG\bD_{\boldsymbol{\Sigma}_p}^{-1/2})^2|\nonumber\\
    & \leq \sum_{i=1}^p|\sigma_{ii}^{-2}||\gamma_{ii}-\sigma_{ii}||\gamma_{ii}+\sigma_{ii}|+|\sum_{i\neq j}^p (\sigma_{ii}\sigma_{jj})^{-1}(\gamma_{ij}-\sigma_{ij})^2+2\sum_{i\neq j}^p (\sigma_{ii}\sigma_{jj})^{-1}\sigma_{ij}(\gamma_{ij}-\sigma_{ij})|\nonumber\\
    &\leq Kp^2\max_{i,j}|\gamma_{ij}-\sigma_{ij}|^2+Kp\max_{i,j}|\gamma_{ij}-\sigma_{ij}|\max_{1\leq j\leq p}\sum_{i=1}^p|\sigma_{ij}|\nonumber\\
      &\leq O(p^{1-2\alpha}(\log p)^{C}).
    \end{align}
By \eqref{r2dgd2} and $\tr\bT_p^2\geq p$, we have as $p\to\infty$
\begin{align}\label{S2lim}
  \frac{S_2}{\tr\bT_p^2} \leq K\frac{|tr\bR_p^2-\tr(\bD_{\boldsymbol{\Sigma}_p}^{-1/2}\bG\boldsymbol{\Sigma}_p\bG\bD_{\boldsymbol{\Sigma}_p}^{-1/2})^2|}{\tr\bT_p^2}\leq K\frac{p^{1-2\alpha}(\log p)^{C}}{p}\to 0.
\end{align}
For $S_3$, by $1/\kappa_1\leq\sigma_{ii}\leq \kappa_1$ for some constant $\kappa_1>0$ and $1/\kappa_2\leq\gamma_{ii}\leq \kappa_2$ for some constant $\kappa_2>0$,
\begin{align}\label{S3bound}
S_3&=|\tr(\bD_{\boldsymbol{\Sigma}_p}^{-1/2}\boldsymbol{\Gamma}_p\bD_{\boldsymbol{\Sigma}_p}^{-1/2})^2-\tr\bT_p^2|
=|\sum_{i,j=1}^p(\sigma_{ii}^{-1}\sigma_{jj}^{-1}\gamma_{ij}^2-\gamma_{ii}^{-1}\gamma_{jj}^{-1}\gamma_{ij}^2)|\nonumber\\
&\leq \max_{i,j}|\frac{1}{\sigma_{ii}\sigma_{jj}}-\frac{1}{\gamma_{ii}\gamma_{jj}}|\sum_{i,j=1}^p\gamma_{ij}^2
\leq K\max_{j}|\gamma_{jj}-\sigma_{jj}|\tr \boldsymbol{\Gamma}_p^2.
\end{align}
Thus, by \eqref{boundT1}, \eqref{gsbound} and \eqref{S3bound}, we have as $p\to \infty$
\begin{align}\label{S3lim}
    \frac{S_3}{\tr\bT_p^2}\leq \frac{K\max_{j}|\gamma_{jj}-\sigma_{jj}|\tr \boldsymbol{\Gamma}_p^2}{\tr\bT_p^2}
    \to 0.
\end{align}
Therefore, by \eqref{VI2}, \eqref{grg2t2}, \eqref{S1lim}, \eqref{S2lim} and \eqref{S3lim}, we have as $p\to\infty$
\begin{align}\label{VI2lim}
    \Var(I_2)\to 0.
\end{align}
 Thus, by \eqref{EI1}-\eqref{EI2}, \eqref{VI1} and \eqref{VI2lim}, \eqref{J1J2limit0} is proved. Hence, by \eqref{diffsum0} and \eqref{J1J2limit0}, we get
\begin{align*}
    \tilde{\tilde{T}}_{\tsum}=\tilde{T}_{\tsum}+o_p(1).
\end{align*}
\end{proof}

\subsection{Proof of Theorem \ref{two-sample-sum}: sum-type test for two-sample case}\label{pftwo-sample-sum}

\begin{proof}
    
The proof shares same spirit as Theorem \ref{one-sample-sum}.
Under the normality assumption and the null hypothesis in  \eqref{hypotest1}, we have $\bar{\bQ}^{(1)}-\bar{\bQ}^{(2)}\sim N(c\bI_p,\frac{n_1+n_2}{n_1n_2}\boldsymbol{\Sigma}_p)$ and then $\Big(\frac{n_1n_2}{n_1+n_2}\Big)^{1/2}(\bar{\bQ}^{(1)}-\bar{\bQ}^{(2)}-c\bI_p)
\sim N(0, \boldsymbol{\Sigma}_p)$. Let $\bD_{\boldsymbol{\Sigma}_p}=\mbox{diag}(\sigma_{11}, \cdots, \sigma_{pp})$ be the diagonal matrix of $\boldsymbol{\Sigma}_p$. Recall that $\bR_p=\bD_{\boldsymbol{\Sigma}_p}^{-1/2}\boldsymbol{\Sigma}_p\bD_{\boldsymbol{\Sigma}_p}^{-1/2}.$  Then
\begin{align*}
\Big(\frac{n_1n_2}{n_1+n_2}\Big)^{1/2}\bD_{\boldsymbol{\Sigma}_p}^{-1/2}(\bar{\bQ}^{(1)}-\bar{\bQ}^{(2)}-c\bI_p)
\sim N(0, \bR_p).
\end{align*}
Let $\bV_2:=\Big(\frac{n_1n_2}{n_1+n_2}\Big)^{1/2}\bG\bD_{\boldsymbol{\Sigma}_p}^{-1/2}(\bar{\bQ}^{(1)}-\bar{\bQ}^{(2)}-c\bI_p)=\Big(\frac{n_1n_2}{n_1+n_2}\Big)^{1/2}\bG\bD_{\boldsymbol{\Sigma}_p}^{-1/2}(\bar{\bQ}^{(1)}-\bar{\bQ}^{(2)}):=(V_{1,2},\cdots,V_{p,2})^T$, then $\bV\sim N(\boldsymbol{0}, \bG\bR_p\bG)$.
By Theorem \ref{theorem_sum}, we have
\begin{align*}
    \tilde{T}_{\tsum,2}=\frac{V_{1,2}^2+\cdots +V_{p,2}^2-\tr(\bG\bR_p\bG)}{\sqrt{2tr(\bG\bR_p\bG)^2}} \overset{d}{\to} N(0,1).
\end{align*}
We will show that
\begin{align}
T_{\tsum,2}\label{twosum3}
=\frac{\frac{n_1n_2}{n_1+n_2}(\bar{\bY}^{(1)}-\bar{\bY}^{(2)})^T \bD_{\boldsymbol{\Gamma}_p}^{-1}(\bar{\bY}^{(1)}-\bar{\bY}^{(2)})-p}{\sqrt{2\tr(\bT_p^2)}}+o_p(1) 
=\tilde{T}_{sum,2}+o_p(1).
\end{align}
Let $\tilde{\tilde{T}}_{\tsum,2}=\frac{\frac{n_1n_2}{n_1+n_2}(\bar{\bY}^{(1)}-\bar{\bY}^{(2)})^T \bD_{\boldsymbol{\Gamma}_p}^{-1}(\bar{\bY}^{(1)}-\bar{\bY}^{(2)})-p}{\sqrt{2\tr(\bT_p^2)}}$, then
\begin{align*}
    &\ \ \tilde{\tilde{T}}_{\tsum,2}- \tilde{T}_{\tsum,2} = \frac{\frac{n_1n_2}{n_1+n_2}(\bar{\bY}^{(1)}-\bar{\bY}^{(2)})^T \bD_{\boldsymbol{\Gamma}_p}^{-1}(\bar{\bY}^{(1)}-\bar{\bY}^{(2)})-p}{\sqrt{2\tr(\bT_p^2)}}- \frac{V_{1,2}^2+\cdots +V_{p,2}^2-\tr(\bG\bR_p\bG)}{\sqrt{2\tr(\bG\bR_p\bG)^2}}\nonumber\\
     &=
    \frac{\sqrt{2\tr(\bG\bR_p\bG)^2}[\frac{n_1n_2}{n_1+n_2}\{(\bar{\bQ}^{(1)}-\bar{\bQ}^{(2)})^T (\bG\bD_{\boldsymbol{\Gamma}_p}^{-1}\bG-\bD_{\boldsymbol{\Sigma}_p}^{-1/2}\bG\bD_{\boldsymbol{\Sigma}_p}^{-1/2})(\bar{\bQ}^{(1)}-\bar{\bQ}^{(2)})\}-\{p-\tr(\bG\bR_p\bG)\}]}{\sqrt{2\tr\bT_p^2}\sqrt{2\tr(\bG\bR_p\bG)^2}}\nonumber\\
    &\ \ \ \ +\frac{(\sqrt{2\tr(\bG\bR_p\bG)^2}-\sqrt{2\tr\bT_p^2})(\frac{n_1n_2}{n_1+n_2}(\bar{\bQ}^{(1)}-\bar{\bQ}^{(2)})^T\bD_{\boldsymbol{\Sigma}_p}^{-1/2}\bG\bD_{\boldsymbol{\Sigma}_p}^{-1/2}(\bar{\bQ}^{(1)}-\bar{\bQ}^{(2)})-\tr(\bG\bR_p\bG))}{\sqrt{2\tr\bT_p^2}\sqrt{2\tr(\bG\bR_p\bG)^2}}\nonumber\\
    &:=J_1+J_2.
\end{align*}
By the same discussion in the proof of Theorem \ref{theorem_sum}, we have
\begin{align}\label{twoJ1J2limit}
    J_1\to 0 \ \text{and} \  J_2\to 0  \ \text{in probability  as} \ p\to\infty.
\end{align}
and hence
\begin{align}\label{twottt}
    \tilde{\tilde{T}}_{\tsum,2} = \tilde{T}_{\tsum,2} +o(1).
\end{align}
Similar to the proof of Theorem 3.1  in \cite{Srivastava2009ATF},
we have
\begin{align}\label{twott}
    T_{\tsum,2} = \tilde{\tilde{T}}_{\tsum,2} +o(1).
\end{align}
Thus, by \eqref{twottt}-\eqref{twott}, \eqref{twosum3} is proved.
Hence, $T_{\tsum,2}\overset{d}{\to} N(0,1)$ as $p\to\infty$.
\end{proof}

\subsection{Proof of Non-Gaussian case}

\subsubsection{Proof of Theorem \ref{onesample-sum-nongaussian}}\label{pfonesample-sum-nongaussian}
\begin{proof}
    We first state Lemma \ref{onesample-mainlemma}, which relaxes the conditions of Theorem 2.1 in \cite{Srivastava2009ATF}.\begin{lemma}\label{onesample-mainlemma}
    Let $\bar{\bU} = (\bar{U}_{1},\dots,\bar{U}_{p})^T$ where $\bar{U}_{j}=1/n \sum_{i=1}^{n} U_{ij}$, $i=1,\dots,n$, $j=1,\dots,p$. $U_{ij}$ be i.i.d. random variables with zero mean, variance 1, and finite fourth moment $\theta$. Then for any $p\times p$ symmetric matrix $\mathbf{A} = (a_{ij})$ for which
$\lim_{p\to\infty}(\tr(\mathbf{A}^2)/p)$ exists,
\begin{equation*}
   \lim_{(p,n)\to\infty}P[(\frac{n \bar{\bU}^{T}\bA \bar{\bU} -\tr\bA}{\sqrt{2p\tau_2}})\leq x] = \Phi(x)
\end{equation*}
Where $\Phi(x)$ is standard normal distribution, and $\tau_2 = \frac{\tr(\mathbf{A}^2)}{p}$.
\end{lemma}
The proof of Lemma \ref{onesample-mainlemma} and remaining of Theorem \ref{onesample-sum-nongaussian} can be fully referenced in \ref{pftwosamplemaintheorem}.
\end{proof}
\subsubsection{Proof of Lemma \ref{twosample-consistencycor}}\label{pftwosample-consistencylemma}
\begin{proof}
To show Lemma \ref{twosample-consistencycor}, it is sufficient to establish the following Lemma \ref{twosample-consistencylemma}. Indeed, Lemma \ref{twosample-consistencylemma} asserts the consistency of the estimators for a general covariance matrix $\boldsymbol{\Sigma}_p$, while Lemma \ref{twosample-consistencycor} specifically delineates this property for $\boldsymbol{\Gamma}_p = \boldsymbol{G}\boldsymbol{\Sigma}_p\boldsymbol{G}$.
\begin{lemma}\label{twosample-consistencylemma}
    Let $\hat{\sigma}^{-1}_{ii}$ be diagonal elements of sample covariance $\hat{\bar{\boldsymbol{\Sigma}}}$,  where $ \hat{\bar{\boldsymbol{\Sigma}}}=\frac{1}{n_1+n_2}\Big[\sum_{i=1}^{n_1}
(\bQ_{i}^{(1)}-\bar{\bQ}^{(1)})(\bQ_{i}^{(1)}-\bar{\bQ}^{(1)})^T+\sum_{i=1}^{n_2}
(\bQ_{i}^{(2)}-\bar{\bQ}^{(2)})(\bQ_{i}^{(2)}-\bar{\bQ}^{(2)})^T\Big],$ and $ \bQ_{i}^{(k)} = \log \bW_{i}^{(k)}, k=1,2 $. If $\lim_{p\to\infty}\frac{\tr(\bR_p^i)}{p}$ exist for $i=2,4$, then 
\begin{itemize}
    \item [(i)] $\hat{\sigma}^{-1}_{ii}$ are consistent estimators of estimators of $\sigma_{ii}^{-1}$,
\item[(ii)]$\frac{1}{p}[\tr({\bD}_{\hat{\bar{\boldsymbol{\Sigma}}}_p}^{-1/2} \hat{\bar{\boldsymbol{\Sigma}}}_p{\bD}_{\hat{\bar{\boldsymbol{\Sigma}}}_p}^{-1/2})^2-(p^2/(N-2))]$ is a consistent estimator of $\tr \bR_{p}^2/p$.
\end{itemize}
where ${\bD}_{\hat{\bar{\boldsymbol{\Sigma}}}_p}$ is diagonal matrix of $\hat{\bar{\boldsymbol{\Sigma}}}_p$.
\end{lemma}
Part (i) can be derived using Corollary 2.6 in \cite{Srivastava2009ATF}. Indeed, part (ii) demonstrates the consistency of the estimator $\frac{1}{p}[\tr({\bD}_{\hat{\bar{\boldsymbol{\Sigma}}}_p}^{-1/2} \hat{\bar{\boldsymbol{\Sigma}}}_p{\bD}_{\hat{\bar{\boldsymbol{\Sigma}}}_p}^{-1/2})^2-(p^2/(N-2))]$ based on the correlation matrix. To establish this, we will now prove the corresponding result for the covariance matrix in Lemma \ref{prooflemma}. Subsequently, by following the methods outlined in the Appendix of \cite{Srivastava2008ATF}, part (ii) can be derived.
\begin{lemma}\label{prooflemma}
    If $N=O(p^{\epsilon})$, $0<\epsilon\leq1$, and $\lim_{p\to\infty} \tr(\Sigma^i_p)/p$ exist for i=2,4, then $\hat{\delta}_1$ and $\hat{\delta}_2$ are consistent estimators of $\delta_1$ and $\delta_2$ as $(N,p)\to\infty$, where $ \delta_1 = \tr(\boldsymbol{\Sigma}_p)/p$, $\hat{\delta}_1 = \tr (\hat{\bar{\boldsymbol{\Sigma}}}) /p$, $\delta_2 = \tr(\boldsymbol{\Sigma}_p)^2 / p$ and $\hat{\delta}_2 = \frac{1}{p}[\tr (\hat{\bar{\boldsymbol{\Sigma}}})^2 - (\tr\hat{\bar{\boldsymbol{\Sigma}}})^2 /(N-2) ]$.
\end{lemma}

\textbf{Proof of Lemma \ref{prooflemma}.} Note that
$$
\begin{aligned}
\hat{\bar{\boldsymbol{\Sigma}}} &= \frac{1}{N}[\sum_{i=1}^{n_1}(\bQ_{i1}-\bar{\bQ}^{(1)})(\bQ_{i1}-\bar{\bQ}^{(1)})^T+\sum_{i=1}^{n_2}(\bQ_{i2}-\bar{\bQ}^{(2)})(\bQ_{i2}-\bar{\bQ}^{(2)})^T]\\
&=\frac{1}{n_1+n_2}[\sum_{i=1}^{n_1}(\bQ_{i1}-\boldsymbol{\mu}_1^W)(\bQ_{i1}-\boldsymbol{\mu}_1^W)^T+\sum_{i=1}^{n_2}(\bQ_{i2}-\boldsymbol{\mu}_2^W)(\bQ_{i2}-\boldsymbol{\mu}_2^W)^T]+\Delta
\end{aligned}
$$
where  
$\Delta = -\frac{n_1}{n}(\bar{\bQ}^{(1)}-\boldsymbol{\mu}_1^W)(\bar{\bQ}^{(1)}-\boldsymbol{\mu}_1^W)^T-\frac{n_2}{n}(\bar{\bQ}^{(2)}-\boldsymbol{\mu}_2^W)(\bar{\bQ}^{(2)}-\boldsymbol{\mu}_2^W)^T$. 
For large $n$, $\Delta \to 0$ in probability. Thus $\hat{\bar{\boldsymbol{\Sigma}}}$ can be approximated by 
$$
\begin{aligned}
\hat{\bar{\boldsymbol{\Sigma}}}^* = \frac{1}{n_1+n_2}[\sum_{i=1}^{n_1}(\bQ_{i1}-\boldsymbol{\mu}_1^W)(\bQ_{i1}-\boldsymbol{\mu}_1^W)^T+\sum_{i=1}^{n_2}(\bQ_{i2}-\boldsymbol{\mu}_2^W)(\bQ_{i2}-\boldsymbol{\mu}_2^W)^T]
= \frac{1}{n_1+n_2}\sum_{i=1}^{n_1+n_2}\boldsymbol{\Sigma_p}^{1/2}\bU_i \bU_{i}^{T}\boldsymbol{\Sigma_p}^{1/2}
\end{aligned}
$$ where $\bU_i, i=1,2,\dots, n_1+n_2$, are i.i.d with mean $\boldsymbol{0}$ and covariance matrix $I_p$.
According to the proof of Theorem 2.2 in \cite{Srivastava2009ATF}, Lemma \ref{prooflemma} holds.

The proof of Lemma \ref{twosample-consistencylemma} (ii) from Lemma \ref{prooflemma} has been shown in detail in \cite{Srivastava2008ATF}, and we will only briefly describe the idea here. Let $c_{ii} = 1+ 1/\hat{\sigma}_{ii}$, $D_c = diag(c_{11},\dots,c_{pp})$ and $D_s = diag(\hat{\sigma}_{11},\dots,\hat{\sigma}_{pp})$. Since $\frac{1}{p}[\tr({\bD}_{\hat{\bar{\boldsymbol{\Sigma}}}_p}^{-1/2} \hat{\bar{\boldsymbol{\Sigma}}}_p{\bD}_{\hat{\bar{\boldsymbol{\Sigma}}}_p}^{-1/2})^2-(p^2/(N-2))]$ is invariant under the scalar transformations of each component of $\bQ^{(k)}_i, i = 1,2,\dots,n_k,k=1,2$, we may as well assume that $\boldsymbol{\Sigma}_p = \bR_p$. We write
$$
\begin{aligned}
\frac{1}{p}[\tr({\bD}_{\hat{\bar{\boldsymbol{\Sigma}}}_p}^{-1/2} \hat{\bar{\boldsymbol{\Sigma}}}_p{\bD}_{\hat{\bar{\boldsymbol{\Sigma}}}_p}^{-1/2})^2-(p^2/(N-2))] &= \frac{1}{p}\left[\tr(D_s^{-1}\hat{\bar{\boldsymbol{\Sigma}}})^2-\frac{\tr(D_s^{-1}\hat{\bar{\boldsymbol{\Sigma}}})^2}{N-2}\right]\\
&=\frac{1}{p}\left[\tr(\hat{\bar{\boldsymbol{\Sigma}}})^2-\frac{\tr(\hat{\bar{\boldsymbol{\Sigma}}})^2}{N-2}\right]-\frac{2}{p}\left[ \tr (D_c(\hat{\bar{\boldsymbol{\Sigma}}})^2) - \frac{\tr(\hat{\bar{\boldsymbol{\Sigma}}})\tr (D_c \hat{\bar{\boldsymbol{\Sigma}}})}{N-2}\right]\\
&\quad + \frac{1}{p}\left[\tr(D_c\hat{\bar{\boldsymbol{\Sigma}}})^2-\frac{\tr(D_c\hat{\bar{\boldsymbol{\Sigma}}})^2}{N-2}\right].
\end{aligned}
$$
Then by assumptions and Lemma \ref{prooflemma}, we can verify that as $(N,p) \to\infty$, the first term converges to $\frac{\tr(\bR_p^2)}{p}$ in probability, and the last two terms converges to zero in probability. So far we have proved that Lemma \ref{twosample-consistencylemma} holds. Thus, by Lemma \ref{twosample-consistencylemma} and \eqref{Gtrb}, Lemma \ref{twosample-consistencycor} holds.
\end{proof}
\subsubsection{Proof of Theorem \ref{twosamplemaintheorem}}\label{pftwosamplemaintheorem}
\begin{proof}
We derive the asymptotic distribution of the test statistic using Lemma \ref{twosample-mainlemma}, which is specifically designed for compositional data. In establishing the asymptotic distribution of the test statistic, we demonstrate that the centered log-ratio vector $\bY^{(k)}_i$ still satisfies properties similar to those of $\log \bW_i^{(k)}$ under Assumptions \ref{assumpa2}, \ref{assumpbb1}, and \ref{addassume1}.

\begin{lemma}\label{twosample-mainlemma}
    Let $\bar{\bU}^{(k)} = (\bar{U}_{1}^{(k)},\dots,\bar{U}_{p}^{(k)})^T$ where $\bar{U}_{j}^{(k)}=1/n_k \sum_{i=1}^{n_k} U_{ij}^{(k)}$, $i=1,\dots,n_k$, $j=1,\dots,p$, $k=1,2$, $U_{ij}^{(k)}$ be i.i.d. random variables with zero mean, variance 1, and finite fourth moment $\theta$. Then for any $p\times p$ symmetric matrix $\mathbf{A} = (a_{ij})$ for which
$\lim_{p\to\infty}(\tr(\mathbf{A}^2)/p)$ exists,
\begin{equation*}
   \lim_{(p,n)\to\infty}P[(\frac{\frac{n_1n_2}{n_1+n_2}(\bar{\bU}^{(1)}-\bar{\bU}^{(2)})^{T}\bA(\bar{\bU}^{(1)}-\bar{\bU}^{(2)})-\tr\bA}{\sqrt{2p\tau_2}})\leq x] = \Phi(x)
\end{equation*}
Where $\Phi(x)$ is standard normal distribution, and $\tau_2 = \frac{\tr(\mathbf{A}^2)}{p}$.
\end{lemma}

\textbf{Proof of Lemma \ref{twosample-mainlemma}.}
Lemma \ref{twosample-mainlemma} is proved by Theorem 2.1 in \cite{wang2014jointclt}. First, some calculations for $\bar{U}_j^{(k)}$ are listed below.
Note that,
\begin{align*}
     \mathbb{E}(\sqrt{n_k}\bar{U}_1^{(k)})^2 = 1,\
     \mathbb{E}(\sqrt{n_1}\bar{U}_1^{(1)})(\sqrt{n_2}\bar{U}_1^{(2)}) = 0,\
     \mathbb{E}(\sqrt{n_k}\bar{U}_1^{(k)})^4 = \frac{\theta}{n_k}+\frac{3(n_k-1)}{n_k}.
\end{align*}
Let 
\begin{align*}
    T^* &= \frac{1}{\sqrt{p}}\left(\frac{n_1n_2}{n_1+n_2}(\bar{\bU}^{(1)}-\bar{\bU}^{(2)})^{T}\bA(\bar{\bU}^{(1)}-\bar{\bU}^{(2)})-tr\bA\right)\\
    &=\frac{1}{\sqrt{p}}\left( \frac{n_2}{n_1+n_2}\left(n_1 (\bar{U}^{(1)})^{T}\bA \bar{U}^{(1)} - tr\bA\right)\right) + \frac{1}{\sqrt{p}}\left( \frac{n_1}{n_1+n_2}\left(n_2 (\bar{U}^{(2)})^{T}\bA \bar{U}^{(2)} - tr\bA\right)\right)\\ &\quad - \frac{2}{\sqrt{p}}\left( \frac{\sqrt{n_1n_2}}{n_1+n_2}\sqrt{n_1} (\bar{U}^{(1)})^{T}\bA \sqrt{n_2}\bar{U}^{(2)} \right)\\
    &:= T_1 + T_2 +T_3.
\end{align*}
By Theorem 2.1 in \cite{wang2014jointclt}, we have $T_1$, $T_2$, and $T_3$ each converge to a zero-mean Gaussian distribution with variances of $v_1$, $v_2$ and $v_3$ respectively, where
$$
v_1 = \lim_{(N,p)\to\infty} \frac{n_1^2}{(n_1+n_2)^2} \frac{2\tr\bA^2}{p}, \quad v_2 = \lim_{(N,p)\to\infty} \frac{n_2^2}{(n_1+n_2)^2} \frac{2\tr\bA^2}{p}, \quad \text{and } v_3 = \lim_{(N,p)\to\infty} \frac{4n_1n_2}{(n_1+n_2)^2} \frac{\tr\bA^2}{p}.$$
Finally, combining the uncorrelation of $T_1$, $T_2$, and $T_3$, Lemma \ref{twosample-mainlemma} can be proved.

\textbf{The remaining proof of Theorem \ref{twosamplemaintheorem}}. We assume that Assumptions \ref{assumpa2}, \ref{assumpbb1} and \ref{addassume1} hold. From Lemma \ref{twosample-consistencycor}, as $p\to \infty$
$$
[\tr( {\bD}_{\hat{\bar{\boldsymbol{\Gamma}}}_p}^{-1/2} \hat{\bar{\boldsymbol{\Gamma}}}_p{\bD}_{\hat{\bar{\boldsymbol{\Gamma}}}_p}^{-1/2})^2-(p^2/(N-2))]/p\overset{p}{\to}\tau_{20},$$
where $\tau_{20} = \lim_{p\to\infty} \frac{\tr(\bT_p^2)}{p}$. Thus, we need only consider the asymptotic distribution of the test statistic
\begin{align*}
T  =\frac{\frac{n_1n_2}{n_1+n_2}(\bar{\bY}^{(1)}-\bar{\bY}^{(2)})^{T} \bD_{\boldsymbol{\Gamma}_p}^{-1}(\bar{\bY}^{(1)}-\bar{\bY}^{(2)})-p}
{\sqrt{2p\tau_{20}}}
 = \frac{\frac{n_1n_2}{n_1+n_2}(\bar{\bU}^{(1)}-\bar{\bU}^{(2)})^{T} \bB (\bar{\bU}^{(1)}-\bar{\bU}^{(2)})-p}
{\sqrt{2p\tau_{20}}},
\end{align*}
where $\bB = \boldsymbol{\Sigma}^{1/2}_p\bG \bD_{\boldsymbol{\Gamma}_p}^{-1} \bG\boldsymbol{\Sigma}^{1/2}_p$.
By Assumptions \ref{assumpbb1}, \ref{addassume1}, similar to the proof of \eqref{GtrT}, we have for $i=2,4$,
\begin{align}\label{Gtrb}
    \lim\limits_{p\to\infty} \frac{\tr (\bT_p^i)}{p} < \infty.
\end{align}
Thus, by \eqref{Gtrb} and  Lemma \ref{twosample-mainlemma}, Theorem \ref{twosamplemaintheorem} is proved.
\end{proof}

\subsubsection{Proof of Proposition \ref{twosample-power}}\label{pftwosample-power}
\begin{proof}
By Theorem \ref{twosamplemaintheorem}, as $(n,p)\to\infty$.
$$
T =\frac{\frac{n_1n_2}{n_1+n_2}(\bar{\bY}_1-\bar{\bY}_2-\boldsymbol{\mu}^{X})^TD_{\boldsymbol{\Gamma}_p}^{-1}(\bar{\bY}_1-\bar{\bY}_2-\boldsymbol{\mu}^X)-p}{[2\tr \boldsymbol{T}_p^2]^{1/2}}
$$
has a standard normal distribution, where $\boldsymbol{\mu}^X = \boldsymbol{\mu}_1^X-\boldsymbol{\mu}_2^X$. For the alternative, 
$$
\begin{aligned}
&\quad \frac{1}{\sqrt{p}}[\frac{n_1n_2}{n_1+n_2}(\bar{\bY}_1-\bar{\bY}_2-\boldsymbol{\mu}^X)^TD_{\boldsymbol{\Gamma}_p}^{-1}(\bar{\bY}_1-\bar{\bY}_2-\boldsymbol{\mu}^X)] \\
& = \frac{1}{\sqrt{p}}\frac{n_1n_2}{n_1+n_2}(\bar{\bY}_1-\bar{\bY}_2)^TD_{\boldsymbol{\Gamma}_p}^{-1}(\bar{\bY}_1-\bar{\bY}_2) + \frac{1}{(N-2)\sqrt{p}}\boldsymbol{\delta}^T\boldsymbol{G}D_{\boldsymbol{\Gamma}_p}^{-1}\boldsymbol{G}\boldsymbol{\delta}\\
&\quad - \frac{2}{\sqrt{p(n_1+n_2-2)}}\sqrt{\frac{n_1n_2}{n_1+n_2}}\boldsymbol{\delta}^T\boldsymbol{G}D_{\boldsymbol{\Gamma}_p}^{-1}(\bar{\bY}_1-\bar{\bY}_2).
\end{aligned}$$
Under \ref{dddm} and \eqref{twosumH1}, 
$$
\begin{aligned}
\frac{1}{\sqrt{p}}[\frac{n_1n_2}{n_1+n_2}(\bar{\bY}_1-\bar{\bY}_2-\boldsymbol{\mu}^X)^T\boldsymbol{G}D_{\boldsymbol{\Gamma}_p}^{-1}\boldsymbol{G}(\bar{\bY}_1-\bar{\bY}_2-\boldsymbol{\mu}^X)]\\ \overset{p}{\to} \frac{1}{\sqrt{p}}[\frac{n_1n_2}{n_1+n_2}(\bar{\bY}_1-\bar{\bY}_2)^TD_{\boldsymbol{\Gamma}_p}^{-1}(\bar{\bY}_1-\bar{\bY}_2)- \frac{\boldsymbol{\delta}^TD_{\boldsymbol{\Gamma}_p}^{-1}\boldsymbol{\delta}}{n_1+n_2-2}].
\end{aligned}$$
Thus, Proposition \ref{twosample-power} holds.
\end{proof}

\section{Proof of Section 4: max-type test}
\subsection{Proof of Lemma \ref{theorem_max}}\label{pftheorem_max}
\begin{proof}
The central essence of the proof hinges on expressing $H_i=Z_i-\bar{Z}$ as $H_i=Z_i-\Expe \bar{Z}+\Expe \bar{Z}-\bar{Z}$ for $i=1,\ldots,p$, where $(Z_1-\Expe \bar{Z},\ldots,Z_p-\Expe \bar{Z})^T$ possesses a positive covariance matrix $\boldsymbol{\Psi}_p$, $\bar{Z}=\frac{1}{p}\sum_{i=1}^p Z_i$ and $\Expe \bar{Z}=c$. The analysis involves considering $H_i=Z_i-\Expe \bar{Z}+\Expe \bar{Z}-\bar{Z}$ ($i=1,\ldots,p$) in two events: $\{H_i=Z_i-\bar{Z}\}\cap\{|\bar{Z}-\Expe \bar{Z}|\leq \epsilon_p\}$ and $\{H_i=Z_i-\bar{Z}\}\cap\{|\bar{Z}-\Expe \bar{Z}|\geq \epsilon_p\}$ for any $\epsilon_p>0$. Here, the event $\{|\bar{Z}-\Expe \bar{Z}|\geq \epsilon_p\}$ occurs with an exponentially small probability, as dictated by concentration inequalities. Building on this foundation, we establish the following bound
\begin{align}\label{Ymaxbound}
I_{lower}(Z,\epsilon_p)+o_p(1) \leq P\Big(\max_{1\leq i \leq p}|H_i|>y\Big)\leq I_{upper}(Z,\epsilon_p)+o_p(1),
\end{align}
where  
\begin{align}\label{y}
y=\big(2\log p - \log \log p+x\big)^{1/2}
\end{align}
for any $x\in \mathbb{R}$, and
$\lim_{p\to\infty}I_{lower}(Z,\epsilon_p)=\lim_{p\to\infty}I_{upper}(Z,\epsilon_p)=\lim_{p\to\infty}P\Big(\max_{1\leq i \leq p}|Z_i-\Expe\bar{Z}|>y\Big)= 1-\exp\Big(-\frac{1}{\sqrt{\pi}}e^{-x/2}\Big)$, and partition the proof of this bound into two steps.

\textbf{Step 1.} 
By Bonferroni inequlity,
\begin{align}\label{bonfer}
I^Y_{lower}:=\sum_{t=1}^{2k}(-1)^{t-1}E_t\leq P\Big(\max_{1\leq i\leq p}|H_i|>y\Big) \leq I^Y_{upper}:=\sum_{t=1}^{2k+1}(-1)^{t-1}E_t
\end{align}
for any $p/2\geq k\geq 1$, where $E_t:=\sum_{1\leq i_1 < \dots< i_t\leq p} P(|H_{i_1}|>y, \cdots, |H_{i_t}|>y) $. 

Note that, $\bar{Z}=\frac{1}{p}\sum_{i=1}^p 
Z_i$, $\Expe \bar{Z}=c$,
\begin{align}\label{proy1pdecom1}
 P(|H_{i_1}|>y, \cdots, |H_{i_t}|>y)
 &=P(|Z_{i_1}-\bar{Z}|>y, \cdots, |Z_{i_t}-\bar{Z}|>y,|\bar{Z}-\Expe \bar{Z}|< \epsilon_p)\nonumber\\
 &\ \ +P(|Z_{i_1}-\bar{Z}|>y, \cdots, |Z_{i_t}-\bar{Z}|>y,|\bar{Z}-\Expe \bar{Z}|\geq \epsilon_p)\nonumber\\
  &\leq P(|Z_{i_1}-\Expe \bar{Z}|>y-\epsilon_p, \cdots, |Z_{i_t}-\Expe \bar{Z}|>y-\epsilon_p)\nonumber\\
 &\ \ +P(|\bar{Z}-\Expe \bar{Z}|\geq \epsilon_p),
\end{align}
and similarly
\begin{align}\label{proy1pdecom2}
 P(|H_{i_1}|>y, \cdots, |H_{i_t}|>y)
 &\geq P(|Z_{i_1}-\Expe \bar{Z}|>y+\epsilon_p, \cdots, |Z_{i_t}-\Expe \bar{Z}|>y+\epsilon_p)\nonumber\\
 &\ \ -P(|\bar{Z}-\Expe \bar{Z}|\geq \epsilon_p).
\end{align}
Therefore, by \eqref{proy1pdecom1} and \eqref{proy1pdecom2}, we get
   \begin{align}\label{yzbound}
&I^Y_{lower}\geq I_{lower}(Z,\epsilon_p)+\sum_{t=1}^{2k}(-1)^{t-1}\sum_{1\leq i_1 < \dots< i_t\leq p}P\Big(|\bar{Z}-\Expe \bar{Z}|>\epsilon_p\Big),\nonumber\\
&I^Y_{upper}\leq I_{upper}(Z,\epsilon_p)+\sum_{t=1}^{2k+1}(-1)^{t-1}\sum_{1\leq i_1 < \dots< i_t\leq p} P\Big(|\bar{Z}-\Expe \bar{Z}|>\epsilon_p\Big),
\end{align}
where
$$I_{lower}(Z,\epsilon_p)=\sum_{t=1}^{2k}(-1)^{t-1}\sum_{1\leq i_1 < \dots< i_t\leq p} P(|Z_{i_1}-\Expe \bar{Z}|>y+(-1)^{t-1}\epsilon_p, \cdots, |Z_{i_t}-\Expe \bar{Z}|>y+(-1)^{t-1}\epsilon_p)$$
and
$$I_{upper}(Z,\epsilon_p) =\sum_{t=1}^{2k+1}(-1)^{t-1}\sum_{1\leq i_1 < \dots< i_t\leq p} P(|Z_{i_1}-\Expe \bar{Z}|>y+(-1)^{t}\epsilon_p, \cdots, |Z_{i_t}-\Expe \bar{Z}|>y+(-1)^{t}\epsilon_p).$$

\textbf{Step 2.} 
Note that by Theorem 2 in \cite{Feng2022AsymptoticIO}, we get
\begin{align*}
    &\ \ \ \ P(|Z_{i_1}-\Expe \bar{Z}|>y+(-1)^{t}\epsilon_p, \cdots, |Z_{i_t}-\Expe Z|>y+(-1)^{t}\epsilon_p)\nonumber\\
    &=(1+o(1))\pi^{-t/2}p^{-t}\exp(-\frac{tx}{2}).
\end{align*}
Furthermore,  
 \begin{align}\label{Qt}
&\lim_{p\rightarrow\infty}\sum_{1\leq i_1 < \dots< i_t\leq p} P(|Z_{i_1}-\Expe \bar{Z}|>y+(-1)^{t-1}\epsilon_p, \cdots, |Z_{i_t}-\Expe \bar{Z}|>y+(-1)^{t-1}\epsilon_p) \nonumber\\
&=\lim_{p\rightarrow\infty}\sum_{1\leq i_1 < \dots< i_t\leq p} P(|Z_{i_1}-\Expe \bar{Z}|>y+(-1)^{t}\epsilon_p, \cdots, |Z_{i_t}-\Expe \bar{Z}|>y+(-1)^{t}\epsilon_p) \nonumber\\
&=\frac{1}{t!}\pi^{-t/2}e^{-tx/2}.
\end{align}
Hence, by \eqref{Qt} and  the same arguement in \cite{Feng2022AsymptoticIO} we have 
\begin{align}\label{boundy}
\lim_{p\rightarrow\infty}I_{upper}(Z,\epsilon_p) 
=\lim_{p\rightarrow\infty}I_{lower}(Z,\epsilon_p) =1-\exp\Big(-\frac{1}{\sqrt{\pi}}e^{-x/2}\Big).
\end{align}
Now we need to estimate $P(|\bar{Z}-\Expe \bar{Z}|\geq \epsilon_p)$. By concentration inequality, we get
\begin{align}\label{concineq}
    P\Big(|\bar{Z}-\Expe \bar{Z}|>\epsilon_p\Big)&= 2\int_{\epsilon_p}^{\infty} \frac{1}{\sqrt{2\pi\Var(\bar{Z})}}\exp(-\frac{|x-\Expe \bar{Z}|^2}{2\Var(\bar{Z})})dx\leq  2\exp(-\frac{\epsilon_p^2}{2\Var(\bar{Z})}).
\end{align}
Note that 
\begin{align}\label{varmeanz}
    \Var(\bar{Z})&=\frac{1}{p^2}\sum_{i,j=1}^p\psi_{ij}\leq \frac{1}{p}||\boldsymbol{\Psi}_p|_2 = \frac{1}{p}\lambda_{\max}(\boldsymbol{\Psi}_p).
\end{align}
Moreover, by \eqref{concineq} and \eqref{varmeanz}, we derive 
\begin{align}\label{smallpro10}
  \sum_{t=1}^{2k+1}(-1)^{t-1}\sum_{1\leq i_1 < \dots< i_t\leq p} P\Big(|\bar{Z}-\Expe \bar{Z}|>\epsilon_p\Big) =& P\Big(|\bar{Z}-\Expe \bar{Z}|>\epsilon_p\Big) \big(\sum_{t=1}^{2k+1}(-1)^{t-1}C_p^t \big)\nonumber\\
  \leq &(2k+1)p^tP\Big(|\bar{Z}-\Expe \bar{Z}|>\epsilon_p\Big) \nonumber\\
   \leq &C_k \exp\big(-[\frac{\sqrt{p}\epsilon_p^2\sqrt{p}(\log p)^{-1}}{2\lambda_{\max}(\boldsymbol{\Psi}_p)}-(2k+1)]\log p\big).
\end{align}
For each fixed $k\geq 1$, take $\epsilon_p=\Big((4k+4)\frac{\lambda_{\max}(\boldsymbol{\Psi}_p)}{\sqrt{p}(\log p)^{-1}}p^{-1/2}\Big)^{1/2}$ and let $p\to\infty$, under the condition $\lambda_{\max}(\boldsymbol{\Psi}_p)<<\sqrt{p}(\log p)^{-1}$ in Assumption \ref{assumption_max} and \eqref{smallpro10}, we have \begin{align}\label{smallpro1}
 \lim_{p\to0}\sum_{t=1}^{2k}(-1)^{t-1}\sum_{1\leq i_1 < \dots< i_t\leq p} P\Big(|\bar{Z}-\Expe \bar{Z}|>\epsilon_p\Big) = 0.
\end{align}
Similarly,
\begin{align}\label{smallpro2}
    \lim_{p\to0}\sum_{t=1}^{2k+1}(-1)^{t-1}\sum_{1\leq i_1 < \dots< i_t\leq p} P\Big(|\bar{Z}-\Expe \bar{Z}|>\epsilon_p\Big) = 0.
\end{align}
Finally, by \eqref{bonfer}, \eqref{yzbound}, and \eqref{smallpro1}-\eqref{smallpro2}, we have \eqref{Ymaxbound}. Hence, by \eqref{Ymaxbound} and \eqref{boundy}, we have
\begin{align}\label{maxdistri}
\lim_{p\to \infty}P\Big(\max_{1\leq i \leq p}|H_i|>y\Big)= 1-\exp\Big(-\frac{1}{\sqrt{\pi}}e^{-x/2}\Big).
\end{align}
\end{proof}

\subsection{Proof of Theorem \ref{one-sample-max}: max-type test for one-sample case}\label{pfone-sample-max}
\begin{proof}

We divide the proof of Theorem \ref{one-sample-max} into the following three steps.\\
    \textbf{Step 1:} We first show that
\begin{align}\label{onemax10}
    \tilde{\tilde{T}}_{\max}=\tilde{T}_{\max}+o_p(1).
\end{align}
where $\tilde{\tilde{T}}_{\max}=n\max_{1\leq i\leq p}\frac{\bar{Y}_i^2}{\gamma_{ii}}$ and $\tilde{T}_{\max}=\max_{1\leq i \leq p}V_i^2$.\\
 \textbf{Step 2:} Secondly, similar to the proof in (S.96) in \cite{Feng2022AsymptoticIO}, we derive
\begin{align}\label{onemax2}
T_{\max}=\tilde{\tilde{T}}_{\max}+o_p(1),
\end{align}
 \textbf{Step 3:} By Steps 1--2, we have
\begin{align}\label{onemax3}
T_{\max}=\tilde{T}_{\max}+o_p(1).
\end{align}
Moreover, by Lemma \ref{theorem_max} we have $\tilde{T}_{\max}$ converges to a Gumbel distribution with cdf $F(x)=\exp\{-\frac{1}{\sqrt{\pi}}e^{-x/2}\}$ as $p\to \infty$. Hence, $T_{\max}$ converges to a Gumbel distribution with cdf $F(x)=\exp\{-\frac{1}{\sqrt{\pi}}e^{-x/2}\}$ as $p\to \infty$.

\textbf{Therefore, the remaining part of the proof is dedicated to proving Step 1.}

\textbf{Step 1.} 
Since $\sqrt{n}(\bar{\bQ}-c\boldsymbol{1}_p)\sim N(\boldsymbol{0},\boldsymbol{\Sigma}_p)$, $\sqrt{n}\bar{\bY}=\sqrt{n}\bG\bar{\bQ}\sim N(\boldsymbol{0},\boldsymbol{\Gamma}_p)$ under $H_0$. Thus, $\sqrt{n}\bD_{\boldsymbol{\Sigma}_p}^{-1/2}(\bar{\bQ}-c\boldsymbol{1}_p)\sim N(\boldsymbol{0},\bR_p)$, $\bV:=\sqrt{n}\bG\bD_{\boldsymbol{\Sigma}_p}^{-1/2}(\bar{\bQ}-c\boldsymbol{1}_p):=(V_1,\cdots,V_p)^T\sim N(\boldsymbol{0},\bG\bR_p\bG)$ under $H_0$.

Recall that $\tilde{T}_{\max}=\max_{1\leq i \leq p}V_i^2$ and $\tilde{\tilde{T}}_{\max}=n\max_{1\leq i\leq p}\frac{\bar{Y}_i^2}{\gamma_{ii}}$, then
\begin{align}\label{onemaxcompara0}
|\tilde{\tilde{T}}_{\max}-\tilde{T}_{\max}|
=& \Big|n\max_{1\leq i\leq p}\frac{\bar{Y}_i^2}{\gamma_{ii}}-\max_{1\leq i\leq p}V_i^2\Big| \nonumber\\
\leq& \Big|n\max_{1\leq i\leq p}\frac{\bar{Y}_i^2}{\gamma_{ii}}\Big| \max_{1\leq i\leq p}\Big| 1-\frac{V_i^2}{n\bar{Y}_i^2/\gamma_{ii}}\Big|\nonumber\\
=&(\max_{1\leq i\leq p}\Big|\frac{\sqrt{n}\bar{Y}_i}{\sqrt{\gamma}_{ii}}\Big|)^2\max_{1\leq i\leq p}\Big|1-\frac{V_i^2}{n\bar{Y}_i^2/\gamma_{ii}}\Big|\nonumber.\\
\end{align}
Note that, $\bar{\bQ}=(\bar{Q}_1,\cdots,\bar{Q}_p)^T$, $\bar{Q}_k=\frac{1}{n}\sum_{j=1}^nQ_{jk}$ and 
\begin{align}
    (\sqrt{n}\frac{\bar{Y}_i}{\sqrt{\gamma_{ii}}})^2=n\bg_i^T(\bar{\bQ}-c\boldsymbol{1}_p)(\bar{\bQ}-c\boldsymbol{1}_p)^T\bg_i/\gamma_{ii}  &= n(-\frac{1}{p}\sum_{k=1}^p(\bar{Q}_k-c)+(\bar{Q}_i-c))^2\gamma_{ii}^{-1}, 
\end{align}
and 
\begin{align*}
    V_i^2=n\bg_i^T\bD_{\boldsymbol{\Sigma}_p}^{-1/2}(\bar{\bQ}-c\boldsymbol{1}_p)(\bar{\bQ}-c\boldsymbol{1}_p)^T \bD_{\boldsymbol{\Sigma}_p}^{-1/2}\bg_i  = n(-\frac{1}{p}\sum_{k=1}^p\sigma_{kk}^{-1/2}(\bar{Q}_k-c)+\sigma_{ii}^{-1/2}(\bar{Q}_i-c))^2.
\end{align*}
Let 
$$\bQ_{sd}=(\frac{\sqrt{n}(\bar{Q}_1-c)}{\sqrt{\sigma_{11}}}\cdots,\frac{\sqrt{n}(\bar{Q}_p-c)}{\sqrt{\sigma_{pp}}}):=(Q_{sd}(1),\cdots,Q_{sd}(p)),$$ 
and
$$\bQ_c=(\sqrt{n}(\bar{Q}_1-c),\cdots,\sqrt{n}(\bar{Q}_p-c)):=(Q_c(1),\cdots,Q_c(p)),$$ 
and
$$\bar{Q}_{sd}=\frac{1}{p}\sum_{i=1}^pQ_{sd}(i),\ \ \bar{Q}_c=\frac{1}{p}\sum_{i=1}^pQ_c(i),$$ 
then $\bQ_{sd}\sim N(\boldsymbol{0},\bR_p)$, $\bQ_c\sim N(\boldsymbol{0}, \boldsymbol{\Sigma}_p)$ and $\Expe \bar{Q}_{sd}=\Expe \bar{Q}_c=0$. Define $a_p=p^{-1/2-\alpha}(\log p)^{C_1/2}$ for $C_1>0$, for any $\epsilon_p>0$
\begin{align}\label{VY1}
    \Prob(\Big|1-\frac{V_i^2}{n\bar{Y}_i^2/\gamma_{ii}}\Big|>a_p)=&\Prob(\Big|1-\frac{V_i^2}{n\bar{Y}_i^2/\gamma_{ii}}\Big|>a_p,|\bar{Q}_{sd}|<\epsilon_p,|\bar{Q}_c|<\epsilon_p)\nonumber\\
    &+\Prob(\Big|1-\frac{V_i^2}{n\bar{Y}_i^2/\gamma_{ii}}\Big|>a_p,|\bar{Q}_{sd}|>\epsilon_p,|\bar{Q}_c|<\epsilon_p)\nonumber\\
    &+\Prob(\Big|1-\frac{V_i^2}{n\bar{Y}_i^2/\gamma_{ii}}\Big|>a_p,|\bar{Q}_{sd}|<\epsilon_p,|\bar{Q}_c|>\epsilon_p)\nonumber\\
    &+\Prob(\Big|1-\frac{V_i^2}{n\bar{Y}_i^2/\gamma_{ii}}\Big|>a_p,|\bar{Q}_{sd}|>\epsilon_p,|\bar{Q}_c|>\epsilon_p)\nonumber\\
    \leq& \Prob(\Big|1-\frac{V_i^2}{n\bar{Y}_i^2/\gamma_{ii}}\Big|>a_p,|\bar{Q}_{sd}|<\epsilon_p,|\bar{Q}_c|<\epsilon_p)+\Prob(|\bar{Q}_{sd}|>\epsilon_p)+2\Prob(|\bar{Q}_c|>\epsilon_p)\nonumber\\
    :=&I_1+I_2+I_3.
\end{align}
By concentration inequality, we have
\begin{align}\label{VY2}
    I_2=\Prob(|\bar{Q}_{sd}|>\epsilon_p)\leq 2\exp(-\frac{\epsilon_p^2}{2\Var(\bar{Q}_{sd})})\to 0 , as \ \ p\to \infty
\end{align}
and
\begin{align}\label{VY3}
    I_3=2\Prob(|\bar{Q}_c|>\epsilon_p)\leq 4\exp(-\frac{\epsilon_p^2}{2\Var(\bar{Q}_c)}) \to 0 , as \ \ p\to \infty
\end{align}
where $\Var(\bar{Q}_{sd})=\frac{1}{p^2}\sum_{i,j=1}^p\rho_{ij}\leq \frac{1}{p}\lambda_{\max}(\boldsymbol{\bR}_p)$, $\Var(\bar{Q}_c)=\frac{1}{p^2}\sum_{i,j=1}^p\sigma_{ij}\leq \frac{1}{p}\lambda_{\max}(\boldsymbol{\Sigma}_p)$, and 
\begin{align}
&\lambda_{\max}(\boldsymbol{\bR}_p)<<\sqrt{p}(\log p)^{-1},\label{ei1}\\
    &\lambda_{\max}(\boldsymbol{\Sigma}_p)=||\bD_{\boldsymbol{\Sigma}_p}^{1/2}\bR_p\bD_{\boldsymbol{\Sigma}_p}^{1/2}||_2
    \leq ||\bD_{\boldsymbol{\Sigma}_p}^{1/2}||_2||\bR_p||_2||\bD_{\boldsymbol{\Sigma}_p}^{1/2}||_2
    <<\sqrt{p}(\log p)^{-1},\label{ei2}
\end{align}
here \eqref{ei1} is derived from Assumption \ref{assumption_max} and \eqref{ei2} is derived from Assumption \ref{assumption_max} and Assumption \ref{addassume1}. Moreover,
\begin{align}\label{VY4}
    I_1 
    &\leq \Prob(1-\frac{(|Q_{sd}(i)|-\epsilon_p)^2}{(|Q_c(i)|+\epsilon_p)^2/\gamma_{ii}}>a_p,|\bar{Q}_{sd}|<\epsilon_p,|\bar{Q}_c|<\epsilon_p)\nonumber\\
&+ \Prob(1-\frac{(|Q_{sd}(i)|+\epsilon_p)^2}{(|Q_c(i)|-\epsilon_p)^2/\gamma_{ii}}<-a_p,|\bar{Q}_{sd}|<\epsilon_p,|\bar{Q}_c|<\epsilon_p)\nonumber\\
&=I_{11}+I_{12}.
\end{align}
And by \eqref{gsbound},
\begin{align}\label{VY5}
    I_{11}
    &=\Prob(1-\frac{(|Q_{sd}(i)|-\epsilon_p)^2}{(|Q_c(i)|+\epsilon_p)^2/\gamma_{ii}}>a_p,|\bar{Q}_{sd}|<\epsilon_p,|\bar{Q}_c|<\epsilon_p)\nonumber\\
    &\leq \Prob(1-\frac{(|Q_{sd}(i)|-\epsilon_p)^2}{(|Q_{sd}(i)|+\frac{\epsilon_p}{\sqrt{\sigma_{ii}}})^2}\frac{\gamma_{ii}}{\sigma_{ii}}>a_p)\to 0 , as \ \ p\to \infty,
\end{align}
and similarly
\begin{align}\label{VY6}
    I_{12}\to 0 , as \ \ p\to \infty.
\end{align}
Therefore, by \eqref{VY1}-\eqref{VY6},
\begin{align}\label{gdz}
     \max_{1\leq i\leq p}\Big| 1-\frac{V_i^2}{(\sqrt{n}\frac{\bar{Y}_i}{\sqrt{\gamma_{ii}}})^2}\Big| = O_p(a_p).
\end{align}
Moreover, since $\sqrt{n}\bar{\bY}\sim N(\boldsymbol{0},\boldsymbol{\Gamma}_p)$ under $H_0$, use the inequality $P(N(0, 1)\geq x)\leq e^{-x^2/2}$ for $x>0$ to see
\begin{align*}
P\big(\max_{1\leq i\leq p}\Big|\frac{\sqrt{n}\bar{Y}_i}{\sqrt{\gamma}_{ii}}\Big|\geq 2\sqrt{\log p}\big)\leq p \cdot P(|N(0, 1)|\geq 2\sqrt{\log p}) \leq \frac{2}{p}.
\end{align*}
Thus,
\begin{align}\label{onemaxcomsma}
(\max_{1\leq i\leq p}\Big|\frac{\sqrt{n}\bar{Y}_i}{\sqrt{\gamma}_{ii}}\Big|)^2=O_p(\log p).
\end{align}
Therefore, by \eqref{onemaxcompara0}, \eqref{gdz} and \eqref{onemaxcomsma}, we get
\begin{align*}
    \tilde{\tilde{T}}_{\max} = \tilde{T}_{\max} + o_p(1).
\end{align*}
\end{proof}

\subsection{Proof of Theorem \ref{two-sample-max}: max-type test for two-sample case}\label{pftwo-sample-max}
\begin{proof}
    
The proof shares same spirit as Theorem \ref{one-sample-max}. Let $\tilde{T}_{\max,2}=\max_{1\leq i \leq p}V_{i,2}^2$, by Theorem \ref{theorem_max}, we have $\tilde{T}_{\max,2}-2\log p +\log\log p$ converges to a Gumbel distribution with cdf $F(x)=\exp\{-\frac{1}{\sqrt{\pi}}e^{-x/2}\}$ as $p\to \infty$. 

Let $\tilde{\tilde{T}}_{\max,2}=\frac{n_1n_2}{n_1+n_2}\max_{1\leq i\leq p}\frac{(\bar{Y}_{i}^{(1)}-\bar{Y}_{i}^{(2)})^2}{\gamma_{ii}}$, then
\begin{align}\label{twomaxcompara0}
&|\tilde{\tilde{T}}_{\max,2}-\tilde{T}_{\max,2}|
= \Big|\frac{n_1n_2}{n_1+n_2}\max_{1\leq i\leq p}\frac{(\bar{Y}_{i}^{(1)}-\bar{Y}_{i}^{(2)})^2}{\gamma_{ii}}-\max_{1\leq i\leq p}V_{i,2}^2\Big| \nonumber\\
\leq & \Big|\frac{n_1n_2}{n_1+n_2}\max_{1\leq i\leq p}\frac{(\bar{Y}_{i}^{(1)}-\bar{Y}_{i}^{(2)})^2}{\gamma_{ii}}\Big| \max_{1\leq i\leq p}\Big| 1-\frac{\bg_i^T\bD_{\boldsymbol{\Sigma}_p}^{-1/2}(\bar{\bQ}^{(1)}-\bar{\bQ}^{(2)})(\bar{\bQ}^{(1)}-\bar{\bQ}^{(2)})^T\bD_{\boldsymbol{\Sigma}_p}^{-1/2}\bg_i}{\bg_i^T(\bar{\bQ}^{(1)}-\bar{\bQ}^{(2)})(\bar{\bQ}^{(1)}-\bar{\bQ}^{(2)})^T\bg_i/\gamma_{ii}}\Big|.
\end{align}
Since $\sqrt{\frac{n_1n_2}{n_1+n_2}}(\bar{\bY}^{(1)}-\bar{\bY}^{(2)})=\sqrt{\frac{n_1n_2}{n_1+n_2}}\bG(\bar{\bQ}^{(1)}-\bar{\bQ}^{(2)})=\sqrt{\frac{n_1n_2}{n_1+n_2}}\bG(\bar{\bQ}^{(1)}-\bar{\bQ}^{(2)}-c\bI_p)\sim N(\boldsymbol{0},\boldsymbol{\Gamma}_p)$, by the same proof of Theorem \ref{one-sample-max-sum}, we have
\begin{align}\label{twomaxcomsma}
(\max_{1\leq i\leq p}\Big|\frac{\sqrt{\frac{n_1n_2}{n_1+n_2}}(\bar{Y}_{i}^{(1)}-\bar{Y}_{i}^{(2)})}{\sqrt{\gamma}_{ii}}\Big|)^2=O_p(\log p),
\end{align}
and $\max_{1\leq i\leq p}\Big| 1-\frac{\bg_i^T\bD_{\boldsymbol{\Sigma}_p}^{-1/2}(\bar{\bQ}^{(1)}-\bar{\bQ}^{(2)})(\bar{\bQ}^{(1)}-\bar{\bQ}^{(2)})^T\bD_{\boldsymbol{\Sigma}_p}^{-1/2}\bg_i}{\bg_i^T(\bar{\bQ}^{(1)}-\bar{\bQ}^{(2)})(\bar{\bQ}^{(1)}-\bar{\bQ}^{(2)})^T\bg_i/\gamma_{ii}}\Big|=O_p(p^{-1/2-\alpha}(\log p)^{C_1/2})$ for $\alpha, C_1>0$. Thus, 
\begin{align}\label{twomaxtttdiff}
\tilde{\tilde{T}}_{\max,2}=\tilde{T}_{\max,2}+o_p(1).
\end{align}
Moreover, notice $\tilde{\tilde{T}}_{\max,2}$ converges to a Gumbel distribution which  implies $\tilde{\tilde{T}}_{\max,2}=O_p(\log p)$, and by the same argument in the proof of Theorem \ref{one-sample-max-sum}, we have $\max_{1\leq i\leq p} |\hat{\gamma}_{ii}^{-1}\gamma_{ii}-1|=O_p(\sqrt{n^{-1}\log p})$. Then as $p\to\infty$,
\begin{align}\label{twot}
|T_{\max,2}-\tilde{\tilde{T}}_{\max,2}|=& \left|\frac{n_1n_2}{n_1+n_2}\max_{1\leq i\leq p}\hat{\gamma}_{ii}^{-2}(\bar{Y}_{i}^{(1)}-\bar{Y}_{i}^{(2)})^2-\frac{n_1n_2}{n_1+n_2}\max_{1\leq i\leq p}\gamma_{ii}^{-1}(\bar{Y}_{i}^{(1)}-\bar{\bY}_{i}^{(2)})^2\right| \nonumber\\
\leq &|\tilde{\tilde{T}}_{\max}^{(2)}|\cdot
\max_{1\leq i\leq p} |\hat{\gamma}_{ii}^{-1}\gamma_{ii}-1| \nonumber\\
=& O\big(n^{-1/2}(\log p)^{3/2}\big)\to 0.
\end{align}
Therefore, by \eqref{twomaxtttdiff}-\eqref{twot}
\begin{align}\label{twomax3}
    T_{\max,2}=\tilde{T}_{\max,2}+o_p(1)
\end{align}
and then $T_{\max,2}$ converges to a Gumbel distribution with cdf $F(x)=\exp\{-\frac{1}{\sqrt{\pi}}e^{-x/2}\}$ as $p\to \infty$. 
\end{proof}

\section{Proof of Section 5: max-sum type test}
\subsection{Proof of Lemma \ref{theorem_max_sum}}\label{pftheorem_max_sum}

\begin{proof}
    
The central idea of the proof is akin to that of Lemma \ref{theorem_max}. Leveraging this concept, we establish the following bound,
\begin{align}\label{indebound}
    J_{lower}(Z,\epsilon_p)+o_p(1)\leq P\Big(\frac{1}{\upsilon_p}(S_p-\tr(\bG\boldsymbol{\Psi}_p\bG))\leq x,\ L_p>l_p\Big)\leq J_{upper}(Z,\epsilon_p)+o_p(1),\ as\  p\to \infty,
\end{align}
where 
$$\lim_{p\to\infty} J_{lower}(Z,\epsilon_p)=\lim_{p\to\infty} J_{upper}(Z,\epsilon_p)=\lim_{p\to\infty} P\Big(\frac{1}{\hat{\upsilon}_p}(S^Z_p-tr\boldsymbol{\Psi}_p)\leq x,\ L^Z_p>l_p\Big)=\Phi(x)\cdot [1-F(y)]$$ 
and 
\begin{align}\label{maxnote}
L_p=\max_{1\leq i \leq p}|H_i|\ \ \mbox{and}\ \ l_p= (2\log p -\log\log p+y)^{1/2}, \ \ \text{for any $y \in \mathbb{R}$},
\end{align}
here, the latter one  makes sense for sufficiently large $p$, and $\upsilon_p=[2\tr(\bG\boldsymbol{\Psi_p}\bG)^2]^{1/2}$, $S_p=H_1^2+\dots+H_p^2$, and $\hat{\upsilon}_p=[2\tr(\boldsymbol{\Psi}_p^2)]^{1/2}$, $S_p^Z=\sum_{i=1}^p|Z_i-\Expe\bar{Z}|^2$ and $L^Z_p=\max_{1\leq i \leq p}|Z_i-\Expe\bar{Z}|$. We partition the proof of this bound into the following two steps.

\textbf{Step 1:} By Theorem \ref{theorem_sum},
\begin{align}
P\Big(\frac{S_p-\mbox{tr}(\bG\boldsymbol{\Psi}_p\bG)}{\upsilon_p}\leq x\Big)=\Phi(x) \label{sumdistri}
\end{align}
as $p\to \infty$ for any $x\in \mathbb{R}$, where $\Phi(x)=\frac{1}{\sqrt{2\pi}}\int_{-\infty}^xe^{-t^2/2}\,dt$.
From Theorem \ref{theorem_max}, we have
\begin{align}
P\big(\max_{1\leq i \leq p}\big\{H_i^2\big\}-2\log p +\log\log p \leq y\big) \to
F(y)=\exp\Big\{-\frac{1}{\sqrt{\pi}}e^{-y/2}\Big\}\  \label{maxdistri1}
\end{align}
as $p\to \infty$ for any $y \in \mathbb{R}$.  To show asymptotic independence, it is enough to prove
\begin{align*}
\lim_{p\to \infty}P\Big(\frac{S_p-\mbox{tr}(\bG\boldsymbol{\Psi}_p\bG)}{\upsilon_p}\leq x,\ \max_{1\leq i \leq p}H_i^2-2\log p +\log\log p\leq y\Big)= \Phi(x)\cdot F(y)
\end{align*}
for any $x\in \mathbb{R}$ and $y \in \mathbb{R}$. 
Due to  \eqref{sumdistri}, the above condition we want to prove is equivalent to
\begin{align}\label{maxsumdistri}
\lim_{p\to \infty}P\Big(\frac{S_p-\mbox{tr}(\bG\boldsymbol{\Psi}_p\bG)}{\upsilon_p}\leq x,\ L_p>l_p\Big)= \Phi(x)\cdot [1-F(y)],
\end{align}
for any $x\in \mathbb{R}$ and $y \in \mathbb{R}$. Let
\begin{align}
A_p=\Big\{\frac{S_p-\mbox{tr}(\bG\boldsymbol{\Psi}_p\bG)}{\upsilon_p}\leq x\Big\}\ \ \ \mbox{and}\ \ \ B_{i}=\big\{|H_i|>l_p\big\} \label{setnote0}
\end{align}
for $1\leq i\leq p$. We can then write
\begin{align}\label{independent1}
P\Big(\frac{1}{\upsilon_p}(S_p-\tr(\bG\boldsymbol{\Psi}_p\bG))\leq x,\ L_p>l_p\Big)=P\Big(\bigcup_{i=1}^pA_pB_{i}\Big).
\end{align}
Here the notation $A_pB_i$ stands for $A_p\cap B_i$. Thus, by Bonferroni inequality, we get
\begin{align}\label{independent2}
J^H_{lower}:=\sum_{t=1}^{2k}(-1)^{t-1}\sum_{1\leq i_1<\dots< i_t\leq p}P(A_pB_{i_1}\dots B_{i_t})\leq P\Big(\bigcup_{i=1}^pA_pB_{i}\Big)&\nonumber\\
\ \ \ \ \leq \sum_{t=1}^{2k+1}(-1)^{t-1}\sum_{1\leq i_1<\dots< i_t\leq p}P(A_pB_{i_1}\dots B_{i_t}):=J^H_{upper}.&
\end{align}

\textbf{Step 2:} We further derive the upper and lower bound of \eqref{independent2} by considering 
$$P(A_pB_{i_1}\dots B_{i_t})=P(A_pB_{i_1}\dots B_{i_t},|\bar{Z}-\Expe \bar{Z}|\leq\epsilon_p)+P(A_pB_{i_1}\dots B_{i_t},|\bar{Z}-\Expe \bar{Z}|\geq\epsilon_p),$$
which investigate $H_i=Z_i-\bar{Z}$ ($i=1,\ldots,p$) in two events $\{H_i=Z_i-\bar{Z}\}\cap\{|\bar{Z}-\Expe \bar{Z}|\leq \epsilon_p\}$ and $\{H_i=Z_i-\bar{Z}\}\cap\{|\bar{Z}-\Expe \bar{Z}|\geq \epsilon_p\}$. Note that,
\begin{align}\label{PAB1}
   &\ \ \ \ \ \  P(A_pB_{i_1}\dots B_{i_t})= P(\frac{H_1^2+\cdots + H_p^2-\tr(\bG\boldsymbol{\Psi}_p\bG)}{\upsilon_p}\leq x, |Z_{i_1}-\bar{Z}|>l_p,\dots,|Z_{i_t}-\bar{Z}|>l_p)\nonumber\\
    &=P(\frac{\sum_{i=1}^p(Z_i-\Expe \bar{Z})^2-\tr\boldsymbol{\Psi}_p}{\upsilon_p}\leq x+\frac{p(\bar{Z}-\Expe \bar{Z})^2-(\tr\boldsymbol{\Psi}_p-\tr(\bG\boldsymbol{\Psi}_p\bG))}{\upsilon_p}, \nonumber\\
    &\ \ \ \ \ \ \ \ \ \ \ \ \ \ \ \ \ \ \ \ \ \ \ \ \ \ \ \ \ \ \ \ \ \ \ \ |Z_{i_1}-\bar{Z}|>l_p,\dots,|Z_{i_t}-\bar{Z}|>l_p,|\bar{Z}-\Expe \bar{Z}|<\epsilon_p)\nonumber\\
    &\ \ \ \ +P(\frac{\sum_{i=1}^p(Z_i-\Expe \bar{Z})^2-\tr\boldsymbol{\Psi}_p}{\upsilon_p}\leq x+\frac{p(\bar{Z}-\Expe \bar{Z})^2-(\tr\boldsymbol{\Psi}_p-\tr(\bG\boldsymbol{\Psi}_p\bG))}{\upsilon_p}, \nonumber\\
    &\ \ \ \ \ \ \ \ \ \ \ \ \ \ \ \ \ \ \ \ \ \ \ \ \ \ \ \ \ \ \ \ \ \ \ \ |Z_{i_1}-\bar{Z}|>l_p,\dots,|Z_{i_t}-\bar{Z}|>l_p,|\bar{Z}-\Expe \bar{Z}|\geq\epsilon_p)\nonumber\\
    &\leq P(\frac{\sum_{i=1}^p(Z_i-\Expe \bar{Z})^2-\tr\boldsymbol{\Psi}_p}{\hat{\upsilon}_p}\leq \frac{\upsilon_p}{\hat{\upsilon}_p}(x+\frac{p\epsilon_p^2-(\tr\boldsymbol{\Psi}_p-\tr(\bG\boldsymbol{\Psi}_p\bG))}{\upsilon_p}),\nonumber\\
    &\ \ \ \ \ \ \ \ \ \ \ \ \ \ \ \ \ \ \ \ \ \ \ \ \ \ \ \ \ \ \ \ \ \ \ \ |Z_{i_1}-\Expe \bar{Z}|>l_p-\epsilon_p,\dots,|Z_{i_t}-\Expe \bar{Z}|>l_p-\epsilon_p)\nonumber\\
    &\ \ \ \ +P(|\bar{Z}-\Expe \bar{Z}|\geq\epsilon_p),
\end{align}
where
\begin{align}\label{limup}
\lim_{p\rightarrow\infty}\frac{\mbox{tr}(\bG\boldsymbol{\Psi}_p\bG)^2}{\mbox{tr}(\boldsymbol{\Psi}_p^2)}&=\lim_{p\rightarrow\infty}\frac{\sum_{i=1}^p\lambda_{i}^2(\bG\boldsymbol{\Psi}_p\bG)}{\sum_{i=1}^p\lambda_{i}^2(\boldsymbol{\Psi}_p)}=1,
\end{align}
and as $p\to \infty$ 
$$\frac{p\epsilon_p^2-(\tr\boldsymbol{\Psi}_p-tr(\bG\boldsymbol{\Psi}_p\bG))}{\upsilon_p}\leq \frac{p\epsilon_p^2}{\sqrt{2p}}\to 0$$
by taking $\epsilon_p=\Big((4k+4)\frac{\lambda_{\max}(\boldsymbol{\Psi}_p)}{\sqrt{p}(\log p)^{-1}}p^{-1/2}\Big)^{1/2}$ and utilizing the fact that $\tr(\bG\boldsymbol{\Psi}_p\bG)\leq \tr\boldsymbol{\Psi}_p$, which is derived from the eigenvalue interlacing theorem.
Here, the last equlity in \eqref{limup} is due to 
\begin{align*}
    \frac{\sum_{i=1}^p\lambda_{i}^2(\bG\boldsymbol{\Psi}_p\bG)}{\sum_{i=1}^p\lambda_{i}^2(\boldsymbol{\Psi}_p)}\leq 1
\end{align*}
and
\begin{align*}
    \frac{\sum_{i=1}^p\lambda_{i}^2(\bG\boldsymbol{\Psi}_p\bG)}{\sum_{i=1}^p\lambda_{i}^2(\boldsymbol{\Psi}_p)}&\geq 1- \frac{\lambda_{\max}^2(\boldsymbol{\Psi}_p)}{\sum_{i=1}^p\lambda_{i}^2(\boldsymbol{\Psi}_p)}
    \geq 1-\frac{\lambda_{\max}^2(\boldsymbol{\Psi}_p)}{p\lambda_{\min}^2(\boldsymbol{\Psi}_p)}\geq 1-p^{-1+2\tau}
\end{align*}
which are derived from the eigenvalue interlacing theorem,
coupled with the use of $\frac{\lambda_{\max}(\boldsymbol{\Psi}_p)}{\lambda_{\min}(\boldsymbol{\Psi}_p)}=O(p^{\tau})$ ($\tau\in(0,1/4)$) in Assumption \ref{assumption_max}.
Similarly, 
\begin{align}\label{PAB2}
     &\ \ \ \ \ \  P(A_pB_{i_1}\dots B_{i_t})\nonumber\\
    &\geq P(\frac{\sum_{i=1}^p(Z_i-\Expe \bar{Z})^2-\tr\boldsymbol{\Psi}_p}{\upsilon_p}\leq x+\frac{p(\bar{Z}-\Expe \bar{Z})^2-(\tr\boldsymbol{\Psi}_p-\tr(\bG\boldsymbol{\Psi}_p\bG))}{\upsilon_p},\nonumber\\
    &\ \ \ \ \ \ \ \ \ \ \ \ \ \ \ \ \ \ \ \ \ \ \ \ \ \ \ \ \ \ \ \ \ \ \ \ |Z_{i_1}-\bar{Z}|>l_p,\dots,|Z_{i_t}-\bar{Z}|>l_p,|\bar{Z}-\Expe \bar{Z}|<\epsilon_p)\nonumber\\
    &\ \ -P(\frac{\sum_{i=1}^p(Z_i-\Expe \bar{Z})^2-\tr\boldsymbol{\Psi}_p}{\hat{\upsilon}_p}\leq \frac{\upsilon_p}{\hat{\upsilon}_p}(x-\frac{(\tr\boldsymbol{\Psi}_p-\tr(\bG\boldsymbol{\Psi}_p\bG))}{\upsilon_p}),\nonumber\\
    &\ \ \ \ \ \ \ \ \ \ \ \ \ \ \ \ \ \ \ \ \ \ \ \ \ \ \ \ \ \ \ \ \ \ \ \ |Z_{i_1}-\Expe \bar{Z}|>l_p+\epsilon_p,\dots,|Z_{i_t}-\Expe \bar{Z}|>l_p+\epsilon_p,|\bar{Z}-\Expe \bar{Z}|\geq\epsilon_p)\nonumber\\
    &\geq P(\frac{\sum_{i=1}^p(Z_i-\Expe \bar{Z})^2-\tr\boldsymbol{\Psi}_p}{\hat{\upsilon}_p}\leq \frac{\upsilon_p}{\hat{\upsilon}_p}(x-\frac{(\tr\boldsymbol{\Psi}_p-\tr(\bG\boldsymbol{\Psi}_p\bG))}{\upsilon_p}), \nonumber\\
    &\ \ \ \ \ \ \ \ \ \ \ \ \ \ \ \ \ \ \ \ \ \ \ \ \ \ \ \ \ \ \ \ \ \ \ \ |Z_{i_1}-\Expe \bar{Z}|>l_p+\epsilon_p,\dots,|Z_{i_t}-\Expe \bar{Z}|>l_p+\epsilon_p,|\bar{Z}-\Expe \bar{Z}|<\epsilon_p)\nonumber\\
    &\ \ -P(|\bar{Z}-\Expe \bar{Z}|\geq\epsilon_p),
\end{align}
where
\begin{align}\label{smallterm0}
    \lim_{p\to \infty}\frac{(\tr\boldsymbol{\Psi}_p-\tr(\bG\boldsymbol{\Psi}_p\bG))}{\hat{\upsilon_p}}=0.
\end{align}
Here \eqref{smallterm0} is due to $\tr\boldsymbol{\Sigma}_p^2\geq p$ and 
\begin{align*}
       \lim_{p\to \infty}\frac{\tr(\bG\boldsymbol{\Psi}_p\bG)-\tr\boldsymbol{\Psi}_p}{\hat{\upsilon_p}} &= \frac{\tr(\bG\boldsymbol{\Psi}_p\bG)-\tr\boldsymbol{\Psi}_p}{\sqrt{2\tr\boldsymbol{\Psi}_p^2}}
=\frac{||\boldsymbol{1}_p'\boldsymbol{\Psi}_p\boldsymbol{1}_p||_2}{p\sqrt{2\tr\boldsymbol{\Psi}_p^2}}
\leq\frac{||\boldsymbol{\Psi}_p||_2}{\sqrt{p}}
<<(\log p)^{-1}\to 0,
\end{align*}
where the last one is derive from $\lambda_{\max}(\boldsymbol{\Psi}_p)<<\sqrt{p}(\log p)^{-1}$ in Assumption \ref{assumption_max}.
Thus, by \eqref{independent2}-\eqref{PAB2}, we get
\begin{align}\label{independent3}
&J_{lower}(Z,\epsilon_p) -\sum_{t=1}^{2k}(-1)^{t-1}\sum_{1\leq i_1<\dots< i_t\leq p}P(|\bar{Z}-\Expe \bar{Z}|\geq\epsilon_p)\nonumber\\
&\leq J^H_{lower}\leq P\Big(\bigcup_{i=1}^pA_pB_{i}\Big)\leq J^H _{upper}\nonumber\\
&\leq J_{upper}(Z,\epsilon_p) +\sum_{t=1}^{2k+1}(-1)^{t-1}\sum_{1\leq i_1<\dots< i_t\leq p}P(|\bar{Z}-\Expe \bar{Z}|\geq\epsilon_p).
\end{align}
where
\begin{align*}
    J_{lower}(Z,\epsilon_p)&=\sum_{t=1}^{2k}(-1)^{t-1}\sum_{1\leq i_1<\dots< i_t\leq p}P\Big(\frac{\sum_{i=1}^p(Z_i-\Expe \bar{Z})^2-\tr\boldsymbol{\Psi}_p}{\hat{\upsilon}_p}\leq \frac{\upsilon_p}{\hat{\upsilon}_p}(x+\frac{p\epsilon_p^2-(\tr\boldsymbol{\Psi}_p-\tr(\bG\boldsymbol{\Psi}_p\bG))}{\upsilon_p}), \nonumber\\
&\ \ \ \ \ \ \ \ \ \ \ \ \ \ \ \ |Z_{i_1}-\Expe \bar{Z}|>l_p+(-1)^{t-1}\epsilon_p,\dots,|Z_{i_t}-\Expe \bar{Z}|>l_p+(-1)^{t-1}\epsilon_p\Big)\nonumber
\end{align*}
and
\begin{align*}
    J_{upper}(Z,\epsilon_p)&=\sum_{t=1}^{2k+1}(-1)^{t-1}\sum_{1\leq i_1<\dots< i_t\leq p}P\Big(\frac{\sum_{i=1}^p(Z_i-\Expe \bar{Z})^2-\tr\boldsymbol{\Psi}_p}{\hat{\upsilon}_p}\leq \frac{\upsilon_p}{\hat{\upsilon}_p}(x-\frac{(\tr\boldsymbol{\Psi}_p-\tr(\bG\boldsymbol{\Psi}_p\bG))}{\upsilon_p}), \nonumber\\
&\ \ \ \ \ \ \ \ \ \ \ \ \ \ \ \ |Z_{i_1}-\Expe \bar{Z}|>l_p+(-1)^{t}\epsilon_p,\dots,|Z_{i_t}-\Expe \bar{Z}|>l_p+(-1)^{t}\epsilon_p\Big).\nonumber
\end{align*}
By \eqref{smallpro1}-\eqref{smallpro2} and \eqref{independent3}, we get \eqref{indebound}. By Theorem 3 in \cite{Feng2022AsymptoticIO}, we have 
\begin{align}\label{zlimit}
    \liminf_{p\to\infty}J_{lower}(Z,\epsilon_p)&=\limsup_{p\to\infty}J_{upper}(Z,\epsilon_p)\nonumber\\
    &=\lim_{p\to\infty} P\Big(\frac{1}{\hat{\upsilon}_p}(S^Z_p-\tr\boldsymbol{\Psi}_p)\leq x,\ L^Z_p>l_p\Big)=\Phi(x)\cdot [1-F(y)].
\end{align}
Hence, by \eqref{indebound} and \eqref{zlimit}, we get
\begin{align*}
\lim_{p\to \infty}P\Big(\frac{S_p-\tr(\bG\boldsymbol{\Psi}_p\bG)}{\upsilon_p}\leq x,\ \max_{1\leq i \leq p}H_i^2-2\log p +\log\log p\leq y\Big)= \Phi(x)\cdot F(y)
\end{align*}
\end{proof}

\subsection{Proof of Theorem \ref{one-sample-max-sum}: max-sum type test for one-sample case}\label{pfone-sample-max-sum}
\begin{proof}

   \textbf{Step 1:} By \eqref{onesum3} and \eqref{onemax3}, we have
    \begin{align*}
    T_{\tsum}=\tilde{T}_{\tsum}+o_p(1) \text{\ and \ }
        T_{\max}=\tilde{T}_{\max}+o_p(1).
    \end{align*}\\
    \textbf{Step 2:} By Lemma \ref{theorem_max_sum}, we have $\tilde{T}_{\tsum}$ and $\tilde{T}_{\max}$ are asymptotically independent.\\ 
    \textbf{Step 3:} By Steps 1-2 and Lemma S.10 in \cite{Feng2022AsymptoticIO}, $T_{\tsum}$ and $T_{\max}$ are asymptotically independent.
\end{proof} 

\subsection{Proof of Theorem \ref{two-sample-max-sum}: max-sum type test for two-sample case}\label{pftwo-sample-max-sum}
\begin{proof}
    
  \textbf{Step 1:} By \eqref{twosum3} and \eqref{twomax3}, we have
    \begin{align*}
    T_{\tsum,2}=\tilde{T}_{\tsum,2}+o_p(1) \text{\ and \ }
        T_{\max,2}=\tilde{T}_{\max,2}+o_p(1).
    \end{align*}\\
\textbf{Step 2:} By Lemma \ref{theorem_max_sum}, we have $\tilde{T}_{\tsum,2}$ and $\tilde{T}_{\max,2}$ are asymptotically independent.\\
    \textbf{Step 3:} By Steps 1-2 and Lemma S.10 in \cite{Feng2022AsymptoticIO}, $T_{\tsum,2}$ and $T_{\max,2}$ are asymptotically independent.

\end{proof}
\end{document}